\documentclass[11pt]{sat}

\newcommand{\sect}[1]{\section{#1}\setcounter{equation}{0}}

\title{A Survey of Weighted Polynomial Approximation
with Exponential Weights}
\def\shorttitle{Weighted Approximation}

\author{D.~S.~Lubinsky}
\def\shortauthor{D.~S.~Lubinsky}

\def\versiondate{1 January 2007}

\def\abstracttext{Let $W:\mathbb{R}\rightarrow (0,1]$ be continuous. Bernstein's approximation
problem, posed in 1924, deals with approximation by polynomials in
the weighted uniform norm $f\rightarrow
\norm{fW}_{L_{\infty}(\mathbb{R})}$. The qualitative form of this
problem was solved by Achieser, Mergelyan, and Pollard, in the
1950's. Quantitative forms of the problem were actively
investigated starting from the 1960's. We survey old and recent
aspects of this topic, including the Bernstein problem, weighted
Jackson and Bernstein Theorems, Markov--Bernstein and Nikolskii
inequalities, orthogonal expansions and Lagrange interpolation. We
present the main ideas used in many of the proofs, and different
techniques of proof, though not the full proofs. The class of
weights we consider is typically even, and supported on the whole
real line, so we exclude Laguerre type weights on $[0,\infty )$.
Nor do we discuss Saff's weighted approximation problem, nor the
asymptotics of orthogonal polynomials. }

\def\MSCnumbers{41A10, 41A17, 41A25, 42C10}

\def\keywords{Weighted approximation, polynomial approximation,
Jackson--Bernstein theorems}

\usepackage{amssymb}

\def\func#1{\mathop{\rm #1}}%
\def\norm#1{\Vert#1\Vert}
\def\dword#1{{\bf #1}}
\def\floor#1{\lfloor #1\rfloor}
\def\textsl#1{{\sl #1}}

\newtheorem{thh}{Theorem}[section]
\newtheorem{cor}[thh]{Corollary}
\newtheorem{lem}[thh]{Lemma}
\newtheorem{deff}[thh]{Definition}

\def\endproofsymbol{\makeblanksquare6{.4}}
\def\endproof{\endproofsymbol\nopf}

\def\meop{~~\endproofsymbol}

\def\nopf{\medskip\goodbreak}

\def\makeblanksquare#1#2{
\dimen0=#1pt\advance\dimen0 by -#2pt
      \vrule height#1pt width#2pt depth0pt\kern-#2pt
      \vrule height#1pt width#1pt depth-\dimen0 \kern-#1pt
      \vrule height#2pt width#1pt depth0pt \kern-#2pt
      \vrule height#1pt width#2pt depth0pt
}

\def\newline{\vskip15pt plus 2pt minus 2pt}
\def\myproof#1#2 {\par
                      \ifdim\lastskip<15pt
                      \removelastskip\penalty-200
                      \vskip15pt plus3pt minus3pt
                      \fi
                       {\def\a{#1}
                       \ifx\a\empty
                       {\noindent\bf #2.}
                       \else
                       {\noindent\bf #2 of #1.}
                       \fi}\enspace}

\def\proof#1.{\myproof{#1}{Proof}}
\def\proofi#1.{\myproof{#1}{Main idea of the proof}}

\def\startpagenumber{1}
\def\volumenumber{3}
\def\year{2007}

\newcommand{\beginddoc}{
\begin{document}
\maketitle
\begin{abstract}
\abstracttext \vskip1pt MSC: \MSCnumbers
\ifx\keywords\empty\else\vskip1pt Keywords: \keywords\fi
\end{abstract}
\insert\footins{\scriptsize
\medskip
\baselineskip 8pt \leftline{Surveys in Approximation Theory}
\leftline{Volume \volumenumber, \year.
pp.~\thepage--\pageref{endpage}.} \leftline{Copyright \copyright\
2006 Surveys in Approximation Theory.} \leftline{ISSN 1555-578X}
\leftline{All rights of reproduction in any form reserved.}
\smallskip
\par\allowbreak}
\tableofcontents}
\renewcommand\rightmark{\ifodd\thepage{\it \hfill\shorttitle\hfill}\else {\it \hfill\shortauthor\hfill}\fi}
\markboth{{\it \shortauthor}}{{\it \shorttitle}} \markright{{\it
\shorttitle}}
\def\endddoc{\label{endpage}\end{document}}
\date{{\small \versiondate}}
\setlength\oddsidemargin{0pc} \setlength\evensidemargin{0pc}
\setlength\topmargin{0in} \setlength\textwidth{6.5in}
\setlength\textheight{8.6in}
\beginddoc



\sect{Bernstein's Approximation Problem}

The first quarter of the twentieth century was a great period for
approximation theory. In that time, Dunham Jackson and Sergei
Bernstein completed their great works on the degree of
approximation. M\"{u}ntz proved his theorem on approximation by
powers $\left\{x^{\lambda _{j}}\right\} _{j=0}^{\infty }$, solving
a problem of Bernstein, Faber introduced Faber polynomials and
Faber series, and Szeg\H{o} was developing the theory of
orthogonal polynomials on the unit circle. Right at the end of
that quarter, in 1924, Bernstein \cite{Bernstein1924} came up with
a problem that became known as \textit{Bernstein's approximation
problem}, and whose ramifications continue to be explored to this
day.

One can speculate that one day Bernstein must have felt confined
in approximating on bounded intervals, and so asked: by
Weierstrass, we know that we can uniformly approximate any
continuous function on a compact interval by polynomials. Are
there analogues on the whole real line? The first thing to deal
with is the unboundedness of polynomials on unbounded intervals.
Clearly we need to damp the growth of a polynomial at infinity, by
multiplying by a weight. For example, consider
\[
P(x)\exp(-x^{2}) ,\qquad x\in
\mathbb{R},
\] %
where $P$ is a polynomial, or more generally,
\[
P(x)W(x).
\] %
Here $W$ must decay sufficiently fast at $\pm \infty $ to
counteract the growth of every polynomial. That is,
\begin{equation}
\lim_{\vert x\vert \rightarrow \infty }x^{n}W(x)=0%
\mbox{, }\qquad n=0,1,2,\ldots\mbox{ .}
\end{equation}

What can be approximated, and in what sense? This problem is known
as \dword{Bernstein's approximation problem}. A more precise
statement is as follows: let $W:\mathbb{R}\rightarrow \lbrack
0,1]$ be measurable. When is it true that for every continuous
$f:\mathbb{R}\rightarrow \mathbb{R}$ with
\[
\lim_{\vert x\vert \rightarrow \infty }(fW)
(x)=0,
\] %
there exists a sequence of polynomials $\left\{ P_{n}\right\}
_{n=1}^{\infty }$ with
\[
\lim_{n\rightarrow \infty }\norm{(f-P_{n}) W} _{L_{\infty }(\mathbb{R}) }=0?
\] %
If true we then say that the polynomials are \dword{dense}, or
that Bernstein's
problem has a positive solution. The restriction that $fW$ has limit $0$ at $%
\pm \infty $ is essential: if $x^{k}W(x)$ is bounded
on the real line for every non-negative $k$, then $x^{k}W(x) $ has limit $0$ at $\pm \infty $ for every such $k$, and
so the same is true of every weighted polynomial $PW$. So we could
not hope to approximate, in the
uniform norm, any function $f$ for which $fW$ does not have limit $0$ at $%
\pm \infty $.

When $W$ vanishes on a set of positive measure, or is not
continuous, extra complications ensue. So in the sequel, we shall
assume that $W$ is both \textit{positive and continuous}. The more
general case is surveyed at length in \cite{Koosis1988}. The case
where we approximate only on a countable set of points is included
in that study.

In his early works on the problem, Bernstein often assumed that
$1/W$ is the restriction to the real line of an even entire
function with positive (even order) Maclaurin series coefficients.
Other major contributors were N.~I.~Akhiezer, K.~I.~Babenko, L.~de
Branges, L.~Carleson, M.~M.~Dzrbasjan, T.~Hall, S.~Izumi,
T.~Kawata, S.~N.~Mergelyan, and V.~S.~Videnskii. Yes, that's
Lennart Carleson, and the same Mergelyan of Mergelyan Theorem fame
(on uniform approximation by polynomials on compact subsets of the
plane).

Bernstein's approximation problem was solved independently by
Achieser, Mergelyan, and Pollard, in the 1950's. Mergelyan
\cite[p.~147]{Koosis1988} \
introduced a regularization of the weight%
\[
\Omega (z) =\sup \left\{ \vert P(z) \vert :P\mbox{ a polynomial and }\sup_{t\in
\mathbb{R}}\frac{\vert P(t) W(t)
\vert }{\sqrt{1+t^{2}}}\leq 1\right\} .
\] %

\begin{thh}  {\bf (Mergelyan)} \textsl{Let
}$W:\mathbb{R}\rightarrow (0,1]$\textsl{\ be continuous and
satisfy $(1.1)$. There is a positive answer to Bernstein's problem
iff}
\[
\int_{-\infty }^{\infty }\frac{\log \Omega (t) }{1+t^{2}}%
dt=\infty .
\]\end{thh}

In another formulation, there is a positive answer iff
\[
\Omega (z) =\infty
\] %
for at least one non-real $z$ (and then $\Omega (z)
=\infty $ for all non-real $z$). We shall outline the proof of
this in the next section.

Akhiezer (and perhaps Bernstein?) \cite[p.~158]{Koosis1988} used
instead the
regularization%
\[
W_{\ast }(z) =\sup \left\{ \vert P(z)
\vert :P\mbox{ a polynomial with }\norm{PW}_{L_{\infty
}( \mathbb{R}) }\leq 1\right\} .
\] %

\begin{thh} {\bf (Akhiezer)} \textsl{Let
}$W:\mathbb{R}\rightarrow (0,1]$\textsl{\ be continuous and
satisfy $(1.1)$. There is a positive answer to Bernstein's problem
iff}
\[
\int_{-\infty }^{\infty }\frac{\log W_{\ast }(t) }{1+t^{2}}%
dt=\infty .
\]\end{thh}

Finally, Pollard \cite{Pollard1953}, \cite[p.~164]{Koosis1988} showed:%

\begin{thh} {\bf  (Pollard)} \textsl{Let
}$W:\mathbb{R}\rightarrow (0,1]$\textsl{\ be continuous and
satisfy $(1.1)$. There is a positive answer to Bernstein's problem iff both}%
\[
\int_{-\infty }^{\infty }\frac{\log ( 1/W(t) ) }{%
1+t^{2}}dt=\infty
\] %
\textsl{and there exists a sequence of polynomials }$\left\{ P_{n}\right\} $%
\textsl{\ such that for each }$x,$
\[
\lim_{n\rightarrow \infty }P_{n}(x)W(x)
=1,
\] %
\textsl{while}
\[
\sup_{n\geq 1}\norm{P_{n}W} _{L_{\infty }( \mathbb{R}) }<\infty .
\]\end{thh}

Pollard later \cite{Pollard1955} reformulated this as
\[
\sup \left\{ \int_{-\infty }^{\infty }\frac{\log
\vert P(x)\vert }{1+x^{2}}dx:P\mbox{ a polynomial
with }\norm{PW} _{L_{\infty }( \mathbb{R}) }\leq
1\right\} =\infty .
\]

Of course, these are not very transparent criteria. When the
weight is in some sense regular, simplifications are possible.

\begin{thh}  \textsl{Let }$W$\textsl{\ be even, and
}$\log \left( 1/W( e^{x}) \right) $\textsl{\ be convex.
There is a positive answer to Bernstein's problem iff}
\begin{equation}
\int_{0}^{\infty }\frac{\log \left( 1/W(x)\right) }{1+x^{2}}%
dx=\infty .
\end{equation}\end{thh}

This result was proved by Lennart Carleson in 1951
\cite{Carleson1951}, although Misha Sodin pointed out to the
author that it appeared in a 1937 paper of Izumi and Kawata
\cite{IzumiKawata1937}. On perusing the latter paper, I agree that
it is clearly implicit in the results there, though not explicitly
stated. I have also seen the result attributed to M. Dzrbasjan.
The reader should conclude that the ``history'' presented in this
survey is by no means authoritative, but merely what could be
deduced (often from secondary sources) in the available time.

\begin{cor}  \textsl{Let }$\alpha >0$\textsl{\ and}
\begin{equation}
W_{\alpha }(x)=\exp ( -\vert x\vert ^{\alpha }).
\end{equation}%
\textsl{There is a positive answer to Bernstein's problem iff}
$\alpha \geq 1 $.\end{cor}

As regards necessary conditions, Hall showed that (1.2) is
necessary for density. When density fails, only a limited class of
entire functions can be approximated \cite{Krooetal1997}. A
comprehensive treatment of this topic is given in Koosis' book
\cite{Koosis1988}. A concise elegant exposition appears in
\cite[p.~28 ff.]{Lorentzetal1996}. A more abstract solution to
Bernstein's problem was given by Louis de Branges (of Bieberbach
fame) in 1959 \cite{deBranges1959}. The ideas in that paper have
had several recent ramifications \cite{Sodin1996B},
\cite{SodinYuditskii1997}.

What about extensions to $L_{p}$? At least when $W$ is continuous,
the answer is the same:

\begin{thh}  \textsl{Let }$W:\mathbb{R}\rightarrow
(0,1]$\textsl{\ be continuous and satisfy $(1.1)$. Let }$p\geq
1$\textsl{. Then any of the conditions of Achieser, Mergelyan, or
Pollard is necessary and sufficient so that for every measurable
}$f:\mathbb{R}\rightarrow \mathbb{R}$\textsl{\ with }
$\norm{fW} _{L_{p}( \mathbb{R}) }<\infty $\textsl{,
there exist polynomials }$\left\{ P_{n}\right\} $\textsl{\ with}
\[
\lim_{n\rightarrow \infty }\left\Vert (f-P_{n})W\right\Vert
_{L_{p}( \mathbb{R}) }=0.
\] \end{thh}

There is even a generalization of this theorem to
weighted M\"{u}ntz polynomials \cite{TieQiang2005}. As regards the
proof of the $L_{p}$ case, let us quote Koosis
\cite[p.~211]{Koosis1988}: \textquotedblleft In general, in the
\textit{kind} of approximation problem considered here, (that of
the \textit{density} of a certain simple class of functions in the
whole space), \textit{it makes very little difference which
}$L_{p}$\textit{\ norm is chosen}. If the proofs vary in
difficulty, they are hardest for the $L_{1}$
norm or the uniform norm.\textquotedblright\ Indeed, in his 1955 paper \cite%
{Pollard1955}, Pollard added about two pages to deal with the
$L_{p}$ case.

In the special case of $L_{2}$, there are connections to the \textit{%
classical moment problem}. Let $\left\{ s_{k}\right\} $ be a
sequence of real numbers. If there is a positive measure $\sigma $
on the real line such that
\begin{equation}
s_{k}=\int_{-\infty }^{\infty }t^{k}d\sigma (t) \mbox{, }%
\qquad k=0,1,2,\ldots\mbox{ ,}
\end{equation}%
then we say that the Hamburger Moment Problem for $\left\{
s_{k}\right\} $ has a solution. If there is only one solution
$\sigma $, of the equations
(1.4), we say the moment problem is \dword{determinate}. Existence of $%
\sigma $ is equivalent to positivity of certain determinants. The
determinacy of the moment problem is often a deeper and more
difficult
issue. If we define%
\[
\Sigma (z) =\sup \left\{ \vert P(z)
\vert ^{2}:P\mbox{ a polynomial with }\int_{-\infty
}^{\infty }\vert P\vert ^{2}d\sigma \leq 1\right\} ,
\] %
then a necessary and sufficient condition for determinacy is \cite[p.~141]%
{Koosis1988}
\[
\int_{-\infty }^{\infty }\frac{\log ^{+}\Sigma (t) }{1+t^{2}}%
dt=\infty .
\] %
Here $\log ^{+}x:=\max \left\{ 0,\log x\right\} $. In particular,
for the measure
\[
d\sigma (t) =\exp ( -\vert t\vert
^{\alpha }) dt,
\] %
one can use this to show that the corresponding moment problem is
determinate iff $\alpha \geq 1$. Recall that Bernstein's
approximation problem also had a positive solution iff $\alpha
\geq 1$.

There are still closer connections to Bernstein's problem
\cite[Thm.~4.3, p.~77]{Freud1971}. It is noteworthy that this was
obtained by M. Riesz in 1923 even before Bernstein's problem was
formulated.

\begin{thh} \textsl{Suppose that }$\sigma $\textsl{\ is a
positive measure on the real line, with finite moments }$\left\{
s_{k}\right\} $\textsl{, as in $(1.4)$.
Suppose, moreover, that }$\sigma $\textsl{\ has no point masses, so that }$%
x\mapsto \int_{-\infty }^{x}d\sigma $\textsl{\ is a continuous
function of }$x$\textsl{. Then the moment problem associated with
}$\left\{
s_{k}\right\} $\textsl{\ is determinate iff for every }$\sigma$\textsl{%
-measurable function }$f$\textsl{\ with }$\int f^{2}d\sigma
$\textsl{\
finite, and for every }$\varepsilon >0$\textsl{, there exists a polynomial }$%
P$\textsl{\ such that}
\[
\int_{-\infty }^{\infty }( f-P) ^{2}d\sigma
<\varepsilon .
\]\end{thh}

Without assuming the continuity of $\int d\sigma $, determinacy
still implies density of the polynomials \cite[p.~74
ff.]{Freud1971}. There is
also an important result on one-sided approximation \cite[Theorem 3.3, p.~73]%
{Freud1971}:

\begin{thh}  \textsl{Suppose that }$\sigma $\textsl{\ is a
positive measure on the real line, with finite moments }$\left\{
s_{k}\right\} $\textsl{, as in (1.4). Suppose moreover, that the
moment problem associated with }$\left\{
s_{k}\right\} $\textsl{\ is determinate. Let }$\varepsilon >0$\textsl{\ and }%
$f:\mathbb{R}\rightarrow \mathbb{R}$ \textsl{be a function that is
Riemann--Stieltjes integrable against }$d\sigma $\textsl{\ over
every finite interval, and improperly Riemann integrable over the
whole real line, and of
polynomial growth at }$\infty $\textsl{. Then there exist polynomials }$R$%
\textsl{\ and }$S$\textsl{\ such that }%
\[
R\leq f\leq S\quad \mbox{ in }\mathbb{R}
\] %
\textsl{and}%
\[
\int_{-\infty }^{\infty }( S-R) d\sigma <\varepsilon .
\]\end{thh}

Some recent work related to Bernstein's approximation problem and
its $L_{2}$
analogue appear in \cite{Borichev2001}, \cite{BorichevSodin1998}, \cite%
{BorichevSodin2001}, \cite{Levin1993}, \cite{PerezQuintana},
\cite{Pitt1983},
 \cite{PortillaQuintanaRodriguezTouris2004},
\cite{PortillaQuintanaRodriguezTouris}, \cite{Quintana},
\cite{Sodin1996}, \cite{Sodin1996B}, \cite{SodinYuditskii1997},
\cite{TieQiang2005}.

\sect{Some Ideas for the Resolution of Bernstein's Problem}

\bigskip In this section, we present some of the ideas used by Akhiezer,
Mergelyan, Pollard, to resolve Bernstein's Approximation problem.
It's relatively easy to derive the necessary parts of the
conditions. We follow \cite[p.~147 ff.]{Koosis1988} and
\cite[p.~28 ff.]{Lorentzetal1996}. Recall Mergelyan's
regularization of $W:$
\[
\Omega (z) =\sup \left\{ \vert P(z)
\vert :P\mbox{ a polynomial and }\sup_{t\in
\mathbb{R}}\frac{\vert P(t) W(t)
\vert }{\sqrt{1+t^{2}}}\leq 1\right\} .
\] %
In the sequel, $C, C_{1}, C_{2}, \ldots$ denote positive
constants independent of $n,x,t$ and polynomials $P$ of degree
$\leq n$. The same symbol does not necessarily denote the same
constant in different occurrences. We write $C=C( \alpha)$
or $C\neq C( \alpha ) $ to respectively show
that $C$ depends on $\alpha $, or does not depend on $\alpha $. We
use the notation $\sim $ in the following sense: given sequences
of real numbers $\left\{ c_{n}\right\} $ and $\left\{
d_{n}\right\} $, we write
\[
c_{n}\sim d_{n}
\] %
if for some positive constants $C_{1},C_{2}$ independent of $n$,
we have
\[
C_{1}\leq c_{n}/d_{n}\leq C_{2}.
\] %

\noindent \begin{lem}  \textsl{In order that Bernstein's
approximation problem has a positive
solution, it is necessary that for all }$z\in \mathbb{C}\backslash \mathbb{R}%
,$
\begin{equation}
\Omega (z) =\infty .
\end{equation}\end{lem}

\proof.
Suppose that Bernstein's problem has a positive solution. Fix $z\in \mathbb{C%
}\backslash \mathbb{R}$. Since the function $f(t)
=( t-z) ^{-1}$ is continuous, and $fW$ has limit $0$ at
$\infty $, we can find polynomials $\left\{ R_{n}\right\} $ such
that
\begin{equation}
\delta _{n}=\norm{ (f-R_{n})W} _{L_{\infty }( \mathbb{R})}
\rightarrow 0,
\end{equation}%
as $n\rightarrow \infty $. Let
\[
P_{n}(t) =\frac{1-( t-z) R_{n}(t) }{%
\delta _{n}}=\frac{t-z}{\delta _{n}}( f-R_{n}) ( t) .
\] %
Then for all $t\in \mathbb{R},$%
\[
\left\vert \frac{P_{n}(t) }{( t-z) }W( t) \right\vert =\frac{1}{\delta _{n}}
\vert f-R_{n}\vert (t) W(t) \leq 1 .
\] %
Next, for some constant $C$ independent of $t$, we have
\[
\left\vert \frac{t-z}{t-i}\right\vert \leq C,\qquad t\in
\mathbb{R}\mbox{.}
\] %
Then
\[
\sup_{t\in \mathbb{R}}\frac{\vert P_{n}(t) W(t)
\vert }{\sqrt{1+t^{2}}}=\sup_{t\in \mathbb{R}}\left\vert \frac{t-z}{t-i%
}\right\vert \left\vert \frac{P_{n}(t) }{( t-z) }%
W(t) \right\vert \leq C.
\] %
It follows that $P_{n}/C$ is one of the polynomials considered in
forming the sup in $\Omega (z) $, so
\[
\Omega (z) \geq \left\vert \frac{P_{n}(z) }{C}%
\right\vert =\frac{1}{\delta _{n}C}\rightarrow \infty ,
\] %
$n\rightarrow \infty $.\endproof

Let's think how we can reverse this to show the condition is also
sufficient. If $\Omega (z) =\infty $, then reversing
the above argument, we see that there are polynomials $\left\{
R_{n}\right\} $
satisfying (2.2), with $f(t) =( t-z) ^{-1}$. If $%
\Omega (z) =\infty $ for all $z$, we can approximate
linear combinations
\[
f(t) =\sum_{j=1}^{m}\frac{c_{j}}{t-z_{j}}
\] %
by polynomials in the norm $\norm{\cdot W}
_{L_{\infty }( \mathbb{R}) }$. For small $\varepsilon
$, we can then approximate
\[
\frac{1}{2\varepsilon }\left[ \frac{1}{( t-z) -\varepsilon }-%
\frac{1}{( t-z) +\varepsilon }\right] =\frac{1}{( t-z) ^{2}-\varepsilon ^{2}}
\] %
and hence also $1/( t-z) ^{2}$. Iterating this, we can
approximate $1/( t-z) ^{m}$ for any non-real $z$, and
$m\geq 1$. We can then approximate polynomials in $1/(t-z) $. The latter
can in turn uniformly approximate on the real line, any continuous function $%
f(t) $ that has limit $0$ at $\infty $. (Use
Weierstrass' Theorem and make a transformation \thinspace
$x=1/(t-c) $.) Since $W\leq 1$, we then obtain a
positive solution to Bernstein's problem, for continuous functions
that have limit $0$ at $\infty $. This class of continuous
functions is big enough to approximate arbitrary ones for the
Bernstein problem. For full details, see \cite[p.~148 ff.
]{Koosis1988}.

Next, we show:

\begin{lem} \textsl{In order that Bernstein's approximation
problem has a positive solution, it is necessary that}
\begin{equation}
\int_{-\infty }^{\infty }\frac{\log \left( 1/W\right) (t) }{%
1+t^{2}}dt=\infty .
\end{equation}\end{lem}

\proof. The basic tool is the inequality%
\begin{equation}
\log \left\vert P(i) \right\vert \leq \frac{1}{\pi }%
\int_{-\infty }^{\infty }\frac{\log \left\vert P(t)
\right\vert }{1+t^{2}}dt,
\end{equation}%
valid for all polynomials $P$. To prove this, suppose first that
$P$ has no zeros in the closed upper half-plane. We can then
choose an analytic branch
of $\log P$ there, and the residue theorem gives%
\[
\log P(i) =\frac{1}{\pi }\int_{-\infty }^{\infty
}\frac{\log P(t) }{1+t^{2}}dt.
\] %
Taking real parts, we obtain (2.4), with equality instead of
inequality. When $P$ has no zeros in the upper half-plane, but
possibly has zeros on the real axis, a continuity argument shows
that (2.4) persists, but with equality. Finally when $P(z) $ contains factors $(z-a) $ with $a$ in the
upper half-plane, we use
\begin{eqnarray*}
\log \left\vert a-i\right\vert \leq \log \left\vert \overline{a}%
-i\right\vert
=\frac{1}{\pi }\int_{-\infty }^{\infty }\frac{\log \left\vert \overline{a}%
-t\right\vert }{1+t^{2}}dt =\frac{1}{\pi }\int_{-\infty }^{\infty
}\frac{\log \left\vert a-t\right\vert }{1+t^{2}}dt
\end{eqnarray*}%
and divide such factors out from $P$. So we have (2.4). Now if
\begin{equation}
\sup_{t\in \mathbb{R}}\frac{\left\vert P(t) W(t) \right\vert }{\sqrt{1+t^{2}}}\leq 1,
\end{equation}%
we obtain%
\[
\frac{1}{\pi }\int_{-\infty }^{\infty }\frac{\log \left\vert
P(t) \right\vert }{1+t^{2}}dt\leq \frac{1}{\pi }\int_{-\infty }^{\infty }%
\frac{\log \left( 1/W(t) \right) }{1+t^{2}}dt+\frac{1}{2\pi }%
\int_{-\infty }^{\infty }\frac{\log \left( 1+t^{2}\right)
}{1+t^{2}}dt
\] %
and hence
\[
\log \left\vert P(i) \right\vert \leq \frac{1}{\pi }%
\int_{-\infty }^{\infty }\frac{\log \left( 1/W(t) \right) }{%
1+t^{2}}dt+\frac{1}{2\pi }\int_{-\infty }^{\infty }\frac{\log
\left( 1+t^{2}\right) }{1+t^{2}}dt.
\] %
Taking sup's over all $P$ satisfying (2.5) gives%
\[
\log \Omega (i) \leq \frac{1}{\pi }\int_{-\infty }^{\infty }%
\frac{\log \left( 1/W(t) \right) }{1+t^{2}}dt+\frac{1}{2\pi }%
\int_{-\infty }^{\infty }\frac{\log \left( 1+t^{2}\right)
}{1+t^{2}}dt.
\] %
As we assumed Bernstein's approximation problem has a positive
solution, Lemma 2.1 shows that $\Omega (i) =\infty $,
and hence (2.3) follows.  \endproof

With a little more work, this proof also gives%
\[
\int_{-\infty }^{\infty }\frac{\log \Omega (t) }{1+t^{2}}%
dt=\infty ,
\] %
see \cite[p.~153]{Koosis1988}.

\sect{Weighted Jackson and Bernstein Theorems}

In the 1950's the search began for quantitative estimates, rather
than theorems establishing the possibility of approximation.
Bernstein and Jackson had provided quantitative forms of
Weierstrass' Theorem in 1911 and 1912, two wonderful years for
approximation theory. Recall first the classical (unweighted)
case. Jackson and Bernstein independently proved that
\begin{equation}
E_{n}\left[ f\right] _{\infty }:=\inf_{\deg (P) \leq
n} \norm{f-P}_{L_{\infty }\left[ -1,1\right] }\leq \frac{C}{n}
\norm{f^{\prime}}_{L_{\infty }\left[ -1,1\right] },
\end{equation}%
with $C$ independent of $f$ and $n$, and the inf being over
(algebraic) polynomials of degree at most $n$. The rate is best
possible amongst absolutely continuous functions $f$ on $\left[
-1,1\right] $ whose derivative is bounded. More generally, if $f$
has a bounded $k$th derivative, then the rate is $O( n^{-k}) $. In addition, Jackson obtained general
results involving moduli of continuity: for example, if $f$ is
continuous, and its modulus of continuity is
\[
\omega \left( f;\delta \right) =\sup \left\{ \left\vert f(x)
-f(y) \right\vert :x,y\in \left[ -1,1\right] \mbox{ and }%
\left\vert x-y\right\vert \leq \delta \right\} ,
\] %
then
\[
E_{n}\left[ f\right] _{\infty }\leq C\omega \left( f;\frac{1}{n}\right) ,
\] %
where $C$ is independent of $f$ and $n$.

Bernstein also obtained remarkable converse theorems, which show
that the rate (or degree) of approximation is determined by the
smoothness of $f$. These are most elegantly stated for
trigonometric polynomial approximation. Let
\[
\mathcal{E}_{n}\left[ g\right] :=\inf_{\deg (R) \leq
n}\left\Vert g-R\right\Vert _{L_{\infty }\left[ 0,2\pi \right] }
\] %
denote the distance from a periodic function $g$ to the set of all
trigonometric polynomials $R$ of degree $\leq n$. Let $0<\alpha
<1$. Bernstein showed that
\[
\mathcal{E}_{n}\left[ g\right] =O( n^{-\alpha }) ,\ \
n\rightarrow
\infty \quad \Longleftrightarrow \quad \omega \left( g;t\right) =O( t^{\alpha }) ,%
\quad \ \ t\rightarrow 0+,
\] %
where $\omega \left( g;\cdot \right) $ is the modulus of
continuity of $g$ on $\left[ 0,2\pi \right] $, defined much as
above. That is, the error of approximation of a $2\pi$-periodic
function $g$ on $\left[ 0,2\pi \right] $
by trigonometric polynomials of degree at most $n$ decays with rate $%
n^{-\alpha }$ iff $g$ satisfies a Lipschitz condition of order
$\alpha $.
Moreover, if $k\geq 1,$%
\[
\mathcal{E}_{n}\left[ g\right] =O( n^{-k-\alpha }),\ \
n\rightarrow \infty \quad \Longleftrightarrow \quad\omega( g^{(k) };t) =O( t^{\alpha }) ,\quad
\mbox{ }t\rightarrow 0+.
\] %
Here in the converse implication, the existence and continuity of
the $k$th derivative $g^{(k) }$ is assured. Bernstein
never resolved the exact smoothness required for a rate of decay
of $n^{-1}$, or more generally $O( n^{-k}) $. The case
$k=1$ was solved much later in 1945 by A.~Zygmund, the father of
the Chicago school of harmonic analysis, and author of the classic
``Trigonometric Series'' \cite{Zygmund2002}. Zygmund used a second
order modulus of continuity.

For approximation by algebraic polynomials, converse theorems are
more complicated, as better approximation is possible near the
endpoints of the interval of approximation. Only in the 1980's
were complete characterizations obtained, with the aid of the
Ditzian--Totik modulus of continuity \cite{DitzianTotik1987}. An
earlier alternative approach is that of Brudnyi--Dzadyk--Timan
\cite{DeVoreLorentz1993}. We shall discuss only the Ditzian--Totik
approach, since that has been adopted in weighted polynomial
approximation. Define the symmetric differences%
\begin{eqnarray*}
\Delta _{h}f(x)&=&f( x+\frac{h}{2}) -f( x-\frac{%
h}{2}) ; \\
\Delta _{h}^{2}f(x)&=&\Delta _{h}\left( \Delta
_{h}f(x) \right) ; \\
&&\vdots \\
\Delta _{h}^{k}f(x)&=&\Delta _{h}\left( \Delta
_{h}^{k-1}f(x)\right)
\end{eqnarray*}%
so that
\begin{equation}
\Delta _{h}^{k}f(x)=\sum_{i=0}^{k}{k\choose i}(-1) ^{i}f( x+k\frac{h}{2}-ih) .
\end{equation}%
If any of the arguments of $f$ lies outside the interval of approximation ($%
\left[ -1,1\right] $ in this setting), we adopt the convention
that the
difference is $0$. The $r$th order Ditzian--Totik modulus of continuity in $%
L_{p}\left[ -1,1\right] $ is
\[
\omega _{\varphi }^{r}( f;h) _{p}=\sup_{0<h\leq t}
\norm{\Delta _{h\sqrt{1-x^{2}}}^{r}f(x)}_{L_{p}\left[ -1,1%
\right] }.
\] %
Note the factor
\[
\varphi (x)=\sqrt{1-x^{2}}
\] %
multiplying the increment $h$. This forces a smaller increment
near the endpoints $\pm 1$ of $\left[ -1,1\right] $, reflecting
the possibility of better approximation rates there.

For $1\leq p\leq \infty $, Ditzian and Totik \cite[Thm.~7.2.1, p.~79]%
{DitzianTotik1987} proved the estimate
\[
E_{n}\left[ f\right] _{p}:=\inf_{\deg (P) \leq n}
\norm{f-P}_{L_{p}\left[ -1,1\right] }\leq C\omega _{\varphi }^{r}( f;%
\frac{1}{n}) _{p},
\] %
with $C$ independent of $f$ and $n$. This implies the Jackson (or
Jackson--Favard) estimate \cite[p.~260]{DeVoreLorentz1993}%
\[
E_{n}\left[ f\right] _{p}\leq Cn^{-r} \norm{\varphi ^{r}f^{(r) }}_{L_{p}\left[ -1,1\right] },
\] %
$n\geq r$, provided $f^{(r-1)}$ is absolutely continuous, and the
norm on the right-hand side is finite. Moreover, they showed that
if $0<\alpha <r$, then \cite[p.~265]{DeVoreLorentz1993}
\begin{equation}
E_{n}\left[ f\right] _{p}=O( n^{-\alpha }) ,\qquad
n\rightarrow \infty ,
\end{equation}%
iff%
\[
\omega _{\varphi }^{r}\left( f;h\right) _{p}=O( h^{\alpha
}) ,\qquad h\rightarrow 0+.
\] %
For example, if (3.3) holds with $\alpha =3\frac{1}{2}$, this implies that $%
f $ has 3 continuous derivatives inside $\left( -1,1\right) $ and
$f^{\prime \prime \prime }$ satisfies a Lipschitz condition of
order 1/2 in each compact subinterval of $\left( -1,1\right) $.

This equivalence is easily deduced from the Jackson inequality
above and
the general converse inequality \cite[Theorem 7.2.4, p.~83]{DitzianTotik1987}%
\[
\omega _{\varphi }^{r}\left( f;t\right) _{p}\leq Mt^{r}\sum_{0<n<\frac{1}{t}%
}(n+1) ^{r-1}E_{n}\left[ f\right] _{p}, \qquad t\in
\left( 0,1\right) .
\] %
The constant $M$ depends on $r$, but is independent of $f$ and
$t$. Of course, this subject has a long and rich history, and all
we are attempting here is to set the background for developments
in weighted approximation. Please forgive the many themes omitted!

For weights on the whole real line, the first attempts at general
Jackson theorems seem due to Dzrbasjan \cite{Dzrbasjan1955}. He
obtained the correct weighted rates, but only when restricting the
approximated function to a finite interval. In the 1960's and
1970's, Freud and Nevai made major strides in this topic
\cite{Nevai1986}. That 1986 survey of Paul Nevai is still
relevant, and a very readable introduction to the subject.

Let us review some of the fundamental features discovered by
Freud, in the case of the weight $W_{\alpha }(x)=\exp
\left( -\vert x\vert ^{\alpha }\right) ,\alpha >1$. A
little calculus shows that the weighted monomial $x^{n}W_{\alpha
}(x)$ attains its maximum modulus on the real line at
\dword{Freud's number}
\[
q_{n}=\left( n/\alpha \right) ^{1/\alpha }.
\] %
Thereafter it decays quickly to zero. With this in mind, Freud and
Nevai proved that there are constants $C_{1}$ and $C_{2}$ such
that for all polynomials $P_{n}$ of degree at most $n$,
\begin{equation}
\norm{P_{n}W_{\alpha }}_{L_{p}\left( \mathbb{R}\right) }\leq
C_{2}\norm{P_{n}W_{\alpha }}_{L_{p}\left[ -C_{1}n^{1/\alpha
},C_{1}n^{1/\alpha }\right] }.
\end{equation}%
The constants $C_{1}$ and $C_{2}$ can be taken independent of
$n,P_{n}$ and
even the $L_{p}$ parameter $p\in \lbrack 1,\infty ]$. Outside the interval $%
\left[ -C_{1}n^{1/\alpha },C_{1}n^{1/\alpha }\right] $,
$P_{n}W_{\alpha }$ decays quickly to zero. This meant that one
cannot hope to approximate $fW$ by $P_{n}W$ outside $\left[
-C_{1}n^{1/\alpha },C_{1}n^{1/\alpha }\right] $. So either a
``tail term'' $\norm{fW_{\alpha }}_{L_{p}\left[ \left\vert
x\right\vert \geq C_{1}n^{1/\alpha }\right] }$ must appear in the
error estimate, or be handled some other way. Inequalities of the
form (3.4)
are called \dword{restricted range inequalities}, or \dword{%
infinite-finite range inequalities}.

The sharp form of these was found later by Mhaskar and Saff, using
potential theory \cite{Mhaskar1996}, \cite{MhaskarSaff1984},
\cite{MhaskarSaff1985},
\cite{MhaskarSaff1987}, \cite{SaffTotik1997}. Let $W=\exp (-Q) $%
, where $Q$ is even, and $xQ^{\prime }(x)$ is
positive and
increasing in $\left( 0,\infty \right) $, with limits $0$ and $\infty $ at $%
0 $ and $\infty $, respectively. For $n\geq 1$, let
$a_{n}=a_{n}(Q) $ denote the positive root of the equation%
\begin{equation}
n=\frac{2}{\pi }\int_{0}^{1}\frac{a_{n}tQ^{\prime }\left( a_{n}t\right) }{%
\sqrt{1-t^{2}}}dt.
\end{equation}%
We call $a_{n}$ the $n$th \dword{Mhaskar--Rakhmanov--Saff number}. For example, if $%
\alpha >0$, and $W_{\alpha }(x)=\exp( -\vert x\vert ^{\alpha }) ,$
\[
a_{n}=\left\{ 2^{\alpha -2}\frac{\Gamma \left( \alpha /2\right)
^{2}}{\Gamma \left( \alpha \right) }\right\} ^{1/\alpha
}n^{1/\alpha }.
\] %
Mhaskar and Saff established the \dword{Mhaskar--Saff identity}:
for polynomials $P_{n}$ of degree at most $n$,
\begin{equation}
\norm{P_{n}W}_{L_{\infty }\left( \mathbb{R}\right) }=
\norm{P_{n}W}_{L_{\infty }\left[ -a_{n},a_{n}\right] }.
\end{equation}%
Moreover, if $P_{n}$ is not the zero polynomial,
\[
\norm{P_{n}W}_{L_{\infty }\left( \mathbb{R}\backslash \left[
-a_{n},a_{n}\right] \right) }<\norm{P_{n}W}_{L_{\infty }\left[
-a_{n},a_{n}\right] }
\] %
and $a_{n}$ is asymptotically the ``smallest'' such number.

There are $L_{p}$ analogues, valid for all $p>0$. For example, if $%
\varepsilon >0$, there exists $C>0$ such that for $n\geq 1$ and polynomials $%
P_{n}$ of degree $\leq n,$%
\begin{equation}
\norm{P_{n}W}_{L_{p}\left( \mathbb{R}\backslash \left[
-a_{n}\left( 1+\varepsilon \right) ,a_{n}\left( 1+\varepsilon
\right) \right]
\right) }\leq e^{-Cn}\norm{P_{n}W}_{L_{p}\left[ -a_{n},a_{n}%
\right] }.
\end{equation}%
We shall discuss these more in Section 4.6.

The next task is to determine what happens on $\left[
-Ca_{n},Ca_{n}\right] $. Now if we had to approximate in the
unweighted setting on this interval, a
scale change in the Jackson--Bernstein estimate (3.1) gives%
\[
\inf_{\deg (P) \leq n}\norm{f-P}_{L_{\infty }\left[
-Ca_{n},Ca_{n}\right] }\leq CC_{1}\frac{a_{n}}{n} \norm{f^{\prime
}}_{L_{\infty }\left[ -Ca_{n},Ca_{n}\right] }.
\] %
Remarkably, the same is true when we insert the weight $W_{\alpha
}$ in both
norms:%
\begin{equation}
\inf_{\deg (P) \leq n}\norm{(f-P)
W_{\alpha
}}_{L_{\infty }\left[ -Ca_{n},Ca_{n}\right] }\leq C_{3}\frac{a_{n}}{%
n}\norm{f^{\prime }W_{\alpha }}_{L_{\infty }\left[
-Ca_{n},Ca_{n}\right] }.
\end{equation}%
Very roughly, this works for the following reason: it seems that
if $C$
is small enough, we can approximate $1/W_{\alpha }$ on $\left[ -Ca_{n},Ca_{n}%
\right] $ by a polynomial $R_{n/2}$ of degree $\leq n/2,$ and then
use the remaining part $n/2$ degree polynomial in $P$ to
approximate $fW_{\alpha }$ itself on $\left[ -Ca_{n},Ca_{n}\right]
$. In real terms, this approach works only for a small class of
weights. Nevertheless, it at least indicated the form that general
results should take. To obtain an estimate over the whole real
line, Freud then proved a ``tail inequality'', such as
\begin{equation}
\norm{fW_{\alpha }}_{L_{p}\left[ \vert x\vert \geq C a_{n}\right]
}\leq C_{4}\frac{a_{n}}{n}\norm{f^{\prime }W_{\alpha
}}_{L_{p}\left( \mathbb{R}\right) },
\end{equation}%
with $C_{4}$ independent of $f$ and $n$. Combining (3.8), (3.9),
and that
suitable weighted polynomials are tiny outside $\left[ -C_{1}a_{n},C_{1}a_{n}%
\right] $ yielded the following:

\begin{thh}
\textsl{Let }$1\leq p\leq \infty ,$\textsl{\ }$\alpha >1$\textsl{\ and }$f:%
\mathbb{R}\rightarrow \mathbb{R}$\textsl{\ be absolutely continuous, with }$%
\norm{f^{\prime }W_{\alpha }}_{L_{p}\left( \mathbb{R}\right)
}<\infty $\textsl{. Then }%
\begin{equation}
E_{n}\left[ f;W_{\alpha }\right] _{p}:=\inf_{\deg (P)
\leq
n}\norm{(f-P) W_{\alpha }}_{L_{p}\left( \mathbb{R}%
\right) }\leq C_{5}\frac{a_{n}}{n}\norm{f^{\prime }W_{\alpha }}
_{L_{p}\left( \mathbb{R}\right) },
\end{equation}%
\textsl{with }$C_{5}$\textsl{\ independent of }$f$\textsl{\ and}
$n$.\end{thh}

While this might illustrate some of the ideas, we emphasize that
the technical details underlying proper proofs of this Jackson (or
Jackson--Favard) inequality are formidable. Some of these ideas
will be illustrated in the next section. Freud and Nevai developed
an original theory of orthogonal polynomials for the weights
$W_{\alpha }^{2}$ partly for use in this approximation theory.

We note that Freud proved (3.10) for $W_{\alpha }$ for $\alpha \geq 2$ \cite%
{Freud1976}, \cite{Freud1977}. The technical estimates required to
extend this to the case $1<\alpha <2$ were provided by Eli Levin
and the author \cite{LevinLubinsky1987A}. What about $\alpha \leq
1$?\ Well, recall that the polynomials are only dense if $\alpha
\geq 1$, so there is no point in considering $\alpha <1$. But
$\alpha =1$ is still worth consideration, and we shall discuss
that below in Section 5.

One consequence of (3.10) is an estimate of the rate of weighted
polynomial approximation of $f$ in terms of that of $f^{\prime }$.
Indeed, if $P_{n}$ is any polynomial of degree $\leq n-1$, then
\[
E_{n}\left[ f;W_{\alpha }\right] _{p}=E_{n}\left[
f-P_{n};W_{\alpha }\right] _{p}\leq
C_{5}\frac{a_{n}}{n}\norm{\left( f-P_{n}\right) ^{\prime
}W_{\alpha }}_{L_{p}\left( \mathbb{R}\right) },
\] %
and since $P_{n}^{\prime }$ may be any polynomial of degree $\leq
n-1$, we
obtain the Favard or Jackson--Favard inequality%
\begin{equation}
E_{n}\left[ f;W_{\alpha }\right] _{p}\leq
C_{5}\frac{a_{n}}{n}E_{n-1}\left[ f^{\prime };W_{\alpha }\right]
_{p}.
\end{equation}%
This can be iterated:

\begin{cor}
\textsl{Let }$1\leq p\leq \infty ,$\textsl{\ }$r\geq 1$, $\alpha >1$\textsl{%
\ and }$f:\mathbb{R}\rightarrow \mathbb{R}$\textsl{\ have }$r-1$ \textsl{%
continuous derivatives. Assume, moreover, that }$f^{(r) }$\textsl{%
\ is absolutely continuous, with }$\norm{f^{(r)
}W_{\alpha }}_{L_{p}\left( \mathbb{R}\right) }<\infty $\textsl{.
Then for some
}$C$\textsl{\ independent of }$f$\textsl{\ and }$n,$%
\begin{equation}
E_{n}\left[ f;W_{\alpha }\right] _{p}\leq C_{5}\left( \frac{a_{n}}{n}\right)
^{r}\norm{f^{(r) }W_{\alpha }}_{L_{p}\left( \mathbb{%
R}\right) }.
\end{equation}\end{cor}

Freud also obtained estimates involving moduli of continuity. Here
one cannot avoid the tail term, and has to build it directly into
the modulus. Partly for this reason, there are many forms of the
modulus, and more than one way to decide which interval is the
principal interval, and over what interval we take the tail.
However it is done, it is awkward. We shall
follow essentially the modulus used by Ditzian and Totik \cite%
{DitzianTotik1987}, Ditzian and the author
\cite{DitzianLubinsky1997}, and Mhaskar \cite{Mhaskar1996}. All
these owe a great deal to earlier work of Freud.

The first order modulus for the weight $W_{\alpha }$ has the form%
\begin{eqnarray*}
\omega _{1,p}(f,W_{\alpha },t) =\sup_{0<h\leq t}\norm{W_{\alpha
}\left( \Delta _{h}f\right)}_{L_{p}[-h^{\frac{1}{1-\alpha }},h^{%
\frac{1}{1-\alpha }}]} +\inf_{c\in \mathbb{R}} \|(f-c)
W_{\alpha }\|
_{L_{p}\left( \mathbb{R}\backslash \lbrack -t^{\frac{1}{1-\alpha }},t^{\frac{%
1}{1-\alpha }}]\right) }.
\end{eqnarray*}%
Why the inf over the constant $c$ in the tail term? It ensures
that if $f$ is constant, then the modulus vanishes identically, as
one expects from a
first order modulus. Why the strange interval $[-h^{\frac{1}{1-\alpha }},h^{%
\frac{1}{1-\alpha }}]$? It ensures that when we substitute
\[
h=\frac{n^{1/\alpha }}{n}=n^{-1+1/\alpha },
\] %
then
\[
\lbrack -h^{\frac{1}{1-\alpha }},h^{\frac{1}{1-\alpha }}]=\left[
-n^{1/\alpha }, n^{1/\alpha }\right] =\left[ -C_{1}a_{n},C_{1}a_{n}%
\right] ,
\] %
for an appropriate constant $C_{1}$ (independent of $n$). More
generally
if $r\geq 1$, the $r$th order modulus is%
\begin{eqnarray}
\omega _{r,p}(f,W_{\alpha },t) &=&\sup_{0<h\leq t}\norm{W_{\alpha
}\left( \Delta _{h}^{r}f\right)}_{L_{p}[-h^{\frac{1}{1-\alpha }%
},h^{\frac{1}{1-\alpha }}]}  \nonumber \\[10pt]
&+&\inf_{\deg (P) \leq r-1} \|(f-P)
W_{\alpha }\|_{L_{p}\left( \mathbb{R}\backslash \lbrack -t^{\frac{1}{%
1-\alpha }},t^{\frac{1}{1-\alpha }}]\right) }.
\end{eqnarray}%
Again the inf in the tail term ensures that if $f$ is a polynomial
of degree $\leq r-1$, then the modulus of continuity vanishes
identically, as is
expected from an $r$th order modulus. The Jackson theorem takes the form%
\begin{equation}
E_{n}\left[ f;W_{\alpha }\right] _{p}\leq C\omega _{r,p}(f,W_{\alpha },\frac{%
a_{n}}{n}).
\end{equation}%
This is valid for $1\leq p\leq \infty $, and the constant $C$ is
independent of $f$ and $n$ (but depends on $p$ and $W_{\alpha }$).

One can consider more general weights than $W_{\alpha }$. Almost
invariably the weight considered has the form $W=\exp (-Q) $, and the rate of growth of $Q$ has a major impact on
the form of the modulus. Let us suppose, for example, that $Q$ is
of polynomial growth at $\infty $, the so-called \dword{Freud case}. The
most general class of such weights for which a Jackson theorem is
known is the following. It includes $W_{\alpha },\alpha
>1 $, but excludes $W_{1}$.

\begin{deff} {\bf  (Freud Weights)}
\textsl{Let }$W=\exp (-Q) $\textsl{, where }$Q:\mathbb{R}%
\rightarrow \mathbb{R}$\textsl{\ is even, }$Q^{\prime }$\textsl{\
exists and
is positive in }$\left( 0,\infty \right) $\textsl{. Moreover, assume that }$%
xQ^{\prime }(x)$\textsl{\ is strictly increasing,
with right
limit }$0$\textsl{\ at }$0$\textsl{, and for some }$\lambda ,A,B>1,C>0,$%
\begin{equation}
A\leq \frac{Q^{\prime }\left( \lambda x\right) }{Q^{\prime }(x)}%
\leq B,\qquad x\geq C.
\end{equation}%
\textsl{Then we write }$W\in \mathcal{F}$.\end{deff}

For such $W$, we take $a_{n}$ to be the positive root of the
equation (3.5) (the existence and uniqueness of $a_{n}$ is
guaranteed by the strict
monotonicity of $xQ^{\prime }(x)$). To replace the function $t^{%
\frac{1}{1-\alpha }}$, we can use the decreasing function of $t,$%
\begin{equation}
\sigma (t) :=\inf \left\{ a_{n}:\frac{a_{n}}{n}\leq
t\right\} ,\qquad t>0.
\end{equation}%
The modulus of continuity becomes%
\begin{eqnarray}
\omega _{r,p}(f,W,t) &=&\sup_{0<h\leq t}\norm{W\left( \Delta
_{h}^{r}f\right)}_{L_{p}[-\sigma (h) ,\sigma (h) ]}  \nonumber+\\[10pt]
&+&\inf_{\deg (P) \leq r-1} \|(f-P) W\|
_{L_{p}\left( \mathbb{R}\backslash \lbrack -\sigma (t)
,\sigma (t) ]\right) }.
\end{eqnarray}

The reader new to this subject will be encouraged to hear that
this strange looking creature has all the main properties of more
familiar moduli of
continuity \cite{DitzianLubinsky1997}, \cite{DitzianTikhonov2005}, \cite%
{DitzianTotik1987}, \cite{HoaKy1991}, \cite{KyHoa1994}, \cite{Mhaskar1996}:%

\begin{thh} {\bf (Properties of }$\omega _{r,p}${\bf )}
\textsl{Let }$W\in \mathcal{F}$, $r\geq 1,0<p\leq \infty $\textsl{. Let }$f:%
\mathbb{R}\rightarrow \mathbb{R}$, \textsl{and} \textsl{if }$p<\infty $%
\textsl{, assume that }$fW\in L_{p}\left( \mathbb{R}\right) $\textsl{. If }$%
p=\infty $\textsl{, assume in addition that }$f$\textsl{\ is
continuous and that }$fW$\textsl{\ has limit }$0$\textsl{\ at
}$\pm \infty $.
\begin{description}
\item[(a)]$\omega _{r,p}\left( f,W,t\right) $\textsl{\ is an increasing
function of }$t\geq 0.$
\item[(b)]
\[
\lim_{t\rightarrow 0+}\omega _{r,p}(f,W,t)=0.
\] %
\item[(c)] {Assume that }$p\geq 1$\textsl{, or assume that
}$W$\textsl{\
admits a Markov--Bernstein inequality }%
\begin{equation}
\Vert P^{\prime }W\Vert _{L_{p}( \mathbb{R}) }\leq C%
\frac{n}{a_{n}}\left\Vert PW\right\Vert _{L_{p}( \mathbb{R}) },
\end{equation}%
\textsl{valid for }$n\geq 1$\textsl{, and polynomials
}$P$\textsl{\ of degree }$\leq n$\textsl{, where }$C\neq C( n,P) $\textsl{. Then
there exists }$C_{1}\neq C_{1}( t,f) $\textsl{\ \ such that }%
\[
\omega _{r,p}( f,W,2t) \leq C_{1}\omega _{r,p}( f,W,t) ,\mbox{ }\qquad t\geq 0.
\] %
\item[(d)]{ Assume that }$p\geq 1$\textsl{. Then for }$t\geq 0,$%
\[
\omega _{r,p}( f,W,t) \leq Ct^{r}\Vert f^{(r) }W\Vert _{L_{p}( \mathbb{R}) },
\] %
\textsl{provided }$f^{(r-1) }$\textsl{\ is absolutely
continuous and the norm on the right-hand side is finite.}\newline
\item[(e)]\textsl{ Let }$r>1$\textsl{. There exists }$C\neq C( f,t) $%
\textsl{, such that}%
\begin{equation}
\omega _{r,p}( f,W,t) \leq C\omega _{r-1,p}( f,W,t) ,%
\mbox{ }\qquad t\geq 0.
\end{equation}%
\item[(f)] \textsl{ Assume the Markov--Bernstein inequality (3.18). Then
if }$q=\min
\left\{ 1,p\right\} $\textsl{, there is the Marchaud inequality\ }%
\[
\omega _{r,p}( f,W,t) \leq C_{1}t^{r}\left\{ \int_{t}^{C_{2}}%
\frac{\omega _{r+1}^{q}( f,W,u)
}{u^{rq+1}}du+\left\Vert fW\right\Vert _{L_{p}( \mathbb{R}) }^{q}\right\} ^{1/q},
\] %
\textsl{where }$C_{j}\neq C_{j}( f,t) ,j=1,2.$\newline
\item[(g)]\textsl{ Assume that }$\beta >1$\textsl{\ and }$0<p<q<\infty
$\textsl{.
Then there is the Ulyanov type inequality}%
\[
\omega _{r,q}( f,W_{\beta },t) \leq C\left\{
\int_{0}^{t}\left[ u^{\frac{1}{q}-\frac{1}{p}}\omega _{r,p}( f,W_{\beta },u) \right] ^{q}\frac{du}{u}\right\}
^{1/q},\mbox{ }\qquad t>0,
\] %
\textsl{where }$C\neq C( f,t) .$\newline
\item[(h)]\textsl{ Assume that }$p\geq 1$ \textsl{and for some }$0<\alpha <r,$%
\begin{equation}
\omega _{r,p}(f,W,t)=O( t^{\alpha }) ,\mbox{ }\qquad
t\rightarrow 0+.
\end{equation}%
\textsl{Let }$k=\floor{\alpha} $\textsl{, the integer part of }$%
\alpha $\textsl{. Then }$f^{(k) }$\textsl{\ exists
a.e.\ in the
real line, and }%
\begin{equation}
\omega _{r-k,p}(f^{(k)},W,t)=O( t^{\alpha -k}) ,\qquad
t\rightarrow 0+.
\end{equation}%
\textsl{If }$p=\infty $\textsl{, }$f^{(k) }$\textsl{\
is continuous on the real line.}
\end{description}\end{thh}

\proof.  (a) This is immediate as $\sigma (t) $ is a decreasing function of $t.$\newline (b) This
follows as for classical moduli of continuity. One first
establishes it for suitably restricted continuous functions, and
then approximates an arbitrary function by a continuous
one.\newline
(c) This is part of Theorem 1.4 in \cite[p.~104]{DitzianLubinsky1997}.%
\newline
(d) This is part of Corollary 1.8 in \cite[p.~105]{DitzianLubinsky1997}.%
\newline
(e) This follows easily from the definition (3.13) and the
recursive definition of the symmetric differences.\newline (f)
This is Corollary 1.7 in
\cite[p.~105]{DitzianLubinsky1997}.\newline
(g) This is part of Theorem 9.1 in \cite[p.~133]{DitzianTikhonov2005}.%
\newline
(h) See \cite[pp.~62--64]{DitzianTotik1987} for the analogous
proofs on a finite interval. \endproof

Note all the strictures for $p\leq 1$. Fundamentally these arise
because we cannot bound norms of functions in terms of their
derivatives when $p<1$. At least the Jackson theorem is the
obvious analogue of the result for finite intervals \cite[Theorem
1.2, p.~102]{DitzianLubinsky1997}:

\begin{thh}
\textsl{Let }$W\in \mathcal{F}$, $r\geq 1,$\textsl{\ and }$0<p\leq \infty $%
\textsl{. Let }$fW\in L_{p}\left( \mathbb{R}\right) $\textsl{. If
}$p=\infty $\textsl{, we also require }$f$\textsl{\ to be
continuous and }$fW$\textsl{\
to have limit }$0$\textsl{\ at }$\pm \infty $\textsl{. Then for }$n\geq r-1,$%
\begin{equation}
E_{n}\left[ f;W\right] _{p}\leq C_{1}\omega
_{r,p}(f,W,C_{2}\frac{a_{n}}{n}),
\end{equation}%
\textsl{where }$C_{1}$\textsl{\ and }$C_{2}$\textsl{\ are independent of }$f$%
\textsl{\ and }$n$\textsl{.}\end{thh}

In the case where $p\geq 1$, or $W$ admits the Markov--Bernstein
inequality (3.18), one can omit the constant $C_{2}$ inside the
modulus. This follows directly from Theorem 3.4(c). We shall say
much more about Markov--Bernstein inequalities in Section
7.

\begin{cor} \textsl{Let }$1\leq p\leq \infty
,$\textsl{\ }$r\geq 1$, $W\in \mathcal{F}$ \textsl{and
}$f:\mathbb{R}\rightarrow \mathbb{R}$\textsl{\ have }$r-1$
\textsl{continuous derivatives. Assume, moreover, that }$f^{(r) }$%
\textsl{\ is absolutely continuous, with }$\norm{f^{(r) }W}_{L_{p}\left( \mathbb{R}\right) }<\infty $\textsl{.
Then for
some }$C$\textsl{\ independent of }$f$\textsl{\ and }$n,$%
\begin{equation}
E_{n}\left[ f;W\right] _{p}\leq C_{5}\left( \frac{a_{n}}{n}\right)
^{r}\norm{f^{(r) }W}_{L_{p}\left( \mathbb{R}\right) }.
\end{equation}\end{cor}

The converse inequality, which can be interpreted as a Bernstein
type converse theorem, has the form
\cite[p.~105]{DitzianLubinsky1997}:

\begin{thh}
\textsl{Let }$0<p\leq \infty ,$\textsl{\ }$r\geq 1$, $W\in \mathcal{F}$.%
\textsl{\ Assume that }$W$\textsl{\ admits the Markov--Bernstein
inequality
$(3.18)$. Let }$q=\min \left\{ 1,p\right\} $\textsl{. For }$t\leq a_{1}$%
\textsl{, define the positive integer }$n=n(t) $\textsl{\ by }%
\begin{equation}
n:=n(t) :=\inf \left\{ k:\frac{a_{k}}{k}\leq t\right\}
.
\end{equation}%
\textsl{Then for some }$C\neq C\left( t,f\right) ,$ \textsl{\ }%
\begin{equation}
\omega _{r,p}(f,W,t)^{q}\leq C\left( \frac{a_{n}}{n}\right)
^{rq}\sum_{j=-1}^{\left[ \log _{2}n\right] }\left( \frac{2^{j}}{a_{2^{j}}}%
\right) ^{rq}E_{2^{j}}\left[ f;W\right] _{p}^{q},
\end{equation}%
\textsl{where we define }$E_{2^{-1}}:=E_{0}$\textsl{\ and }$\floor{\log _{2}n}%
$\textsl{\ denotes the largest integer }$\leq \log _{2}n$.\end{thh}

From this one readily deduces:

\begin{cor}
\textsl{Assume the hypotheses of the previous theorem, and let }$0<\alpha <r$%
\textsl{. Then}
\[
\omega _{r,p}(f,W,t)=O( t^{\alpha }) ,\mbox{ }\ \
t\rightarrow 0+\quad
\Longleftrightarrow \quad  E_{n}\left[ f;W\right] _{p}=O\left( \left( \frac{a_{n}}{n%
}\right) ^{\alpha }\right) ,\mbox{ }\ \ n\rightarrow \infty .
\] \end{cor}

A related smoothness theorem is given by Damelin
\cite{Damelin1998D}. \newline

One of the important tools in establishing this is $K$-functionals
and the concept of realization. This is a topic on its own. In the
setting of weighted polynomial approximation, it has been explored
by Freud and
Mhaskar, and later Ditzian and Totik, Damelin and the author. See \cite%
{Damelin1998}, \cite{DamelinLubinsky1998}, \cite{DitzianLubinsky1997}, \cite%
{FreudMhaskar1980}, \cite{FreudMhaskar1983}, \cite{Mhaskar1996}, \cite%
{Mhaskar2004}. In our context, an appropriate $K$-functional is%
\[
K_{r,p}( f,W,t^{r}) :=\inf_{g}\left\{ \Vert (f-g) W\Vert _{L_{p}( \mathbb{R})
}+t^{r}\Vert
g^{(r) }W\Vert _{L_{p}( \mathbb{R}) }\right\} ,%
\mbox{ }\qquad t\geq 0,
\] %
where the inf is taken over all $g$ whose $(r-1) $st
derivative is locally absolutely continuous. It works only for
$p\geq 1$, again because of the problems of estimating functions
in $L_{p}$, $p<1$, in terms of their derivatives. The
$K$-functional is equivalent to the modulus of continuity
in the following sense \cite[Thm.~1.4, Cor.~1.9, pp.~104--105]%
{DitzianLubinsky1997}:

\begin{thh}  \textsl{Let }$W\in \mathcal{F}$\textsl{\ and
}$p\geq 1$\textsl{. Then for
some }$C_{1},C_{2}>0$\textsl{\ independent of }$f,t,$%
\[
C_{1}\omega _{r,p}(f,W,t)\leq K_{r,p}\left( f,W,t^{r}\right) \leq
C_{2}\omega _{r,p}(f,W,t).
\]\end{thh}

The appearance of fairly general functions $g$ in the inf defining
$K_{r,p}$ helps to explain its usefulness. In fact, many of the
properties of the
modulus described above, go via the $K$-functional. When $p<1$, the $K$%
-functional is identically zero, so instead we use the \dword{realization
functional}, introduced by Hristov and Ivanov
\cite{HristovIvanov1990} and
analyzed by those authors and Ditzian \cite{Ditzianetal1995}:%
\[
\overline{K}_{r,p}( f,W,t^{r}) :=\inf_{P}\left\{
\norm{(f-P) W}_{L_{p}( \mathbb{R})
}+t^{r}\norm{P^{(r) }W}_{L_{p}( \mathbb{R}%
) }\right\} ,\quad \mbox{ }t>0.
\] %
Here the inf is taken over all polynomials $P$ of degree $\leq
n(t) $ and $n(t) $ is defined by (3.24). For $p\geq 1$, $%
K_{r,p}$ and $\overline{K}_{r,p}$ are equivalent
\cite{Ditzianetal1995}. For all $p$, one can prove \cite[Thm.~1.4,
p.~104]{DitzianLubinsky1997}:

\begin{thh}  \textsl{Assume }$W\in \mathcal{F}$
\textsl{and }$W$ \textsl{admits the Markov--Bernstein inequality
$(3.18)$. Then for some }$C_{1},C_{2}>0$\textsl{\
independent of }$f,t,$%
\[
C_{1}\omega _{r,p}(f,W,t)\leq \overline{K}_{r,p}\left( f,W,t^{r}\right) \leq C_{2}\omega _{r,p}(f,W,t).
\]\end{thh}

Observe that if we choose $t=a_{n}/n$, the above result actually
gives more than the Jackson estimate Theorem 3.5, at least when
$a_{k}/k$ decreases strictly with $k$. We obtain%
\[
\inf_{\deg (P) \leq n}\left\{ \left\Vert (f-P) W\right\Vert _{L_{p}\left( \mathbb{R}\right) }+\left(
\frac{a_{n}}{n}\right)
^{r}\norm{P^{(r) }W}_{L_{p}( \mathbb{R}%
) }\right\} \leq C_{2}\omega _{r,p}(f,W,\frac{a_{n}}{n}),
\] %
so there is an automatic bound on the $r$th derivatives of best
approximating polynomials.

The Freud weights above have the form $W=\exp (-Q) $,
where $Q$ is of polynomial growth at $\infty $, with $Q(x) $ growing at least as fast as $\vert x\vert
^{\alpha }$ for some $\alpha >1$. For the special weight $\exp
\left( -\vert x\vert \right) $, the polynomials are
still dense, but we have not established anything about the degree
of approximation. We shall devote a separate section to this. The
case where $Q$ is of faster than polynomial growth is often called the
\dword{Erd\H{o}s case}, and also has received some attention. The main
difference is that the modulus of continuity becomes more
complicated. Here is a suitable class of Erd\H{o}s
weights:

\begin{deff} {\bf (Erd\H{o}s Weights)}
\textsl{Let }$W=\exp (-Q) $\textsl{, where }$Q:\mathbb{R}%
\rightarrow \mathbb{R}$\textsl{\ is even, }$Q^{\prime }$\textsl{\
exists and is positive in }$\left( 0,\infty \right) $\textsl{.
Assume that }$xQ^{\prime
}(x)$\textsl{\ is strictly increasing, with right limit }$0$%
\textsl{\ at }$0$\textsl{, and the function}%
\begin{equation}
T(x):=\frac{xQ^{\prime }(x)}{Q(x) }
\end{equation}%
\textsl{is quasi-increasing in the sense that for some }$C>0,$%
\begin{equation}
0\leq x<y\;\Rightarrow\; T(x)\leq CT(y) ,
\end{equation}%
\textsl{while}
\[
\lim_{x\rightarrow \infty }T(x)=\infty .
\] %
\textsl{Assume, moreover, that for some }$C_{1},C_{2}$\textsl{\ and }$%
C_{3}>0, $%
\[
\frac{yQ^{\prime }(y) }{xQ^{\prime }(x)
}\leq C_{1}\left( \frac{Q(y) }{Q(x)
}\right) ^{C_{2}},\qquad y\geq x\geq C_{3}.
\] %
\textsl{Then we write }$W\in \mathcal{E}$.\end{deff}

Examples of weights in this class include
\[
W(x)=\exp( -\exp _{k}(\vert
x\vert ^{\alpha }) +\exp _{k}(0)) ,
\] %
where $\alpha >0$ and for $k\geq 1,$
\[
\exp _{k}=\exp( \exp(\cdots \exp( {})
))
\] %
denotes the $k$th iterated exponential. We set $\exp _{0}(x) =x$. For this weight,
\[
T(x)\sim \alpha x^{\alpha
}\prod\limits_{j=1}^{k-1}\exp _{j}\left( x^{\alpha }\right) ,
\] %
for large $x$ \cite[p.~9]{LevinLubinsky2001}, while the
Mhaskar--Rakhmanov--Saff number has the asymptotic \cite[p.~29]%
{LevinLubinsky2001}
\begin{eqnarray*}
a_{n} =\left\{ \log _{k-1}\left( \log
n-\frac{1}{2}\sum_{j=1}^{k+1}\log _{j}n+O(1) \right)
\right\} ^{1/\alpha } =\left( \log _{k}n\right) ^{1/\alpha }\left( 1+o(1) \right) ,
\end{eqnarray*}%
$n\rightarrow \infty $; $\log _{k}$ denotes the $k$th iterated
logarithm.

For general $W\in \mathcal{E}$, the modulus of continuity involves
the
function%
\[
\Phi _{t}(x):=\sqrt{1-\frac{\vert x\vert
}{\sigma (t) }}+T( \sigma (t))
^{-1/2},
\] %
where $\sigma (t) $ is as in (3.16). One may think of
this as an analogue of the function $\sqrt{1-x^{2}}$ which appears
in the Ditzian--Totik modulus of continuity, and it appears for
the same reason: in approximating by polynomials with Erd\H{o}s
weights, the rate of approximation improves towards the endpoints
of the \dword{Mhaskar--Rakhmanov--Saff interval}. The modulus is
\begin{eqnarray}
\omega _{r,p}(f,W,t) &=&\sup_{0<h\leq t}\norm{W(x)
\left( \Delta _{h\Phi _{t}(x)}^{r}f(x)
\right)}
_{L_{p}[-\sigma (2t) ,\sigma (2t) ]} \nonumber\\[10pt]
&+&\inf_{\deg (P) \leq r-1} \| (f-P) W\|
_{L_{p}\left( \mathbb{R}\backslash \lbrack -\sigma (4t) ,\sigma (4t) ]\right) }.
\end{eqnarray}%
Once one has this modulus, the Jackson estimate goes through
\cite[Thm.~1.2, p.~337]{DamelinLubinsky1998}:

\begin{thh}
\textsl{Let }$W\in \mathcal{E}$, $r\geq 1,$\textsl{\ and }$0<p\leq \infty $%
\textsl{. Let }$fW\in L_{p}\left( \mathbb{R}\right) $\textsl{. If
}$p=\infty $\textsl{, we also require }$f$\textsl{\ to be
continuous and }$fW$\textsl{\
to have limit }$0$\textsl{\ at }$\pm \infty $\textsl{. Then for }$n\geq 1,$%
\begin{equation}
E_{n}\left[ f;W\right] _{p}\leq C_{1}\omega
_{r,p}(f,W,C_{2}\frac{a_{n}}{n}),
\end{equation}%
\textsl{where }$C_{1}$\textsl{\ and }$C_{2}$\textsl{\ are independent of }$f$%
\textsl{\ and }$n$\textsl{.}\end{thh}

There are also converse estimates \cite{Damelin1998},
\cite{Damelin1999B}, \cite{DamelinLubinsky1998}. The details are
more difficult than for Freud weights, because of the more
complicated modulus. Analogous results for
exponential weights on $\left( -1,1\right) $ are given in \cite{Damelin1999}%
, \cite{Lubinsky1997A}, \cite{Lubinsky1997B}. For exponential
weights multiplied by a generalized Jacobi weight or other factor
having
singularities, see \cite{HorvathSzabados1998}, \cite{MastroianniSzabados2001}%
, \cite{MastroianniSzabados2002}, \cite{MastroianniSzabados2006}, \cite%
{Szabados2003}. For Laguerre and other exponential weights on
$\left( 0,\infty \right) $, see \cite{DeBonisetal2002}, \cite{KasugaSakai2004}, \cite%
{Mastroianni2002}. Geometric rates of approximation in the
weighted setting
have been explored in \cite{Krooetal1997}, \cite{Mhaskar1981}, \cite%
{Mhaskar1982}, \cite{Mhaskar1993}, \cite{Mhaskar1996}.

\sect{Methods for Proving Weighted Jackson Theorems}

In this section, we shall outline various methods to prove
weighted Jackson Theorems, but will not provide complete
expositions. In the unweighted case, on $\left[ -1,1\right] $,
many of the elegant methods involve convolution operators.
However, unfortunately these depend heavily on translation
invariance, so fail for the weighted case. We begin with two of
the oldest, used by Freud and Nevai.\newline

\subsection{Freud and Nevai's One-sided $L_{1}$ Method}

The $L_{1}$ method is primarily a tool to obtain estimates on the
degree of approximation for special functions such as
characteristic functions. Once one has it, one can use duality and
other tricks to go to $L_{\infty }$, and then interpolation for
$1<p<\infty $, and this will be done in the next section. The
method is based on the theory of orthogonal polynomials, and Gauss
quadratures. Under the tutelage of G\'{e}za Freud and Paul Nevai,
the two subjects of weighted approximation and orthogonal
polynomials for weights on the real line, developed in tandem
throughout the 1970's and 1980's. Mhaskar's monograph
\cite{Mhaskar1996} provides an excellent treatment of the material
in this and the next section.

Corresponding to the weight $W$, we define orthonormal polynomials%
\[
p_{n}(x)=p_{n}( W^{2},x) =\gamma
_{n}x^{n}+\cdots ,
\] %
where $\gamma _{n}>0$, satisfying
\[
\int_{-\infty }^{\infty }p_{n}p_{m}W^{2}=\delta _{mn.}
\] %
Note that the weight is $W^{2}$, not $W$. This convention
simplifies some formulations later on. Let us denote the zeros of
$p_{n}$ by
\[
x_{nn}<x_{n-1,n}< \cdots <x_{2n}<x_{1n}.
\] %
The $n$th Christoffel function for $W^{2}$ is%
\begin{equation}
\lambda _{n}( W^{2},x) =\inf_{\deg (P) \leq n-1}\frac{%
\int_{-\infty }^{\infty }(PW) ^{2}}{P^{2}(x) }.
\end{equation}%
It also satisfies%
\[
\lambda _{n}( W^{2},x)
=1\Bigl/\;\sum_{j=0}^{n-1}p_{j}^{2}(x).
\] %
Try this as an exercise: expand an arbitrary polynomial of degree
$\leq n-1$ in $\left\{ p_{j}\right\} _{j=0}^{n-1}$, and then use
Cauchy--Schwarz, and orthonormality.

Many readers unfamiliar with the detailed theory of orthogonal
polynomials, will nevertheless have seen the Christoffel functions
in the Gauss quadrature formula:\newline
\[
\int_{-\infty }^{\infty }PW^{2}=\sum_{j=1}^{n}\lambda _{n}( W^{2},x_{jn}) P( x_{jn}) ,
\] %
valid for all polynomials $P$ of degree $\leq 2n-1$. For fixed
$\xi $, one can also develop similar quadrature formulae, based on
the zeros of
\[
p_{n}(x)p_{n-1}( \xi ) -p_{n-1}(x) p_{n}( \xi ) .
\] %
These zeros are real and simple and include $\xi $. There are
important inequalities, the Posse--Markov--Stieltjes inequalities,
that are used to analyze these quadratures. In the simple Gauss
case, they assert that if $f$ is a function with its first $2n-1$
derivatives positive in $\left( -\infty ,x_{kn}\right) $, then
\[
\sum_{j=k+1}^{n}\lambda _{n}( W^{2},x_{jn}) f( x_{jn}) \leq \int_{-\infty }^{x_{kn}}fW^{2}\leq
\sum_{j=k}^{n}\lambda _{n}( W^{2},x_{jn}) f( x_{jn}) .
\] %
See \cite[p.~33]{Freud1971}, \cite[p.~13]{Mhaskar1996}. This may
be used to prove \cite[p.~17]{Mhaskar1996}:

\begin{lem}  \textsl{Let }$n\geq 1$\textsl{, and }$\xi \in
(x_{k+1,n},x_{kn}]$\textsl{.
Then there exist upper and lower polynomials~}$R_{\xi }$\textsl{\ and }$%
r_{\xi }$\textsl{\ such that }%
\[
r_{\xi }\leq \chi _{(-\infty ,\xi ]}\leq R_{\xi }\quad \mbox{ in
}\mathbb{R}
\] %
\textsl{and}%
\[
\int_{-\infty }^{\infty }\left( R_{\xi }-r_{\xi }\right) W^{2}\leq
\lambda _{n}( W^{2},x_{kn}) +\lambda_{n}(W^{2},x_{k+1,n}).
\]\end{lem}

Thus once we have upper estimates on the Christoffel functions, we
have bounds on the error of one-sided polynomial approximation.
One could write a survey on methods to estimate Christoffel
functions, there are so many. Paul Nevai paid homage to them and
their applications in his still relevant
survey article \cite{Nevai1986}, as well as in his earlier memoir \cite%
{Nevai1979}. We shall present a very simple method of Freud in a
very special case:

\begin{lem}  \textsl{Assume that }$W=\exp (-Q)
\in \mathcal{F}$\textsl{, and
in addition that }$Q$\textsl{\ is convex. Then there exists }$C_{1},C_{2}>0$%
\textsl{\ such that for }$n\geq 1$\textsl{\ and }$\vert \xi \vert
\leq C_{1}a_{n},$%
\begin{equation}
\lambda _{n}( W^{2},\xi ) \leq
C_{2}\frac{a_{n}}{n}W^{2}( \xi ) .
\end{equation}\end{lem}

\proof.  Let
\[
e_{m}(x)=\sum_{j=0}^{m}\frac{x^{j}}{j!}
\] %
denote the $m$th partial sum of the exponential function. Fix $\xi
,n$ and let
\[
R(x)=W^{-1}( \xi ) e_{\floor{n/2}
}\left( ( x-\xi ) Q^{\prime }( \xi ) \right)
.
\] %
Here, $\floor{n/2} $ denotes the integer part of $n/2$. Now
\[
\left\vert e_{\floor{n/2} }(u) \right\vert \leq
C_{2}\exp (u) ,
\] %
for $\vert u\vert \leq \frac{n}{8}$. Moreover, given
$\varepsilon
>0$, we have if $\vert x\vert \leq 2a_{n}$, and $\vert \xi
\vert \leq \varepsilon a_{n}$,%
\begin{equation}
\left\vert ( x-\xi ) Q^{\prime }( \xi )
\right\vert \leq 3a_{n}Q^{\prime }( \varepsilon a_{n})
\leq \frac{n}{8},
\end{equation}%
if only $\varepsilon $ is small enough. To see this, we use the
definition of $a_{n}$ and the monotonicity of $t\mapsto
tQ^{\prime }(t)
:$%
\begin{eqnarray*}
n =\frac{2}{\pi }\int_{0}^{1}\frac{a_{n}tQ^{\prime }( a_{n}t) }{%
\sqrt{1-t^{2}}}dt \geq \frac{2}{\pi }\int_{1/2}^{1}\frac{dt}{\sqrt{1-t^{2}}}\frac{a_{n}}{2}%
Q^{\prime }\left( \frac{a_{n}}{2}\right) .
\end{eqnarray*}%
Thus%
\[
a_{n}Q^{\prime }\left( \frac{a_{n}}{2}\right) \leq C_{1}n.
\] %
Using the lower bound in (3.15) of Definition 3.3, we obtain for
$m\geq 1,$
\[
a_{n}Q^{\prime }\left( \lambda ^{-m}\frac{a_{n}}{2}\right) \leq
A^{-m}a_{n}Q^{\prime }\left( \frac{a_{n}}{2}\right) \leq
C_{1}A^{-m}n.
\] %
There we had $\lambda ,A>1$, so choosing $\varepsilon =\lambda
^{-m}/2$ with large enough $m$, gives (4.3). Next,
\begin{eqnarray*}
\vert R(x)\vert \leq C_{2}W^{-1}( \xi
) \exp ( ( x-\xi ) Q^{\prime }( \xi
)) =C_{2}\exp( Q( \xi ) +( x-\xi ) Q^{\prime }( \xi )) \leq C_{2}\exp
( Q(x)) ,
\end{eqnarray*}%
by convexity of $Q$. Thus
\[
\vert R(x)\vert W(x)\leq
C_{2}\mbox{ for all }\vert x\vert \leq 2a_{n}\mbox{.}
\] %
Moreover,
\[
(RW) ( \xi ) =1.
\] %
Now, we use the extremal property (4.1) of Christoffel functions. We set $%
P=R\rho $, there, where $\rho $ is a polynomial of degree $< \floor{n/2}%
$. We see that
\begin{eqnarray*}
\lambda _{n}( W^{2},\xi ) /W^{2}( \xi )
&=&\inf_{\deg (P) < n}\left[ \int_{-\infty
}^{\infty }(PW)
^{2}\right] /(PW) ^{2}( \xi ) \\[10pt]
&\leq &\inf_{\deg ( \rho ) < \floor{ n/2}
}\left[ \int_{-\infty }^{\infty }( \rho RW)
^{2}\right] /( \rho
RW) ^{2}( \xi ) \\[10pt]
&\leq &C\inf_{\deg ( \rho) < \floor{n/2}
}\left[ \int_{-2a_{n}}^{2a_{n}}( \rho RW) ^{2}\right]
/( \rho RW) ^{2}( \xi) ,
\end{eqnarray*}%
by the restricted range inequality (3.7). Now we use the upper
bound on $RW$
and $(RW)( \xi) =1$ to continue this as%
\begin{eqnarray*}
\leq CC_{2}^{2}\inf_{\deg( \rho ) < \floor{
n/2}}\left( \int_{-2a_{n}}^{2a_{n}}\rho ^{2}\right) /\rho
^{2}( \xi) =CC_{2}^{2}2a_{n}\inf_{\deg (S)
< \floor{ n/2} }\left( \int_{-1}^{1}S^{2}\right)
/S^{2}\left( \frac{\xi }{2a_{n}}\right) .
\end{eqnarray*}%
The latter is the Christoffel function of order $\floor{n/2}
$ for the Legendre weight on $\left[ -1,1\right] $, evaluated at
$\xi /( 2a_{n}) $. Using classical estimates for these \cite[p.~103]{Freud1971}%
, we continue this as
\[
\lambda _{n}( W^{2},\xi) /W^{2}( \xi) \leq C\frac{%
a_{n}}{n}. \qquad\meop
\]

Now we present a very special case of the $L_{1}$
approximation:

\begin{thh}
\textsl{Assume the hypotheses of Lemma $4.2$ on }$W$\textsl{. Assume that }$%
f^{\prime }$\textsl{\ is continuous, and }$f^{\prime }W^{2}\in
L_{1}( \mathbb{R}) $\textsl{, while }$f$\textsl{\ is of polynomial growth at }%
$\infty $\textsl{. Then there exist upper and lower polynomials }$S_{n}$%
\textsl{\ and }$s_{n}$\textsl{\ such that }%
\begin{equation}
s_{n}\leq f\leq S_{n}\mbox{ in }\mathbb{R}
\end{equation}%
\textsl{and}%
\[
\int_{-\infty }^{\infty }( S_{n}-s_{n}) W^{2}\leq C\frac{a_{n}}{n}%
\left( \int_{-\infty }^{\infty }\left\vert f^{\prime }\right\vert
W^{2}+\norm{f^{\prime }W^{2}}_{L_{\infty }( \vert x\vert \geq Ca_{n}) }\right).
\]\end{thh}

\proofi. We assume that $f^{\prime
}\geq 0$, and that $f^{\prime }=0$ outside $\left( x_{nn},x_{1n}\right) $ to simplify the proof. Recall that $x_{nn}$ and $%
x_{1n}$ are the largest and smallest zeros of $p_{n}( W^{2},x) $.
We write the fundamental theorem of calculus%
\[
f(x)=f(0) +\int_{0}^{x}f^{\prime }( \xi ) d\xi
\] %
in the form%
\begin{eqnarray*}
f(x)=f(0) +\int_{0}^{\infty }(1-\chi
_{(-\infty ,\xi ]}(x))f^{\prime }( \xi)
d\xi -\int_{-\infty }^{0}\chi _{(-\infty ,\xi ]}(x)
f^{\prime }( \xi) d\xi .
\end{eqnarray*}%
To check this, consider separately $x\geq 0$ and $x<0$. We use the
upper and lower polynomials $R_{\xi }$ and $r_{\xi }$ of Lemma 4.1
and define
\begin{eqnarray*}
S_{n}(x)=f(0) +\int_{0}^{\infty
}(1-r_{\xi }(x))f^{\prime }( \xi) d\xi
-\int_{-\infty }^{0}r_{\xi }(x)f^{\prime }( \xi
) d\xi
\end{eqnarray*}%
and%
\begin{eqnarray*}
s_{n}(x)=f(0) +\int_{0}^{\infty
}(1-R_{\xi }(x))f^{\prime }(\xi) d\xi
-\int_{-\infty }^{0}R_{\xi }(x)f^{\prime }( \xi
) d\xi .
\end{eqnarray*}%
As $f^{\prime }\geq 0$ and $r_{\xi }\leq \chi _{(-\infty ,\xi
]}\leq R_{\xi } $, we obtain (4.4). Moreover,
\begin{eqnarray*}
( S_{n}-s_{n}) (x)=\int_{-\infty }^{\infty
}(R_{\xi }-r_{\xi })(x)f^{\prime }(\xi)
d\xi
\end{eqnarray*}%
so%
\[
\int_{-\infty }^{\infty }\left( S_{n}-s_{n}\right) W^{2}\leq
\int_{-\infty }^{\infty }\left[ \int_{-\infty }^{\infty }(R_{\xi
}-r_{\xi })(x)W^{2}(x)dx\right]
f^{\prime }(\xi) d\xi .
\] %
Now by Lemmas 4.1 and 4.2, and some other minor considerations,
\[
\int_{-\infty }^{\infty }(R_{\xi }-r_{\xi })(x)
W^{2}(x)dx\leq C\frac{a_{n}}{n}W^{2}( \xi
) ,
\] %
at least for $\vert \xi \vert \leq Ca_{n}$. Thus,
\begin{eqnarray*}
\int_{-\infty }^{\infty }( S_{n}-s_{n}) W^{2}
&\leq& C\frac{a_{n}}{n}\int_{-Ca_{n}}^{Ca_{n}}f^{\prime }W^{2} \\
&+&\left\Vert f^{\prime }W^{2}\right\Vert _{L_{\infty }( \vert x\vert \geq Ca_{n}) }\int_{\left\vert \xi
\right\vert \geq Ca_{n}}\left[ W^{-2}(\xi)
\int_{-\infty }^{\infty }(R_{\xi }-r_{\xi })(x)
W^{2}(x)dx\right] d\xi .
\end{eqnarray*}%
The first term on the right-hand side has the correct form. The
second tail term does not. One can still use estimates for
Christoffel functions, and some other tricks. See
\cite{Freud1974B}, \cite{FreudNevai1973} or \cite[p.~83
ff.]{Mhaskar1996} for more details.  \endproof

This is what Freud proved in 1974 \cite[p.~297]{Freud1974B}, using
these ideas:

\begin{thh} \textsl{Let }$Q$\textsl{\ be an even and
convex function on the real line. Assume that }$Q$\textsl{\ is
increasing and }$Q^{\prime }$\textsl{\ exists
in }$\left( 0,\infty \right) $\textsl{, and that for some }$C_{1},C_{2}>0,$%
\[
1+C_{1}<\frac{Q^{\prime }(2x) }{Q^{\prime }(x)}%
<1+C_{2}.
\] %
\textsl{Let }$r\geq 0$\textsl{\ and }$f^{(r)
}$\textsl{\ be of
bounded variation over every finite interval, and of polynomial growth at }$%
\infty $\textsl{, satisfying for some }$A,B>0$\textsl{\ and integer }$m,$%
\[
\left\vert f(x)\right\vert \leq A+Bx^{2m}.
\] %
\textsl{\ Then there exist upper and lower polynomials
}$S_{n}$\textsl{\ and
}$s_{n}$\textsl{\ of degree }$\leq n$ \textsl{such that (4.4) holds and}%
\[
\int_{-\infty }^{\infty }( S_{n}-s_{n}) W^{2}\leq C\left( \frac{%
a_{n}}{n}\right) ^{r+1}\left( \int_{-\infty }^{\infty
}W^{2}\vert df^{(r) }\vert +A+B\right) .
\] %
\textsl{The constant }$C$\textsl{\ is independent of }$f,n,A,B.$\end{thh}

While this is a one-sided $L_{1}$ result, it also implies
\[
E_{n}[f;W^{2}]_{1}\leq C\left( \frac{a_{n}}{n}\right) ^{r+1}
\left( \int_{-\infty }^{\infty }W^{2}\vert df^{(r) }\vert +A+B\right) .
\] %
This can be extended to other $L_{p}$ spaces, using duality. See \cite%
{MastroianniSzabados2001} for an extension to weights with inner
singularities.

\subsection{Freud's Method involving de la Vall\'ee Poussin Means}

This is based on orthogonal expansions, and for the finite
interval, was used before Freud. Freud was the first to make it
work on infinite intervals
\cite{Freud1976}, \cite{Freud1977}. For a measurable function for which $%
\int_{-\infty }^{\infty }(fW) ^{2}$ is finite, we can
form the
orthonormal expansion%
\[
f\sim \sum_{j=0}^{\infty }c_{j}p_{j},
\] %
where%
\[
c_{j}=\int_{-\infty }^{\infty }fp_{j}W^{2}\mbox{, }\qquad j\geq 0.
\] %
Define the partial sums of the orthonormal expansion,%
\[
S_{n}\left[ f\right] =\sum_{j=0}^{n-1}c_{j}p_{j}.
\] %
The classic minimum property asserts that%
\[
\left\Vert \left( f-S_{n}\left[ f\right] \right) W\right\Vert
_{L_{2}( \mathbb{R}) }=E_{n-1}\left[ f;W\right]
_{2}=\inf_{\deg (P)
\leq n-1}\left\Vert (f-P) W\right\Vert _{L_{2}( \mathbb{R}%
) },
\] %
and the classic Bessel's inequality asserts that
\[
\left\Vert S_{n}\left[ f\right] W\right\Vert _{L_{2}( \mathbb{R}) }=\left( \sum_{j=0}^{n-1}\vert
c_{j}\vert ^{2}\right) ^{1/2}\leq \left\Vert fW\right\Vert
_{L_{2}( \mathbb{R}) }.
\] %
In particular, $\left\{ S_{n}\right\} $ is a uniformly bounded
sequence of operators in a weighted $L_{2}$ space.

Sometimes it is more convenient to use the de la Vall\'ee Poussin
operators
\[
V_{n}\left[ f\right] =\frac{1}{n}\sum_{j=n+1}^{2n}S_{j}\left[
f\right] .
\] %
We see that still
\[
\left\Vert V_{n}\left[ f\right] W\right\Vert _{L_{2}( \mathbb{R}) }\leq \left\Vert fW\right\Vert _{L_{2}(
\mathbb{R}) }
\] %
so the $\{ V_{n}\} $ are uniformly bounded in this weighted $%
L_{2} $ setting. They also reproduce polynomials:%
\[
V_{n}\left[ P\right] =P
\] %
if $P$ is a polynomial of degree $\leq n-1$. However, they are not
projection operators, as $V_{n}\left[ P\right] $ is a polynomial of degree $%
2n-1$ in general. The real advantage is that in spaces other than $L_{2}$, $%
\left\{ V_{n}\right\} $ is often uniformly bounded when $\left\{
S_{n}\right\} $ is not.

Another crucial ingredient are estimates for the Christoffel
functions defined by (4.1). Let us outline some of the main ideas
of Freud's method to
prove the Jackson--Favard inequality%
\begin{equation}
E_{n}\left[ f;W\right] _{p}\leq C\frac{a_{n}}{n}E_{n-1}\left[ f^{\prime };W%
\right] _{p}.
\end{equation}%
(a) We show that $\left\{ V_{n}\right\} $ is a uniformly bounded
sequence of linear operators, first in $L_{\infty }$, and then in
$L_{1}$ (using duality) and then in $L_{p},1<p<\infty $ (using
interpolation).
\newline\noindent (b) We approximate a very special
function, something like a characteristic function, in $L_{1}$,
using the method of the previous section.
\newline\noindent (c) We
prove the estimate in $L_{\infty }$, using $\left\{ V_{n}\right\}
$ and the special function.
\newline\noindent (d) We extend to
$L_{p}$, $1\leq p<\infty .$\newline

In our outline, we shall do this only for $p=\infty $.\newline
The first step is to show, for some $C\neq C(n,f)$, that
\begin{equation}
E_{2n}\left[ f;W\right] _{\infty }\leq \left\Vert \left( f-V_{n}\left[ f%
\right] \right) W\right\Vert _{L_{\infty }( \mathbb{R})
}\leq CE_{n}\left[ f;W\right] _{\infty }.
\end{equation}%

\begin{lem}  \textsl{The inequality $(4.6)$ holds, provided we
assume the estimate $(4.13)$ below.}
\end{lem}

\proof.
The left inequality here is immediate, as $V_{n}[f]$ is a polynomial of degree $%
\leq 2n-1$. It is the right-hand inequality that requires work. We
follow Freud \cite{Freud1977}. The idea goes back at least to
Torsten Carleman. The
partial sum $s_{m}\left[ f\right] $ admits the representation%
\begin{equation}
s_{m}\left[ f\right] (x)
=\sum_{j=0}^{m-1}c_{j}p_{j}(x)=\int_{-\infty
}^{\infty }f(t) K_{m}( x,t) W^{2}(t) dt,
\end{equation}%
where
\[
K_{m}( x,t) :=\sum_{j=0}^{m-1}p_{j}(x)
p_{j}(t) .
\] %
The Christoffel--Darboux formula \cite{Freud1971} asserts that
\begin{equation}
K_{m}( x,t) =\frac{\gamma _{m-1}}{\gamma
_{m}}\frac{p_{m}(x)p_{m-1}(t)
-p_{m-1}(x)p_{m}(t) }{x-t}.
\end{equation}%
We now fix $x$ and $n$, and define
\[
g=f(x)\chi _{\left( x-\frac{a_{n}}{n},x+\frac{a_{n}}{n}\right) }
\] %
and
\[
h=f-g=f\chi _{\mathbb{R}\backslash \left( x-\frac{a_{n}}{n},x+\frac{a_{n}}{n}%
\right) }.
\] %
Here $\chi _{S}$ denotes the characteristic function of a set $S$.
Thus $g=f$ in an interval of radius $\frac{a_{n}}{n}$ about $x$,
and is $0$ elsewhere, while $h$ is the \textquotedblleft
rest\textquotedblright\ of $f$. We have
\begin{equation}
s_{m}\left[ f\right] (x)=s_{m}\left[ g\right] (x) +s_{m}\left[ h\right] (x).
\end{equation}%
Now%
\begin{eqnarray*}
\left\vert s_{m}\left[ g\right] \right\vert (x)
&=&\left\vert \int_{x-\frac{a_{n}}{n}}^{x+\frac{a_{n}}{n}}f(t) K_{m}( x,t) W^{2}(t) dt\right\vert
\leq \left\Vert fW\right\Vert _{L_{\infty }( \mathbb{R})
}\int_{x-\frac{a_{n}}{n}}^{x+\frac{a_{n}}{n}}\left\vert
K_{m}( x,t) \right\vert W(t) dt \\[10pt]
&\leq &\left\Vert fW\right\Vert _{L_{\infty }( \mathbb{R})
}\left( \frac{2a_{n}}{n}\right) ^{1/2}\left( \int_{x-\frac{a_{n}}{n}}^{x+%
\frac{a_{n}}{n}}K_{m}^{2}( x,t) W^{2}(t)
dt\right) ^{1/2},
\end{eqnarray*}%
by the Cauchy--Schwarz inequality. Here, using orthogonality,
\begin{eqnarray*}
\int_{x-\frac{a_{n}}{n}}^{x+\frac{a_{n}}{n}}K_{m}^{2}\left( x,t\right) W^{2}(t) dt &\leq &\int_{-\infty }^{\infty
}K_{m}^{2}( x,t) W^{2}(t) dt \\
&=&\sum_{j=0}^{m-1}p_{j}^{2}(x)\leq
\sum_{j=0}^{2n-1}p_{j}^{2}(x)=\lambda
_{2n}^{-1}( W^{2},x) ,
\end{eqnarray*}%
provided $m\leq 2n$. Thus for $m\leq 2n,$%
\[
\left\vert s_{m}\left[ g\right] (x)\right\vert \leq
\left\Vert
fW\right\Vert _{L_{\infty }( \mathbb{R}) }\left( \frac{2a_{n}}{n}%
\lambda _{2n}^{-1}( W^{2},x) \right) ^{1/2},
\] %
and hence, averaging over $s_{m},$ $n+1\leq m\leq 2n,$
\begin{equation}
\left\vert V_{n}\left[ g\right] (x)\right\vert \leq
\left\Vert
fW\right\Vert _{L_{\infty }( \mathbb{R}) }\left( \frac{2a_{n}}{n}%
\lambda _{2n}^{-1}( W^{2},x) \right) ^{1/2}.
\end{equation}%
Now comes the clever idea. Let
\[
H(t) =\frac{h(t) }{x-t},\qquad t\in
\mathbb{R},
\] %
and denote its Fourier coefficients with respect to $\left\{
p_{j}\right\} $ by $\left\{ c_{j}\left[ H\right] \right\} $, so
that
\[
H\sim \sum_{j=0}^{\infty }c_{j}\left[ H\right] p_{j}.
\] %
We use the Christoffel--Darboux formula (4.8) to write%
\begin{eqnarray*}
\left\vert s_{m}\left[ h\right] (x)\right\vert
&=&\frac{\gamma _{m-1}}{\gamma _{m}}\left\vert p_{m}(x)
\int_{-\infty }^{\infty }H(t) p_{m-1}(t) W^{2}(t) dt-p_{m-1}(x)
\int_{-\infty }^{\infty }H(t)
p_{m}(t) W^{2}(t) dt\right\vert \\[10pt]
&=&\frac{\gamma _{m-1}}{\gamma _{m}}\left\vert p_{m}(x)c_{m-1}%
\left[ H\right] -p_{m-1}(x)c_{m}\left[ H\right]
\right\vert .
\end{eqnarray*}%
Then summing and using Cauchy--Schwarz,
\begin{equation}
\left\vert V_{n}\left[ h\right] (x)\right\vert \leq \frac{1}{n}%
\left[ \max_{m\leq 2n}\frac{\gamma _{m-1}}{\gamma _{m}}\right]
2\left( \sum_{m=0}^{2n}p_{m}^{2}(x)\right)
^{1/2}\left( \sum_{m=0}^{2n}c_{m}^{2}(H) \right)
^{1/2}.
\end{equation}%
Yet another clever idea: use Bessel's inequality on the Fourier
coefficients of $H$ (recall, this holds for any orthogonal system
in any inner product space),
\begin{eqnarray*}
\sum_{m=0}^{2n}c_{m}^{2}(H) &\leq &\sum_{m=0}^{\infty
}c_{m}^{2}(H) \leq \int_{-\infty }^{\infty }(HW) ^{2} =\int_{\left\{ t:\left\vert t-x\right\vert \geq
\frac{a_{n}}{n}\right\} }\left( \frac{f(t) W(t) }{t-x}\right) ^{2}dt
\leq \left\Vert fW\right\Vert _{L_{\infty }( \mathbb{R}) }^{2}2%
\frac{n}{a_{n}}.
\end{eqnarray*}%
Then (4.11) becomes%
\[
\left\vert V_{n}\left[ h\right] (x)\right\vert \leq \frac{2%
\sqrt{2}}{n}\left[ \max_{m\leq 2n}\frac{\gamma _{m-1}}{\gamma
_{m}}\right] \left( \frac{n}{a_{n}}\lambda _{2n+1}^{-1}( W^{2},x)
\right) ^{1/2}.
\] %
Finally, combining this with (4.9) and (4.10) gives%
\begin{eqnarray*}
&&\left\vert V_{n}\left[ f\right] (x)\right\vert \leq \left\Vert
fW\right\Vert _{L_{\infty }( \mathbb{R}) }\left( \lambda _{2n+1}(
W^{2},x) \right) ^{-1/2}\left\{ \left( \frac{2a_{n}}{n}\right)
^{1/2}+\frac{1}{n}\left[ \max_{m\leq 2n}\frac{\gamma
_{m-1}}{\gamma _{m}}\right] 2\sqrt{2}\left( \frac{n}{a_{n}}\right)
^{1/2}\right\} .
\end{eqnarray*}%
\begin{equation}
\end{equation}%
Up to this stage, we have not used any properties of the weight,
it's
completely general. But now we use lower bounds for the Christoffel function%
\begin{equation}
\lambda _{n}( W^{2},x) \geq C\frac{a_{n}}{n}W^{2}(x), \qquad x\in \mathbb{R},\ n\geq 1,
\end{equation}%
and an upper bound
\begin{equation}
\max_{m\leq 2n}\frac{\gamma _{m-1}}{\gamma _{m}}\leq Ca_{n}
\end{equation}%
to deduce that uniformly for $n\geq 1$, and $x\in \mathbb{R}$,
\[
|V_{n}\left[ f\right] (x)|W(x)\leq C\left\Vert fW\right\Vert
_{L_{\infty }( \mathbb{R}) }.
\] %
Thus, for some $C_{\ast }\neq C( n,f)$,%
\[
\left\Vert V_{n}\left[ f\right] W\right\Vert _{L_{\infty }( \mathbb{R}%
) }\leq C_{\ast }\left\Vert fW\right\Vert _{L_{\infty }( \mathbb{R%
}) }.
\] %
Then for any polynomial $P$ of degree $\leq n$, the reproducing property of $%
V_{n}$ gives
\begin{eqnarray*}
\left\Vert \left( f-V_{n}\left[ f\right] \right) W\right\Vert
_{L_{\infty }( \mathbb{R}) } &=&\left\Vert \left( f-P-V_{n}\left[ f-P\right] \right) W\right\Vert
_{L_{\infty }( \mathbb{R}) } \\[10pt]
&\leq &\left\Vert (f-P) W\right\Vert _{L_{\infty }( \mathbb{%
R}) }+\left\Vert V_{n}\left[ f-P\right] W\right\Vert
_{L_{\infty
}( \mathbb{R}) } \\[10pt]
&\leq &( 1+C_{\ast }) \left\Vert (f-P)
W\right\Vert _{L_{\infty }( \mathbb{R}) }.
\end{eqnarray*}%
Taking the inf over all $P$ of degree $\leq n$ gives the
right-hand inequality in (4.6).

We turn to the discussion of (4.13) and (4.14). For the
latter, we can just use restricted range inequalities: if $m\leq 2n$%
\begin{eqnarray*}
\frac{\gamma _{m-1}}{\gamma _{m}} &=&\int_{-\infty }^{\infty
}xp_{m-1}(x) p_{m}(x)W^{2}(x)dx \\[10pt]
&\leq &C\int_{-2a_{m}}^{2a_{m}}\left\vert xp_{m-1}(x)
p_{m}(x)\right\vert W^{2}(x)dx \leq
C2a_{m},
\end{eqnarray*}%
by Cauchy--Schwarz. The lower bound (4.13) is more difficult; see
for example \cite{LevinLubinsky1992}, \cite{Mhaskar1996},
\cite{Nevai1986}. \endproof

As an aside, let us see how duality also gives this bound in
$L_{1}$. We use
duality of $L_{p}$ norms:%
\begin{equation}
\left\Vert V_{n}\left[ f\right] W\right\Vert _{L_{1}( \mathbb{R}) }=\sup_{g}\int_{-\infty }^{\infty }V_{n}\left[
f\right] gW^{2},
\end{equation}%
where the sup is over all measurable $g$ with $\left\Vert
gW\right\Vert _{L_{\infty }( \mathbb{R}) }\leq 1$.
(Sorry, we used $g$ above in
a different sense.) Now we use the self-adjointness of $V_{n}:$%
\begin{equation}
\int_{-\infty }^{\infty }V_{n}\left[ f\right] gW^{2}=\int_{-\infty
}^{\infty }fV_{n}\left[ g\right] W^{2}.
\end{equation}%
This follows easily, once we prove the self-adjointness of $s_{m}:$%
\[
\int_{-\infty }^{\infty }s_{m}\left[ f\right] gW^{2}=\int_{-\infty
}^{\infty }fs_{m}\left[ g\right] W^{2}.
\] %
We leave the latter as an exercise (just substitute in the definitions of $%
s_{m}\left[ f\right] $ and $s_{m}\left[ g\right] $). Combining
(4.15) and
(4.16) gives%
\begin{eqnarray}
\left\Vert V_{n}\left[ f\right] W\right\Vert _{L_{1}( \mathbb{R}) } &=&\sup_{g}\int_{-\infty }^{\infty
}fV_{n}\left[ g\right] W^{2} \leq \sup_{g}\left\Vert fW\right\Vert
_{L_{1}( \mathbb{R})
}\left\Vert V_{n}\left[ g\right] W\right\Vert _{L_{\infty }( \mathbb{R}%
) }  \nonumber\\[10pt]
&\leq& C\left\Vert fW\right\Vert _{L_{1}( \mathbb{R})
},
\end{eqnarray}%
using the case $p=\infty $. With this bound, we can finish off the
rest of the proof as we did for the case $p=\infty $. For
$1<p<\infty $, one can use interpolation of operators. It is
possible to improve the estimate (4.6) to include factors that
vanish near $\pm a_{n}$, reflecting improved approximation there,
see \cite{LubinskyMache2000}. A partially successful
attempt to extend (4.5) to general exponential weights was given in \cite%
{LubinskyMashele2002}.

An obvious question is for which weights, we have the lower bound
(4.13) for
the Christoffel functions. They have not been established for the class $%
\mathcal{F}$ of Definition 3.3. However, they were established in \cite%
{LevinLubinsky1992} for the following class:

\begin{deff}
\textsl{Let }$W=\exp (-Q) $\textsl{, where }$Q^{\prime \prime }$%
\textsl{\ exists and is positive in }$\left( 0,\infty \right)
$\textsl{, while }$Q^{\prime }$\textsl{\ is positive there, with
limit }$0$\textsl{\ at
}$0$\textsl{, and for some }$A,B>1,$%
\begin{equation}
A-1\leq \frac{xQ^{\prime \prime }(x)}{Q^{\prime
}(x)}\leq B-1,\mbox{ }\qquad x\in \left( 0,\infty
\right) .
\end{equation}%
\textsl{Then we write }$W\in \mathcal{F}^{\ast }$.\end{deff}

For weaker (but difficult to formulate) hypotheses, these
estimates were
proved in \cite{LevinLubinsky2001}. An excellent exposition is given in \cite%
[Chapter 3]{Mhaskar1996}. The Christoffel function bound (4.13)
there is
proved assuming $Q^{\prime \prime }$ is increasing. This is true for $%
W=W_{\alpha }$ if $\alpha \geq 2$, while (4.18) holds for $\alpha
>1$.

The next step is an inequality that is a distant relative of
Hardy's inequality:

\begin{lem}
\textsl{Let }$n\geq 1$\textsl{\ and }$g$\textsl{\ be a function such that }%
\begin{equation}
\int_{-\infty }^{\infty }gPW^{2}=0,
\end{equation}%
\textsl{for all polynomials }$P$ \textsl{of degree }$\leq
n$\textsl{. Then there exists }$C\neq C( n,g)
$\textsl{\ such that}
\begin{equation}
\sup_{x\in \mathbb{R}}W(x)\left\vert
\int_{0}^{x}g\right\vert
\leq C\frac{a_{n}}{n}\left\Vert gW\right\Vert _{L_{\infty }( \mathbb{R}%
) }.
\end{equation}\end{lem}

\proofi. Fix $x>0$, and let
\[
\phi _{x}(t) =W^{-2}(t) \chi _{\left[
0,x\right] }(t) .
\] %
For $g$ as above, and any polynomial $P$ of degree $\leq n,$%
\begin{eqnarray*}
\left\vert \int_{0}^{x}g\right\vert &=&\left\vert \int_{-\infty
}^{\infty }g(t) \phi _{x}(t) W^{2}(t) dt\right\vert
\\[10pt]
&=&\left\vert \int_{-\infty }^{\infty }g(t) \left[
\phi _{x}(t) -P(t) \right] W^{2}(t)
dt\right\vert \\[10pt]
&\leq &\left\Vert gW\right\Vert _{L_{\infty }( \mathbb{R}) }\int_{-\infty }^{\infty }\left\vert \phi
_{x}-P\right\vert W.
\end{eqnarray*}%
Taking the inf over all such $P$ gives%
\[
\left\vert \int_{0}^{x}g\right\vert \leq \left\Vert gW\right\Vert
_{L_{\infty }( \mathbb{R}) }E_{n}\left[ \phi
_{x};W\right] _{1}.
\] %
Once we have the estimate
\begin{equation}
W(x)E_{n}\left[ \phi _{x};W\right] _{1}\leq
C\frac{a_{n}}{n},
\end{equation}%
for some $C\neq C( n,x) $, the result follows. To prove
this, one
can use one-sided approximation as in Section 4.1. See, for example, \cite[p.~34]%
{Freud1977}. \endproof

Now we can give the Jackson--Favard inequality in the case $p=\infty $:%

\begin{thh} \textsl{Let }$W\in \mathcal{F}^{\ast
}$\textsl{. Let }$f$\textsl{\ be absolutely continuous in each
finite interval, with }$\left\Vert f^{\prime }W\right\Vert
_{L_{\infty }( \mathbb{R}) }$\textsl{\ finite. Then
for some }$C\neq C( n,f) ,$\textsl{\ }%
\begin{equation}
E_{n}\left[ f;W\right] _{\infty }\leq C\frac{a_{n}}{n}\left\Vert
f^{\prime }W\right\Vert _{L_{\infty }( \mathbb{R}) }
\end{equation}%
\textsl{and}%
\begin{equation}
E_{n}\left[ f;W\right] _{\infty }\leq
C\frac{a_{n}}{n}E_{n-1}\left[ f^{\prime };W\right] _{\infty }.
\end{equation}\end{thh}

\proof. Let
\[
g(x)=f^{\prime }(x)-V_{n}\left[
f^{\prime }\right] (x)\mbox{.}
\] %
This does satisfy (4.19). Indeed, if $P$ is a polynomial of degree $\leq $ $%
n $, self-adjointness of $V_{n}$ gives
\begin{eqnarray*}
\int_{-\infty }^{\infty }gPW^{2} &=&\int_{-\infty }^{\infty
}f^{\prime
}PW^{2}-\int_{-\infty }^{\infty }V_{n}\left[ f^{\prime }\right] PW^{2} \\[10pt]
&=&\int_{-\infty }^{\infty }f^{\prime }PW^{2}-\int_{-\infty
}^{\infty }f^{\prime }V_{n}\left[ P\right] W^{2}=0,
\end{eqnarray*}%
as $V_{n}\left[ P\right] =P$. Next, let
\[
U_{n}(x)=f(0) +\int_{0}^{x}V_{n}\left[ f^{\prime }%
\right] (t) dt,
\] %
so that
\[
( f-U_{n}) (x)=\int_{0}^{x}\left( f^{\prime }-V_{n} \left[ f^{\prime }\right] \right)
=\int_{0}^{x}g.
\] %
Then
\begin{eqnarray*}
\left\Vert W\left( f-U_{n}\right) \right\Vert _{L_{\infty }\left( \mathbb{R}%
\right) } &=&\left\Vert W\int_{0}^{x}g\right\Vert _{L_{\infty
}(\mathbb{R}) }
\leq C\frac{a_{n}}{n}\left\Vert gW\right\Vert _{L_{\infty }\left( \mathbb{R%
}\right) } \\[10pt]
&=&C\frac{a_{n}}{n}\left\Vert \left( f^{\prime }-V_{n}\left[ f^{\prime }%
\right] \right) W\right\Vert _{L_{\infty }(\mathbb{R})
} \leq \frac{a_{n}}{n}E_{n}\left[ f^{\prime };W\right] _{\infty }.
\end{eqnarray*}%
Here we applied Lemma 4.7 and then Lemma 4.5. This also gives, as
$U_{n}$
has degree $\leq 2n,$%
\[
E_{2n}\left[ f;W\right] _{\infty }\leq C\frac{a_{n}}{n}E_{n}\left[
f^{\prime };W\right] _{\infty }\leq C\frac{a_{n}}{n}\left\Vert
f^{\prime }W\right\Vert _{L_{\infty }(\mathbb{R}) }.
\] %
Replacing $n$ by $n/2$, and using the fact that $a_{n/2}\leq
a_{n}$, we
obtain%
\[
E_{n}\left[ f;W\right] _{\infty }\leq C\frac{a_{n}}{n}\left\Vert
f^{\prime }W\right\Vert _{L_{\infty }(\mathbb{R}) }.
\] %
Finally, we observe that for any polynomial $P$ of degree $\leq n,$%
\begin{eqnarray*}
E_{n}\left[ f;W\right] _{\infty } =E_{n}\left[ f-P;W\right]
_{\infty } \leq C\frac{a_{n}}{n}\left\Vert (f-P)
^{\prime }W\right\Vert _{L_{\infty }(\mathbb{R}) }.
\end{eqnarray*}%
Choosing $P^{\prime }$ suitably gives (4.23).  \endproof

For the extension to all $1\leq p\leq \infty $, see \cite[Chapter 4]%
{Mhaskar1996}.

\subsection{The Kro\'{o}--Szabados Method}

The idea here \cite{KrooSzabados1995} is to make use of the
alternation/equioscillation for best polynomial approximants,
together with some clever tricks. On a finite interval, the idea
was used by Bojanov \cite{Bojanov1996}, Babenko and Shalaev. While
it works quite generally, it does not yield the correct Jackson
rate. Maybe some clever tweaking can repair that?

Let us fix a function $f$, and $n\geq 1$, and let $P_{n}^{\ast }$
be its best polynomial approximation of degree $\leq n$ in the
weighted uniform
norm. Thus%
\[
\left\Vert (f-P_{n}^{\ast })W\right\Vert _{L_{\infty }\left( \mathbb{R}%
\right) }=\inf_{\deg (P) \leq n}\left\Vert
(f-P)W\right\Vert _{L_{\infty }(\mathbb{R})
}=E_{n}\left[ f;W\right] _{\infty }.
\] %
Then there exist equioscillation points $\left\{ y_{j}\right\}
_{j=0}^{n+1}$ such that for $0\leq j\leq n+1,$
\[
\left[ (f-P_{n}^{\ast })W\right] \left( y_{j}\right) =\varepsilon
(-1) ^{j}E_{n}\left[ f;W\right] _{\infty }.
\] %
The number $\varepsilon \in \left\{ -1,1\right\} $ is independent
of $j$. In
terms of these, there is the determinant expression%
\[
E_{n}\left[ f;W\right] _{\infty }=\left\vert \frac{U\left( \begin{array}{rrrr}
f & 1 & \ldots & x^{n} \\
y_{0} & y_{1} & \ldots & y_{n+1}%
\end{array}%
\right) }{U\left( \begin{array}{rrrr}
g & 1 & \ldots & x^{n} \\
y_{0} & y_{1} & \ldots & y_{n+1}%
\end{array}%
\right) }\right\vert ,
\] %
where for functions $\left\{ \phi _{j}\right\} $%
\[
U\left( \begin{array}{rrrr}
\phi _{0} & \phi _{1} & \ldots & \phi _{n+1} \\
y_{0} & y_{1} & \ldots & y_{n+1}%
\end{array}%
\right) :=\det \left( \phi _{i}\left( y_{j}\right) \right)
_{i,j=0}^{n+1}
\] %
and $g$ is a function such that
\[
(gW) \left( y_{j}\right) =(-1) ^{j},\qquad
0\leq j\leq n+1.
\] %
A proof of this formula may be found in \cite[p.~28]{Shapiro1971}.
By some elementary determinantal manipulations, one can show that
\[
U\left( \begin{array}{rrrr}
f & 1 & \ldots & x^{n} \\
y_{0} & y_{1} & \ldots & y_{n+1}%
\end{array}%
\right) =(-1) ^{n}\sum_{k=0}^{n}(-1)
^{k}\left[ f\left( y_{k+1}\right) -f\left( y_{k}\right) \right]
B_{k}
\] %
and%
\[
U\left( \begin{array}{rrrr}
g & 1 & \ldots & x^{n} \\
y_{0} & y_{1} & \ldots & y_{n+1}%
\end{array}%
\right) =(-1) ^{n+1}\sum_{k=0}^{n}\left[ W^{-1}\left( y_{k+1}\right) +W^{-1}\left( y_{k}\right) \right] B_{k}.
\] %
Here $B_{k}$ is the determinant of a matrix with entries involving only $%
\left\{ y_{j}\right\} $. Then one obtains
\begin{equation}
E_{n}\left[ f;W\right] _{\infty }=\left\vert \sum_{k=0}^{n}(-1) ^{k}\left[ f\left( y_{k+1}\right) -f\left( y_{k}\right)
\right] d_{k}\right\vert ,
\end{equation}%
where%
\[
d_{k}=\frac{B_{k}}{\sum_{j=0}^{n}\left[ W^{-1}\left( y_{j+1}\right) +W^{-1}\left( y_{j}\right) \right] B_{j}}.
\] %
Suppose we define the (unusual!) modulus
\[
\omega _{\gamma }(f;t)=\sup_{x,y\in \mathbb{R}\mbox{, }\left\vert
x-y\right\vert \leq t}\frac{\left\vert f(x)-f(y) \right\vert }{W(x)^{-\gamma }+W\left(
y\right) ^{-\gamma }},
\] %
where $\gamma \in \left( 0,1\right) $. Then from (4.24),
\[
E_{n}\left[ f;W\right] _{\infty }\leq \sum_{k=0}^{n}\left\vert
d_{k}\right\vert \left( W\left( y_{k+1}\right) ^{-\gamma }+W\left( y_{k}\right) ^{-\gamma }\right) \omega _{\gamma }\left(
f;\left\vert y_{k+1}-y_{k}\right\vert \right) .
\] %
Now one uses properties of the modulus $\omega _{\gamma }$ and
then has to estimate the $\left\{ d_{k}\right\} $. This involves
tricks such as needle polynomials. Here is a sample of what can be
achieved:

\begin{thh} \textsl{Let }$Q:\mathbb{R}\rightarrow
\mathbb{R}$\textsl{\ be an even
continuous function, which is positive and differentiable for large }$x$%
\textsl{, and with }%
\begin{equation}
0<\liminf_{x\rightarrow \infty }\frac{xQ^{\prime }(x)
}{Q(x)}\leq \limsup_{x\rightarrow \infty
}\frac{xQ^{\prime }(x)}{Q(x)}<\infty .
\end{equation}%
\textsl{Let }$Q^{\left[ -1\right] }$\textsl{\ denote the inverse of }$Q$%
\textsl{, defined for sufficiently large positive }$x$\textsl{. For small }$%
n $\textsl{, take }$Q^{\left[ -1\right] }(n) $\textsl{\ to be }$%
1 $\textsl{. Assume that }$0<\gamma <1$\textsl{, and }$f:\mathbb{R}%
\rightarrow \mathbb{R}$\textsl{\ is continuous, with }%
\[
\lim_{\vert x\vert \rightarrow \infty }f(x)
W^{\gamma }(x)=0.
\] %
\begin{description}
\item[(a)] \textsl{Then for some }$C_{j}\neq C_{j}\left( f,n\right) ,$%
\begin{equation}
E_{n}\left[ f;W\right] _{\infty }\leq C_{1}\omega _{\gamma }(f;\frac{Q^{%
\left[ -1\right] }(n) \log n}{\left( 1-\gamma \right) n}%
)+e^{-C_{2}n}\left\Vert fW^{\gamma }\right\Vert _{L_{\infty }\left( \mathbb{R%
}\right) }.
\end{equation}%
\item[(b)]\textsl{Let }$0<\varepsilon <\frac{1-\gamma }{2-\gamma
}$\textsl{. Then
for some }$C\neq C\left( f,n\right) ,$%
\begin{equation}
E_{n}\left[ f;W\right] _{\infty }\leq C_{1}\omega _{\gamma }(f;I_{n}^{-\frac{%
1-\gamma }{2-\gamma }+\varepsilon }),
\end{equation}%
\textsl{where for large enough }$n$,%
\begin{equation}
I_{n}:=\int_{1}^{Q^{\left[ -1\right] }(n)
}\frac{Q(t) }{t^{2}}dt.
\end{equation}
\end{description}\end{thh}

Note that (4.25) requires less than we required for $W\in
\mathcal{F}$, but
the restrictions on $f$ are more severe than in Theorem 3.5, and instead of $%
\frac{a_{n}}{n}$ inside the modulus, we obtain essentially $\frac{a_{n}\log n%
}{n}$. The boundary case $Q(x)=\left\vert
x\right\vert $
satisfies the above conditions, and in this case we obtain%
\[
E_{n}\left[ f;W\right] _{\infty }\leq C_{1}\omega _{\gamma
}(f;\left( \log n\right) ^{-\frac{1-\gamma }{2-\gamma
}+\varepsilon }).
\] %
Here by choosing $\gamma $ small enough, the exponent of $\log n$
can be made arbitrarily close to $-1$.

This method is interesting, and general. The challenge is how to
tweak it, if possible, to get the correct Jackson rate.

\subsection{The Piecewise Polynomial Method}

This is undoubtedly the most direct and general method, and
Ditzian and the author used it to prove Theorem 3.5. However, it
does pose substantial technical challenges. For finite intervals,
it has been used for a long time, and in spirit goes back to
Lebesgue's proof of Weierstrass' Theorem. Lebesgue first
approximated by a piecewise linear function, and then polynomials.
It has served as a powerful tool on finite intervals, for example
in investigating shape preserving polynomial approximation, the
degree of spline approximation, and even approximation by rational
functions \cite{DeVoreetal1992}, \cite{PetrushevPopov1987}.

The function $f$ is first approximated by a piecewise polynomial
(or spline). Each of the piecewise polynomials is generated via
Whitney's Theorem. Then special polynomials that approximate
characteristic functions are used to turn the spline approximation
into a polynomial approximation. We illustrate the method as it is
used to prove Theorem 3.5 for Freud weights.
\smallskip

\noindent
\textbf{Step 1: Partition} $\left[
-a_{n},a_{n}\right]$. Recall that our modulus $\omega _{r,p}$
involves a tail piece that will take care of the behavior of
$f(x)$ for very large $x$. So we fix $n$ and
concentrate on approximation on the Mhaskar--Rakhmanov--Saff
interval $\left[ -a_{n},a_{n}\right] $. We partition this interval
into
small intervals, all of length $\frac{a_{n}}{n}:$%
\[
-a_{n}=\tau _{0}<\tau _{1}<\tau _{2}<\cdots<\tau _{2n}=a_{n}.
\] %
Set
\[
I_{j}:=\left[ \tau _{j},\tau _{j+1}\right], \qquad 0\leq j\leq
2n-1.
\] %
\smallskip

\noindent
\textbf{Step 2: Use Whitney's Theorem to develop a piecewise
polynomial approximation.}
We fix $r$ --- this will be the order of the modulus. We approximate $f$ on $%
I_{j}$ by a polynomial $p_{j}$ of degree $\leq r$, and then form
the
piecewise polynomial%
\[
S\left[ f\right] :=p_{0}+\sum_{j=1}^{2n-1}\left( p_{j}-p_{j-1}\right) \chi _{%
\left[ \tau _{j},a_{n}\right] },
\] %
where, as usual, $\chi _{\left[ \tau _{j},a_{n}\right] }$ denotes
the characteristic function of the interval $\left[ \tau
_{j},a_{n}\right] $. When restricted to $I_{k}$, the sum becomes a
telescopic sum: we see that in
the interior of $I_{k}$, namely in $(\tau _{k},\tau _{k+1}),$%
\[
S\left[ f\right] =p_{0}+\sum_{j=1}^{k}\left( p_{j}-p_{j-1}\right)
=p_{k}.
\] %
Then, if $p<\infty ,$
\begin{eqnarray}
\left\Vert \left( f-S\left[ f\right] \right) W\right\Vert
_{L_{p}\left[ -a_{n},a_{n}\right] }^{p}
=\sum_{j=0}^{2n-1}\int_{I_{j}}\left\vert \left( f-p_{j}\right)
W\right\vert ^{p} &\leq &C\sum_{j=0}^{2n-1}W\left( \tau
_{j}\right) ^{p}\left\Vert f-p_{j}\right\Vert _{L_{p}\left( I_{j}\right) }^{p}.
\end{eqnarray}%
In this step, we use the fact that because the intervals $I_{j}$
have length $\frac{a_{n}}{n}$, $W$ does not grow or decay by more
than a constant. The
idea is that for $x,y\in I_{k},$%
\begin{eqnarray}
W(x)/W(y) =\exp( Q(y)
-Q(x)) =\exp( Q^{\prime }( \xi
) (x-y)) \leq \exp \left( Q^{\prime
}(\xi) \frac{a_{n}}{n}\right) \leq C,
\end{eqnarray}%
since $Q^{\prime }$ is bounded by $C\frac{n}{a_{n}}$ throughout
$\left[
-a_{n},a_{n}\right] $. (This is true for Freud weights in the class $%
\mathcal{F}^{\ast }$ but fails for Erd\H{o}s weights.) Now comes
the application of Whitney. Let
\[
I_{j}^{\ast }=I_{j}\cup I_{j+1}.
\] %
By Whitney's Theorem on the interval $I_{j}^{\ast }$ \cite[p.~195, p.~191]%
{PetrushevPopov1987}, we can choose $p_{j}$ of degree $\leq r$
such that
\[
\left\Vert f-p_{j}\right\Vert _{L_{p}\left( I_{j}^{\ast }\right) }^{p}\leq C%
\frac{n}{a_{n}}\int_{0}^{a_{n}/n}\int_{I_{j}^{\ast }}\left\vert
\Delta _{s}^{r}f(x)\right\vert ^{p}dx\mbox{
}ds=:C\Omega _{j}^{p}.
\] %
Here $C$ is independent of $f,n,j$. The strange creature on the
right is really the $p$th power of an $r$th order integral modulus
of continuity on the interval $I_{j}^{\ast }$. In forming the
$r$th difference in this integral, one uses the convention that
the difference is taken as $0$ if any of the arguments of the
function are outside the interval $I_{j}^{\ast }$.
Substituting this in (4.29) and using (4.30) gives%
\begin{eqnarray}
\left\Vert \left( f-S\left[ f\right] \right) W\right\Vert
_{L_{p}\left[
-a_{n},a_{n}\right] }^{p} &\leq &C\frac{n}{a_{n}}\int_{0}^{a_{n}/n}%
\int_{-a_{n}}^{a_{n}}\left\vert W\Delta _{s}^{r}f(x)
\right\vert
^{p}dx\mbox{ }ds  \nonumber \\[10pt]
&\leq &C\sup_{0<h\leq \frac{a_{n}}{n}}\left\Vert W\Delta
_{h}^{r}f\right\Vert _{L_{p}\left[ -a_{n},a_{n}\right] }^{p}.
\end{eqnarray}%
\smallskip

\noindent\textbf{Step 3: Approximate the characteristic function }$\chi
_{\left[ \tau _{j},a_{n}\right]}$.
Now comes the difficult part. We approximate $\chi _{\left[ \tau _{j},a_{n}%
\right] }$ by a polynomial $R_{j}$ of degree $\leq Ln$ giving the polynomial%
\[
P\left[ f\right] =p_{0}+\sum_{j=1}^{2n-1}\left( p_{j}-p_{j-1}\right) R_{j},
\] %
which is of degree at most $Ln+r$. Here the constant $L$ is independent of $%
n $, and arises in the challenging task of generating $R_{j}$. We
see that then
\begin{equation}
S(f) -P\left[ f\right] =\sum_{j=1}^{2n-1}\left( p_{j}-p_{j-1}\right) \left( \chi _{\left[ \tau _{j},a_{n}\right]
}-R_{j}\right) .
\end{equation}%
To estimate this, we compare $p_{j}$ and $p_{j-1}$ on the common
interval
where they approximate $f$, namely $I_{j}$. (That is why we used Whitney on $%
I_{j}^{\ast }=I_{j}\cup I_{j+1}$.) We obtain, even for $p<1,$%
\[
\left\Vert p_{j}-p_{j-1}\right\Vert _{L_{p}\left( I_{j}\right)
}\leq C\left( \Omega _{j-1}+\Omega _{j}\right) .
\] %
Now we use Nikolskii inequalities, which compare the norms of
polynomials in $L_{p}$ and $L_{q}$, and the Bernstein--Walsh
inequality, which bounds the growth of polynomials outside an
interval, once we know their size on the
interval. Together they yield for all real $x,$%
\begin{equation}
\left\vert p_{j}-p_{j-1}\right\vert (x)\leq C\left( \frac{n}{%
a_{n}}\right) ^{1/p}\left( 1+\frac{n}{a_{n}}\left\vert x-\tau
_{j}\right\vert \right) ^{r}\left( \Omega _{j-1}+\Omega
_{j}\right) .
\end{equation}%
Again, the constant $C$ does not depend on $x$, or $n$, or $f$. It
does however depend on $r$, which crucially remains fixed. The
$\frac{n}{a_{n}}$ factor arises from the length of $I_{j}$.
Suppose now that for a given $\ell $,
\begin{equation}
\left\vert \chi _{\left[ \tau _{j},a_{n}\right] }-R_{j}\right\vert
(x) \frac{W(x)}{W\left( \tau _{j}\right) }\leq C\left( 1+%
\frac{n}{a_{n}}\left\vert x-\tau _{j}\right\vert \right) ^{-\ell },\qquad \mbox{ }%
x\in \mathbb{R},
\end{equation}%
where $C$ is independent of $f,n,j,x$. Then substituting this and
(4.33)
into (4.32) gives%
\begin{eqnarray*}
\left\vert S\left[ f\right] -P\left[ f\right] \right\vert (x) W(x)\
\leq C\left( \frac{n}{a_{n}}\right) ^{1/p}\sum_{j=1}^{2n-1}\left( 1+\frac{n%
}{a_{n}}\left\vert x-\tau _{j}\right\vert \right) ^{r-\ell
}W\left( \tau _{j}\right) \left( \Omega _{j-1}+\Omega _{j}\right)
.
\end{eqnarray*}%

From here on, we need to proceed a little differently for $p\leq 1$ and $p>1$%
. Let us suppose $p>1$. By H\"{o}lder's inequality, with
$q=\frac{p}{p-1},$
\begin{eqnarray}
\left\vert S\left[ f\right] -P\left[ f\right] \right\vert
^{p}(x)W^{p}(x)
&\leq &C\frac{n}{a_{n}}\left[ \sum_{j=1}^{2n-1}\left( 1+\frac{n}{a_{n}}%
\left\vert x-\tau _{j}\right\vert \right) ^{( r-\ell)
p/2}W^{p}( \tau _{j}) \Omega _{j}^{p}\right] \times  \nonumber \\
&&\times \left[ \sum_{j=1}^{2n-1}\left( 1+\frac{n}{a_{n}}\vert x-\tau _{j}\vert \right)
^{( r-\ell) q/2}\right] ^{p/q}.
\end{eqnarray}%
We also use the fact that $1+\frac{n}{a_{n}}\left\vert x-\tau
_{j}\right\vert $ is bounded by a constant times $1+\frac{n}{a_{n}}%
\left\vert x-\tau _{j-1}\right\vert $ throughout the real line. Next, if $%
\ell >r,$ the function $u\mapsto \left( 1+\frac{n}{a_{n}}\left\vert
x-u\right\vert \right) ^{\left( r-\ell \right) q/2}$ is increasing in $%
\left( -\infty ,x\right) $ and decreasing in $\left( x,\infty
\right) $, so
we can bound the second sum by an integral:%
\begin{eqnarray*}
\sum_{j=1}^{2n-1}\left( 1+\frac{n}{a_{n}}\left\vert x-\tau
_{j}\right\vert \right) ^{\left( r-\ell \right) q/2}
\;\leq\; 2\frac{n}{a_{n}}\int_{-\infty }^{\infty }\left( 1+\frac{n}{a_{n}}%
\left\vert x-u\right\vert \right) ^{\left( r-\ell \right) q/2}du+1
\leq C,
\end{eqnarray*}%
with $C\neq C\left( n,x\right) $, provided only $\left( r-\ell
\right) q/2<-1 $. So all we need is that $\ell $ is large enough.
Now we integrate
(4.35):%
\begin{eqnarray}
\left\Vert (S\left[ f\right] -P\left[ f\right] )W\right\Vert
_{L_{p}(\mathbb{R}) }^{p} &\leq
&C\frac{n}{a_{n}}\sum_{j=1}^{2n-1}W^{p}\left( \tau _{j}\right)
\Omega _{j}^{p}\int_{-\infty }^{\infty }\left( 1+\frac{n}{a_{n}}\left\vert x-\tau
_{j}\right\vert \right) ^{\left( r-\ell \right) p/2}dx  \nonumber \\
&\leq &C\sum_{j=1}^{2n-1}W^{p}\left( \tau _{j}\right) \Omega
_{j}^{p},
\end{eqnarray}%
again, provided $\ell $ is so large that $\left( r-\ell \right)
p/2<-1$. Finally,
\begin{eqnarray*}
\sum_{j=1}^{2n-1}W^{p}\left( \tau _{j}\right) \Omega _{j}^{p}
&=&\frac{n}{a_{n}}\int_{0}^{a_{n}/n}\left[
\sum_{j=1}^{2n-1}W^{p}\left( \tau _{j}\right) \int_{I_{j}^{\ast
}}\left\vert \Delta _{s}^{r}f(x)
\right\vert ^{p}dx\right] \mbox{ }ds \\
&\leq &C\sup_{0<h\leq
\frac{a_{n}}{n}}\int_{-a_{n}}^{a_{n}}\left\vert W(x) \Delta
_{h}^{r}f(x)\right\vert ^{p}dx=C\sup_{0<h\leq
\frac{a_{n}}{n}}\left\Vert W\Delta _{h}^{r}f\right\Vert
_{L_{p}\left[ -a_{n},a_{n}\right] }^{p}.
\end{eqnarray*}%
Combining this, (4.31), and (4.36), gives%
\begin{eqnarray*}
\left\Vert (f-P\left[ f\right] )W\right\Vert _{L_{p}\left( \mathbb{R}%
\right) }^{p} \leq C\left\{ \sup_{0<h\leq
\frac{a_{n}}{n}}\left\Vert W\Delta _{h}^{r}f\right\Vert
_{L_{p}\left[ -a_{n},a_{n}\right] }^{p}+\left\Vert fW\right\Vert
_{L_{p}\left( \vert x\vert \geq a_{n}\right)
}^{p}\right\} .
\end{eqnarray*}%
That's it! Reformulating this in terms of the modulus of
continuity is relatively straightforward. Of course we
assumed:
\newline\noindent
\textbf{Step 4: Construction of the }$\left\{ R_{j}\right\}$.
Recall we want $R_{j}$ of degree $\leq Ln$, satisfying
\begin{equation}
\left\vert \chi _{\left[ \tau ,a_{n}\right] }-R_{j}\right\vert
(x) \frac{W(x)}{W(\tau) }\leq C\left( 1+\frac{%
n}{a_{n}}\left\vert x-\tau \right\vert \right) ^{-\ell },\quad
\mbox{ }x\in \mathbb{R},\mbox{ }\tau \in \left[
-a_{n},a_{n}\right] ,
\end{equation}%
with constants independent of $\tau ,n,x$. The problem here is the
$W(\tau) $ in the denominator. It's tiny for $\tau $
close to $a_{n}$, and we want $R_{j}$ to approximate $1$ in
$\left[ \tau ,a_{n}\right] $, and to approximate $0$ elsewhere in
$\mathbb{R}$. One starts with an even entire function
\[
G(x)=\sum_{j=0}^{\infty }g_{2j}x^{2j}
\] %
with all $g_{2j}\geq 0$ such that
\[
C_{1}\leq (GW) (x)\leq C_{2}\mbox{, }\qquad x\in \mathbb{R}%
.
\] %
Such functions were constructed in \cite{Lubinsky1986}. We let
$G_{n}$ denote the $n$th partial sum. One can show that
\[
\left( G_{n}W\right) (x)\leq C_{1},\mbox{ }\qquad
x\in \mathbb{R},
\] %
and
\[
C_{1}\leq \left( G_{n}W\right) (x)\leq C_{2}\mbox{,
}\qquad \vert x\vert \leq C_{1}a_{n}.
\] %
Next, we need needle or peaking polynomials $V_{n,\xi }$, built
from
Chebyshev polynomials, satisfying%
\[
\left\Vert V_{n,\xi }\right\Vert _{L_{\infty }\left[ -1,1\right]
}=V_{n,\xi }(\xi) =1;
\] %
\[
\left\vert V_{n,\xi }(t) \right\vert \leq \frac{B\sqrt{%
1-\left\vert \xi \right\vert }}{n\left\vert t-\xi \right\vert }\mbox{, }\qquad t\in %
\left[ -1,1\right] \backslash \left\{ \xi \right\};
\] %
and
\[
V_{n,\xi }(t) \geq \frac{1}{2},\mbox{ }\qquad
\left\vert t-\xi \right\vert \leq C\frac{\sqrt{1-\left\vert \xi
\right\vert }}{n}.
\] %
The constants $B,C$ are independent of $n,\xi ,x$. We define \cite[p.~121]%
{DitzianLubinsky1997}
\[
R_{j}(x)=\frac{\int_{0}^{x}G_{Ln/4}(s)
V_{n,\xi }( s/a_{2Ln}) ^{L}ds}{\int_{0}^{\tau ^{\ast
}}G_{Ln/4}(s) V_{n,\xi }( s/a_{2Ln})
^{L}ds}
\] %
where $L$ and $\tau ^{\ast }$ are appropriately chosen, and $\xi
=\tau /a_{Ln}$. To prove this works, we split the integral into
various pieces, consider several ranges of $x$, and reduce other
ranges to the main range, where $\tau \in \left[ S,a_{n}\right] $,
for some fixed $S$. All the details appear in
\cite{DitzianLubinsky1997}.

We note that for Erd\H{o}s weights, or exponential weights on $\left[ -1,1%
\right] $, the technical details are still more complicated ---
the subintervals $I_{j}$ of the $\left[ -a_{n},a_{n}\right] $ are
no longer of
equal length. See \cite{DamelinLubinsky1998}, \cite{Lubinsky1997A}, \cite%
{Lubinsky1997B}.

\sect{Weights Close to exp$\left( -\vert x\vert \right)
$}

The weight $\exp( -\vert x\vert) $ sits
``on'' the boundary of the class of weights admitting a positive
solution to
Bernstein's problem. That boundary is fuzzy, but if you recall that $%
W_{\alpha }(x)=\exp( -\vert x\vert
^{\alpha }) $ admits a positive solution iff $\alpha \geq
1$, this makes sense. So it is not surprising that results like
Jackson's theorem, tend to take a different form. Freud was
interested in this boundary case right throughout his research on
weighted polynomial approximation. In 1978, Freud, Giroux and
Rahman \cite[p.~360]{Freudetal1978} proved that
\begin{eqnarray*}
E_{n}[f;W_{1}]_{1} =\inf_{\deg (P) \leq n}\norm{(f-P) W_{1}}_{L_{1}(\mathbb{R}) } \leq C\left[
\omega \left( f,\frac{1}{\log n}\right) +\int_{\left\vert
x\right\vert \geq \sqrt{n}}\left\vert fW_{1}\right\vert (x)dx%
\right] ,
\end{eqnarray*}%
where
\[
\omega \left( f,\varepsilon \right) =\sup_{\vert h\vert
\leq \varepsilon }\int_{-\infty }^{\infty }\left\vert ( fW_{1}) (x+h) -( fW_{1})
(x) \right\vert dx+\varepsilon \int_{-\infty }^{\infty
}\left\vert fW_{1}\right\vert .
\] %
Here $C$ is independent of $f$ and $n$, and $\sqrt{n}$ could be replaced by $%
n^{1-\delta }$ for any fixed $\delta \in \left( 0,1\right) $. Compare the $%
\frac{1}{\log n}$ inside the modulus to the $n^{-1+1/\alpha }$ we
obtained for $W_{\alpha },\alpha >1$. This suggests that
\begin{equation}
\lim_{\alpha \rightarrow 1+}n^{-1+1/\alpha }=\frac{1}{\log n},
\end{equation}%
at least in the sense of Jackson rates! Ditzian, the author, Nevai
and Totik
\cite{Ditzianetal1987} later extended this result to a characterization in $%
L_{1}$. The technique used by Freud, Giroux and Rahman was essentially an $%
L_{1}$ technique, using the relation between one-sided weighted
approximation, Gauss quadratures, and Christoffel functions --- as
we discussed in Section 4.1.

Only recently has it been possible to establish the analogous results in $%
L_{p},$ $p>1$ \cite{Lubinsky2005}. The author modified the spline
method discussed in Section 4.4 above. The challenge is that
however you tweak the construction of $\left\{ R_{j}\right\} $
there, it does not work. So a new
idea was needed. As the peaking polynomials used there do not work for $%
W_{1} $, they were replaced by the reproducing kernel for
orthogonal polynomials for $W_{1}^{2}$, and in the proofs, the
author needed bounds for these orthogonal polynomials, implied by
recent work of Kriecherbauer and McLaughlin
\cite{KriecherbauerMcLaughlin1999}.

If we examine the modulus used in (3.13) for $W_{\alpha },\alpha
>1$, we see that the interval $[-h^{\frac{1}{1-\alpha
}},h^{\frac{1}{1-\alpha }}]$ is no longer meaningful for $\alpha
=1$. It turns out to be replaced by $\left[
-\exp \left( \frac{1-\varepsilon }{h}\right) ,\exp \left( \frac{%
1-\varepsilon }{h}\right) \right] $, for some fixed $\varepsilon
\in \left( 0,1\right) $. The modulus becomes%
\begin{eqnarray}
\omega _{r,p}(f,W_{1},t) &=&\sup_{0<h\leq t}\norm{W_{1}\left( \Delta
_{h}^{r}f\right)}_{L_{p}\left[ -\exp \left( \frac{1-\varepsilon }{t%
}\right) ,\exp \left( \frac{1-\varepsilon }{t}\right) \right] }  \nonumber \\[10pt]
&+&\inf_{\deg (P) \leq r-1} \norm{ (f-P)
W_{1}}_{L_{p}\left( \mathbb{R}\backslash \lbrack -\exp \left( \frac{%
1-\varepsilon }{t}\right) +1,\exp \left(
\frac{1-\varepsilon}{t}\right) -1]\right) }.
\end{eqnarray}%
The author proved \cite{Lubinsky2005}:

\begin{thh}
\textsl{For }$0<p\leq \infty ,$\textsl{\ and }$n\geq C_{3},$%
\begin{equation}
E_{n}\left[ f;W_{1}\right] _{p}\leq C_{1}\omega
_{r,p}(f,W_{1},\frac{1}{\log \left( C_{2}n\right) }).
\end{equation}%
\textsl{Here }$C_{1},C_{2}$\textsl{\ are independent of }$f$\textsl{\ and }$%
n $\textsl{.}\end{thh}

While this may be a technical achievement, it is scarcely
surprising, given
that Freud, Giroux and Rahman already had the rate $O\left( \frac{1}{\log n}%
\right) $ in the $L_{1}$ case. What about a Jackson--Favard
inequality? The \textquotedblleft limit\textquotedblright\ (5.1)
suggests that an analogue of (3.10) should have the form
\[
E_{n}\left[ f;W_{1}\right] _{p}\leq \frac{C}{\log
n}\norm{f^{\prime }W_{1}}_{L_{p}(\mathbb{R}) }.
\] %
Remarkably enough this is false, and there is no Jackson--Favard
inequality for $W_{1}$, not even if we replace $\frac{1}{\log n}$
by a sequence
decreasing arbitrarily slowly to $0$. More generally, we answered in \cite%
{Lubinsky2005B} the question: which weights admit a Jackson type
theorem, of the form (3.10), with $\left\{ a_{n}/n\right\}
_{n=1}^{\infty }$ replaced by some sequence $\left\{ \eta
_{n}\right\} _{n=1}^{\infty }$ with limit $0?$ We proved
\cite{Lubinsky2005B}:

\begin{thh}  \textsl{Let }$W:\mathbb{R}\rightarrow \left( 0,\infty \right) $\textsl{\ be continuous. The following are
equivalent:}
\begin{description}
\item[(a)] {There exists a sequence }$\left\{ \eta _{n}\right\} _{n=1}^{\infty }$%
\textsl{\ of positive numbers with limit }$0$\textsl{\ and with
the following property. For each }$1\leq p\leq \infty $\textsl{,
and for all absolutely continuous }$f$\textsl{\ with
}$\norm{f^{\prime }W} _{L_{p}(\mathbb{R}) }$\textsl{\
finite, we have }
\begin{equation}
\inf_{\deg (P) \leq n}\norm{(f-P) W}
_{L_{p}(\mathbb{R}) }\leq \eta _{n}\norm{f^{\prime
}W}_{L_{p}(\mathbb{R}) },\qquad n\geq 1.
\end{equation}%
\item[(b)] \textsl{ Both}
\begin{equation}
\lim_{x\rightarrow \infty }W(x)\int_{0}^{x}W^{-1}=0
\end{equation}%
\textsl{and}
\begin{equation}
\lim_{x\rightarrow \infty }W(x)^{-1}\int_{x}^{\infty
}W=0,
\end{equation}%
\textsl{with analogous limits as }$x\rightarrow -\infty
$\textsl{.}\end{description}\end{thh}

Two fairly direct corollaries of this are:

\begin{cor}  \textsl{Let }$W:\mathbb{R}\rightarrow
\left( 0,\infty \right) $\textsl{\ be continuous, with
}$W=e^{-Q}$\textsl{, where }$Q(x)$\textsl{\ is
differentiable for large }$\vert x\vert $\textsl{, and
}
\begin{equation}
\lim_{x\rightarrow \infty }Q^{\prime }(x)=\infty \mbox{ and }%
\lim_{x\rightarrow -\infty }Q^{\prime }(x)=-\infty .
\end{equation}%
\textsl{Then there exists a sequence }$\left\{ \eta _{n}\right\}
_{n=1}^{\infty }$\textsl{\ of positive numbers } \textsl{with
limit }$0$\ \textsl{such that for each }$1\leq p\leq \infty
$\textsl{, and for all absolutely continuous }$f$\textsl{\ with
}$\norm{f^{\prime }W} _{L_{p}(\mathbb{R}) }$\textsl{\
finite, we have $(5.4)$.}\end{cor}

\begin{cor}  \textsl{Let }$W:\mathbb{R}\rightarrow
\left( 0,\infty \right) $\textsl{\ be continuous, with
}$W=e^{-Q}$\textsl{, where }$Q(x)$\textsl{\ is
differentiable for large }$\vert x\vert $\textsl{, and }$%
Q^{\prime }(x)$ \textsl{is bounded for large
}$\left\vert
x\right\vert $\textsl{. Then for both }$p=1$\textsl{\ and }$p=\infty $%
\textsl{, there does not exist a sequence }$\left\{ \eta
_{n}\right\} _{n=1}^{\infty }$\textsl{\ of positive numbers with
limit }$0$\textsl{\
satisfying (5.4) for all absolutely continuous }$f$\textsl{\ with }$%
\norm{f^{\prime }W}_{L_{p}(\mathbb{R}) }$\textsl{\
finite.}\end{cor}

In particular for $W_{1}$, there is no Jackson--Favard inequality,
since both (5.5) and (5.6) are false. Thus there is a real
difference between density of weighted polynomials, and weighted
Jackson--Favard theorems. It is possible to have the former
without the latter.

A key ingredient in the above theorem is an estimate for tails \cite%
{Lubinsky2005B}:

\begin{thh} \textsl{Assume that }$W:\mathbb{R}\rightarrow \left( 0,\infty \right) $%
\textsl{\ is continuous.}
\begin{description}
\item[(a)]\textsl{ Assume }$W$
\textsl{satisfies $(5.5)$ and $(5.6)$, with analogous limits at
}$-\infty $\textsl{. Then there exists a decreasing positive
function }$\eta :[0,\infty )\rightarrow \left( 0,\infty \right)
$\textsl{\ with limit }$0$\textsl{\ at }$\infty $\textsl{\ such
that for }$1\leq p\leq \infty $\textsl{\ and }$\lambda \geq 0$,
\begin{equation}
\norm{fW}_{L_{p}\left( \mathbb{R}\backslash \left[ -\lambda
,\lambda \right] \right) }\leq \eta \left( \lambda \right)
\norm{f^{\prime }W}_{L_{p}(\mathbb{R}) }
\end{equation}%
\textsl{for all absolutely continuous functions
}$f:\mathbb{R}\rightarrow \mathbb{R}$\textsl{\ for which }$f(0)=0$ \textsl{and the right-hand side is finite.}
\item[(b)]{ Conversely assume that $(5.8)$ holds for }$p=1$\textsl{\ and for }$%
p=\infty $\textsl{, for large enough }$\lambda $\textsl{. Then the
limits $(5.5)$ and $(5.6)$ in Theorem $1.1$ are valid, with analogous
limits at }$-\infty $.\end{description}\end{thh}

The above results deal with $L_{p}$ for all $1\leq p\leq \infty $.
What
happens if we focus on a single $L_{p}$ space? We recently proved \cite%
{Lubinsky2006}:

\begin{thh}  \textsl{Let }$W:\mathbb{R}\rightarrow \left( 0,\infty \right) $\textsl{\ be
continuous and let }$1\leq p\leq \infty $ \textsl{and} $\frac{1}{p}+\frac{1}{%
q}=1$\textsl{. The following are equivalent:}
\begin{description}
\item[(a)] {There exists a sequence }$\left\{ \eta _{n}\right\} _{n=1}^{\infty }$%
\textsl{\ of positive numbers with limit }$0$\textsl{\ and with
the
following property. For all absolutely continuous }$f$\textsl{\ with }$%
\norm{f^{\prime }W}_{L_{p}(\mathbb{R}) }$\textsl{\
finite, we have }
\begin{equation}
\inf_{\deg (P) \leq n}\norm{(f-P) W}
_{L_{p}(\mathbb{R}) }\leq \eta _{n}
\norm{f^{\prime}W}_{L_{p}(\mathbb{R}) },\qquad n\geq
1.
\end{equation}%
\item[(b)]
\begin{equation}
\lim_{x\rightarrow \infty }\left\Vert W^{-1}\right\Vert _{L_{q}\left[ 0,x%
\right] }\left\Vert W\right\Vert _{L_{p}[x,\infty )}=0,
\end{equation}%
\textsl{with analogous limits as }$x\rightarrow -\infty .$
\end{description}\end{thh}

As a consequence one can construct weights that admit a Jackson theorem in $%
L_{p}$ but not in $L_{r}$ for any $1\leq p,r\leq \infty $ with $p\neq r$:%

\begin{thh}  \textsl{Let }$1\leq p,r\leq \infty $\textsl{\
with} $p\neq r$. \textsl{There exists }$W:\mathbb{R}\rightarrow
\left( 0,\infty \right) $\textsl{\ such
that }%
\[
\frac{1}{1+x^{2}}\leq W(x)/\exp( -x^{2})
\leq 1+x^{2},\qquad x\in \mathbb{R}\mbox{,}
\] %
\textsl{and }$W$ \textsl{admits an }$L_{r}$\textsl{\ Jackson
theorem , but not an }$L_{p}$\textsl{\ Jackson theorem. That is,
there exists }$\left\{
\eta _{n}\right\} _{n=1}^{\infty }$\textsl{\ with limit }$0$\textsl{\ at }$%
\infty $\textsl{\ satisfying $(5.9)$ in the }$L_{r}$ \textsl{norm,
but there
does not exist such a sequence satisfying $(5.9)$ in the }$L_{p}$ \textsl{norm.}\end{thh}

This is a highly unusual occurrence in weighted approximation ---
in fact the first occurrence of this phenomenon known to this
author. Density of polynomials, and the degree of approximation is
almost invariably the same for any $L_{p}$ space (suitably
weighted of course). Recall Koosis' remark quoted after Theorem
1.6 \cite[pp.~210--211]{Koosis1988}.

\sect{Restricted Range Inequalities}

We have already seen the role played by \textit{restricted range inequalities%
} (often called \dword{infinite-finite range inequalities}) in
weighted Jackson theorems. Paul Nevai rated Freud's discovery of
these as one of his most important contributions to weighted
approximation theory. One is reminded of their import when one
recalls that Dzrbasjan did not have these, and could only prove
estimates in a fixed finite interval. This should not detract from
admiration for the generality of Dzrbasjan's results, and the
sophistication of his ideas, which are still being used.

Let's recap a little. Recall that Freud and Nevai proved that
there are constants $C_{1}$ and $C_{2}$ such that for all
polynomials $P_{n}$ of degree at most $n$,
\begin{equation}
\norm{P_{n}W_{\alpha }}_{L_{p}(\mathbb{R}) }\leq
C_{1}\norm{P_{n}W_{\alpha }}_{L_{p}\left[ -C_{1}n^{1/\alpha
},C_{1}n^{1/\alpha }\right] }.
\end{equation}%
The constants $C_{1}$ and $C_{2}$ can be taken independent of
$n,P_{n}$ and
even the $L_{p}$ parameter $p\in \lbrack 1,\infty ]$. Outside the interval $%
\left[ -C_{1}n^{1/\alpha },C_{1}n^{1/\alpha }\right] $,
$P_{n}W_{\alpha }$ decays quickly to zero. For more general
$W=\exp (-Q) $, one replaces $n^{1/\alpha }$, with
Freud's number $q_{n}$, the positive root of
the equation%
\[
n=q_{n}Q^{\prime }\left( q_{n}\right).
\] %
If $xQ^{\prime }(x)$ is positive, continuous, and
strictly increasing in $\left( 0,\infty \right) $, with limit $0$
at $0$, and $\infty $ at $\infty $, then $q_{n}$ exists and is
unique.

\begin{thh}  \textsl{Assume that }$xQ^{\prime }(x) $\textsl{\ is positive, continuous, and strictly
increasing in }$\left( 0,\infty \right) $\textsl{,
with limit }$0$\textsl{\ at }$0$\textsl{, and }$\infty $\textsl{\ at }$%
\infty$.\textsl{\ Then for }$%
n\geq 1$\textsl{\ and polynomials }$P$\textsl{\ of degree} $\leq n,$%
\begin{equation}
\left\Vert PW\right\Vert _{L_{\infty }(\mathbb{R})
}=\left\Vert PW\right\Vert _{L_{\infty }\left[
-4q_{2n},4q_{2n}\right] }.
\end{equation}\end{thh}

\proof. This type of proof was used many times by
Freud and Nevai; see also the monograph \cite[p.~66]{Mhaskar1996}.
Recall that if $T_{n}$ is the classical Chebyshev polynomial, then
for $\vert x\vert \geq 1$ and
polynomials $P$ of degree $\leq n,$%
\begin{equation}
\left\vert P(x)\right\vert \leq T_{n}\left( \vert x\vert \right) \left\Vert P\right\Vert
_{L_{\infty }\left[ -1,1\right] }\leq \left( 2\left\vert
x\right\vert \right) ^{n}\left\Vert P\right\Vert _{L_{\infty
}\left[ -1,1\right] }.
\end{equation}%
Scaling this to $\left[ -q_{2n},q_{2n}\right] $, and using the
fact that $Q$
is increasing in $\left( 0,\infty \right) $ gives%
\begin{eqnarray*}
\left\vert P(x)\right\vert \leq \left( \frac{2\vert x\vert }{q_{2n}}\right) ^{n}\left\Vert
P\right\Vert _{L_{\infty }\left[ -q_{2n},q_{2n}\right] } \leq
\left( \frac{2\vert x\vert }{q_{2n}}\right)
^{n}W^{-1}\left( q_{2n}\right) \left\Vert PW\right\Vert
_{L_{\infty }\left[ -q_{2n},q_{2n}\right] }.
\end{eqnarray*}%
So%
\[
\left\vert PW\right\vert (x)\leq
2^{n}\frac{\left\vert
x\right\vert ^{n}W(x)}{q_{2n}^{n}W\left( q_{2n}\right) }%
\left\Vert PW\right\Vert _{L_{\infty }\left[ -q_{2n},q_{2n}\right]
}.
\] %
Now if $x\geq 4q_{2n},$%
\begin{eqnarray*}
\log \frac{x^{n}W(x)}{q_{2n}^{n}W\left( q_{2n}\right)
} =\int_{q_{2n}}^{x}\frac{n-uQ^{\prime }(u) }{u}du
\leq -n\int_{q_{2n}}^{4q_{2n}}\frac{du}{u}=-n\log 4,
\end{eqnarray*}%
since $u\geq q_{2n}\Rightarrow uQ^{\prime }(u) \geq
2n$.
Substituting in above, gives for $x\geq 4q_{2n},$%
\[
\vert PW\vert (x)\leq 2^{-n}\left\Vert
PW\right\Vert _{L_{\infty }\left[ -q_{2n},q_{2n}\right] }.
\] %
This certainly implies (6.2).  \endproof

The proof clearly shows the exponential decay of $PW$ outside the interval $%
\left[ -q_{2n},q_{2n}\right] $, which depends only on $n$. In $L_{p}$, any $%
p>0$, a typical analogue is%
\[
\left\Vert PW\right\Vert _{L_{p}(\mathbb{R}) }\leq
(1+e^{-Cn})\left\Vert PW\right\Vert _{L_{p}\left[
-Bq_{2n},Bq_{2n}\right] }
\] %
where $C>0$, and $B$ is large enough.

These inequalities are sufficient for weighted Jackson and
Bernstein theorems, but not for some of the deeper results such as
Bernstein inequalities with endpoint effects, or convergence of
orthogonal expansions, or Lagrange interpolation, let alone
asymptotics of orthogonal polynomials. With such questions in
mind, Mhaskar and Saff asked in the early 1980's: \textit{where
does the sup or }$L_{p}$\textit{\ norm of a weighted polynomial
really live?} Clearly it is inside the intervals $\left[
-Bq_{2n},B_{q_{2n}}\right] $, with appropriate $B$, but what is the best $B$%
? In seminal papers \cite{MhaskarSaff1984}, \cite{MhaskarSaff1985}, \cite%
{MhaskarSaff1987}, they used potential theory to obtain the
answer. In
another seminal development, E.~A.~Rakhmanov a little earlier \cite%
{Rakhmanov1984} developed the same potential theory in order to
investigate asymptotics of orthogonal polynomials.

One of the best ways to understand their work is to recall
Bernstein's bound
for growth of polynomials in the complex plane. For $n\geq 1$, polynomials $%
P $ of degree $\leq n$, and $z\notin \left[ -1,1\right] $,
\[
\left\vert P(z) \right\vert \leq \left\vert z+\sqrt{z^{2}-1}%
\right\vert ^{n}\left\Vert P\right\Vert _{L_{\infty }\left[
-1,1\right] }.
\] %
Thus once we have a bound on $P$ on an interval, we can estimate
its growth outside. Of course, it is related to the inequality
involving $T_{n}$ that we used above, but that works only for
real $x$.

Now let us look for a weighted analogue for even weights $W=\exp
(-Q) $, where, say, $x\mapsto xQ^{\prime }(x) $ is increasing, positive,
and continuous in $\left(0,\infty \right) $. Suppose that $a>0$, and we have a function
$G=G_{n,a}$ with the following properties:
\begin{description}
\item[(I)] $G$ is analytic and non-vanishing in $\mathbb{C}\backslash [ -a,a%
] $;
\item[(II)] $G(z) z^{n}$ has a finite
limit as $\vert z\vert \rightarrow \infty ;$
\item[(III)] $\vert G(z) \vert \rightarrow W(x) $ as $z\rightarrow x\in \left( -a,a\right) $.
\end{description}

Then given a polynomial $P$ of degree $\leq n$, the function $PG$
is analytic in $\mathbb{C}\backslash [ -a,a] $,
including at $\infty
$, where it has a finite limit. Moreover, for $x\in \left( -a,a\right) ,$%
\[
\lim_{z\rightarrow x}\vert PG\vert (z)
=\vert PW\vert (x)\leq \Vert
PW\Vert _{L_{\infty } [ -a,a] }.
\]
By the maximum modulus principle, we then obtain%
\[
\vert PG\vert (z) \leq \Vert
PW\Vert _{L_{\infty }[ -a,a] }\mbox{, }\qquad
z\notin [ -a,a] .
\] %
In particular, for real $x$ with $\vert x\vert >a,$
\[
\vert PW\vert (x)\leq \Vert
PW\Vert _{L_{\infty }[ -a,a] }W(x)
/G_{n,a}(x).
\]

For a given $n$, Mhaskar and Saff found the smallest $a=a_{n}$ for
which
\begin{equation}
W(x)/G_{n,a}(x)<1\quad \mbox{ for all
}\vert x\vert >a.
\end{equation}%
Clearly, we then have%
\[
\vert PW\vert (x)<\Vert PW\Vert
_{L_{\infty }\left[ -a_{n},a_{n}\right] },\qquad \left\vert
x\right\vert >a_{n},
\] %
and the \textit{Mhaskar--Saff identity}
\[
\Vert PW\Vert _{L_{\infty }(\mathbb{R})
}=\Vert PW\Vert _{L_{\infty }\left[ -a_{n},a_{n}\right]
}.
\] %
Recall that $a_{n}$ is called the \textit{Mhaskar--Rakhmanov--Saff
}number,
and is the positive root of the equation%
\begin{equation}
n=\frac{2}{\pi }\int_{0}^{1}\frac{a_{n}tQ^{\prime }\left( a_{n}t\right) }{%
\sqrt{1-t^{2}}}dt.
\end{equation}%
For the weight $W_{\alpha }(x)=\exp( -\vert x\vert ^{\alpha }) ,\alpha >0$,
\[
a_{n}=\left\{ 2^{\alpha -2}\frac{\Gamma \left( \alpha /2\right)
^{2}}{\Gamma \left( \alpha \right) }\right\} ^{1/\alpha
}n^{1/\alpha }.
\]

How does this arise and how does one find the function $G_{n,a}?$
We solve
the integral equation%
\begin{equation}
\int_{-1}^{1}\log \left\vert x-t\right\vert \mu _{n,a}(t) dt=%
\frac{Q(ax) }{n}+c_{n,a}\mbox{, }\qquad x\in \left[
-1,1\right],
\end{equation}%
for the function $\mu _{n,a}$, and some constant $c_{n,a}$,
subject to the
condition%
\[
\int_{-a}^{a}\mu _{n,a}=1.
\] %
By (formal) differentiation, we see that this is really equivalent
to solving a singular integral equation
\[
\int_{-1}^{1}\frac{\mu _{n,a}(t)
}{x-t}dt=\frac{aQ^{\prime }(ax) }{n}\mbox{, }\qquad
x\in \left( -a,a\right) \mbox{.}
\] %
The integral on the left is a Cauchy principal value integral.
There is a
well developed theory of such equations \cite{Henrici1986}, \cite%
{SaffTotik1997}. One can show that \cite[p.~37]{LubinskySaff1988},
\cite[p.~124 ff.]{Mhaskar1996}
\begin{equation}
\mu _{n,a}(x)=\frac{2}{\pi ^{2}}\int_{0}^{1}\frac{\sqrt{%
1-x^{2}}}{\sqrt{1-s^{2}}}\frac{asQ^{\prime }(as)
-axQ^{\prime
}(ax) }{n\left( s^{2}-x^{2}\right) }ds+\frac{B_{n,a}}{\pi \sqrt{%
1-x^{2}}},
\end{equation}%
where%
\begin{equation}
B_{n,a}=1-\frac{2}{n\pi }\int_{0}^{1}\frac{atQ^{\prime }(at) }{%
\sqrt{1-t^{2}}}dt.
\end{equation}%
Then the function $G_{n,a}$ is given by
\[
G_{n,a}(z) =\exp \left( -n\int_{-1}^{1}\log \left( \frac{z}{a}%
-t\right) \mu _{n,a}(t) dt+nc_{n,a}\right)
\] %
so that for $x\in \left[ -a,a\right] ,$
\[
\left\vert G_{n,a}(x)\right\vert =\exp \left( -n\int_{-1}^{1}\log \left\vert \frac{x}{a}-t\right\vert \mu
_{n,a}(t) dt+nc_{n,a}\right) =W(x),
\] %
by (6.6). So at least we have verified (III) above. It turns out
that whenever $B_{n,a}>0$, then $W/G_{n,a}>1$ in some right
neighborhood of $1$. Thus we look for the smallest $a$ for which
$B_{n,a}=0$. We see from (6.8) that this requires
\[
1-\frac{2}{n\pi }\int_{0}^{1}\frac{atQ^{\prime }(at) }{\sqrt{%
1-t^{2}}}dt=0\quad \Longrightarrow \quad a=a_{n}.
\]

This scratches the surface of an extensive theory. The monograph
of Saff and Totik \cite{SaffTotik1997} contains a detailed and
deep development. A more elementary treatment is provided in
\cite{Mhaskar1996}, though that will be more than sufficient for
those interested only in the topics of this survey. For general
exponential weights, restricted range inequalities in quite
precise form are investigated in \cite{LevinLubinsky2001}. Here is
a typical
result \cite[Theorem 1.8, p.~15]{LevinLubinsky2001}, \cite[Thm.~6.2.4, p.~142%
]{Mhaskar1996}:

\begin{thh}  \textsl{Let} $Q:\mathbb{R}\rightarrow \lbrack
0,\infty )$ \textsl{be even
and convex, with limit }$\infty $\textsl{\ at }$\infty $\textsl{\ and}%
\[
0=Q(0) <Q(x)\mbox{, }\qquad x\neq 0.
\] %
\begin{description}
\item[(a)]\textsl{ For not identically zero polynomials }$P$\textsl{\ of degree }$%
\leq n,$%
\[
\left\Vert PW\right\Vert _{L_{\infty }(\mathbb{R})
}=\left\Vert PW\right\Vert _{L_{\infty }\left[ -a_{n},a_{n}\right]
}
\] %
\textsl{and}%
\[
\left\Vert PW\right\Vert _{L_{\infty }\left( \mathbb{R}\backslash
\lbrack -a_{n},a_{n}]\right) }<\left\Vert PW\right\Vert
_{L_{\infty }\left[ -a_{n},a_{n}\right] }.
\] %
\item[(b)]\textsl{ Let }$0<p<\infty $\textsl{\ and }$P$\textsl{\
be a not identically zero polynomial of degree }$\leq
n-\frac{2}{p}$\textsl{. Then}
\[
\left\Vert PW\right\Vert _{L_{p}(\mathbb{R})
}<2^{1/p}\left\Vert PW\right\Vert _{L_{p}\left[
-a_{n},a_{n}\right] }
\] %
\textsl{and}%
\[
\left\Vert PW\right\Vert _{L_{p}\left( \mathbb{R}\backslash
\lbrack
-a_{n},a_{n}]\right) }<\left\Vert PW\right\Vert _{L_{p}\left[ -a_{n},a_{n}%
\right] }.
\]\end{description}\end{thh}

Quite often, we need a smaller estimate for the tail, and this is
possible provided we omit a little more than the Mhaskar--Saff
interval. Here is what can be achieved \cite[Thm.~1.9,
p.~15]{LevinLubinsky2001} for the class of weights
$\mathcal{F}^{\ast }$, specified in Definition 4.6:

\begin{thh} \textsl{Let} $W\in \mathcal{F}^{\ast }$ a\textsl{nd }$0<p\leq \infty $%
\textsl{. For }$n\geq 1,\kappa \in (0,1]$\textsl{\ and polynomials }$P$%
\textsl{\ of degree }$\leq n,$%
\[
\left\Vert PW\right\Vert _{L_{p}\left( \mathbb{R}\backslash
\lbrack -a_{n}\left( 1+\kappa \right) ,a_{n}\left( 1+\kappa
\right) ]\right) }\leq
C_{1}\exp( -nC_{2}\kappa ^{3/2}) \left\Vert PW\right\Vert _{L_{p}%
\left[ -a_{n},a_{n}\right] }.
\] %
\textsl{Here }$C_{1}$\textsl{\ and }$C_{2}$\textsl{\ are independent of }$%
n,P,\kappa $.\end{thh}

More generally, it was proved there

\begin{thh}  \textsl{Let }$W=e^{-Q}$\textsl{, where
}$Q:\mathbb{R}\rightarrow \lbrack
0,\infty )$\textsl{\ is even and convex, with }$Q$\textsl{\ having limit }$%
\infty $\textsl{\ at }$\infty ,$\textsl{\ and }%
\[
T(t) =\frac{tQ^{\prime }(t) }{Q(t) }
\] %
\textsl{is quasi-increasing in the sense of $(3.27)$ of Definition
$3.11$.
Assume furthermore that for some }$\Lambda >1$\textsl{, }$T\geq \Lambda $%
\textsl{\ in }$\left( 0,\infty \right) $\textsl{, while }%
\begin{equation}
T(y) \sim T\left( y\left[ 1-\frac{1}{T(y)
}\right] \right) ,\mbox{ }\qquad y\in \left( 0,\infty \right) .
\end{equation}%
\textsl{Then for }$n\geq 1,\kappa \in (0,\frac{1}{T\left( a_{n}\right) }]$%
\textsl{\ and polynomials }$P$\textsl{\ of degree }$\leq n,$%
\begin{equation}
\left\Vert PW\right\Vert _{L_{p}\left( \mathbb{R}\backslash
\lbrack -a_{n}\left( 1+\kappa \right) ,a_{n}\left( 1+\kappa
\right) ]\right) }\leq C_{1}\exp( -C_{2}nT( a_{n}) \kappa ^{3/2}) \left\Vert PW\right\Vert
_{L_{p}\left[ -a_{n},a_{n}\right] }.
\end{equation}%
\textsl{Here }$C_{1}$\textsl{\ and }$C_{2}$\textsl{\ are independent of }$%
n,P,\kappa $.\end{thh}

Finally, we can go back a little inside the Mhaskar--Saff interval
if we
allow a cruder estimate on the tail \cite[Theorem 4.2(a), p.~96]%
{LevinLubinsky2001}:

\begin{thh} \textsl{Assume in addition to the hypotheses of Theorem $6.4$ that }$%
Q^{\prime \prime }$\textsl{\ exists and is positive in }$\left( 0,\infty
\right) $\textsl{, while for some }$C>0$\textsl{, and large enough }$x,$%
\[
\frac{Q^{\prime \prime }(x)}{Q^{\prime }(x)}\leq C%
\frac{Q^{\prime }(x)}{Q(x)}.
\] %
\textsl{Let }$0<p\leq \infty $\textsl{\ and} $\lambda >0$\textsl{.
Then
there exists }$n_{0}$\textsl{\ and }$C$\textsl{\ such that for }$n\geq n_{0}$%
\textsl{\ and polynomials }$P$\textsl{\ of degree }$\leq n,$%
\begin{equation}
\left\Vert PW\right\Vert _{L_{p}(\mathbb{R}) }\leq
C\left\Vert PW\right\Vert _{\left[ -a_{n}\left( 1-\lambda \eta
_{n}\right) ,a_{n}\left( 1-\lambda \eta _{n}\right) \right] },
\end{equation}%
\textsl{where}%
\begin{equation}
\eta _{n}=\left( nT\left( a_{n}\right) \right) ^{-2/3}.
\end{equation}\end{thh}

See \cite{JungSakai2006}, \cite{KasugaSakai2003},
\cite{LevinLubinsky2001},
\cite{Mhaskar1991}, \cite{Mhaskar1996}, \cite{Nevai1973}, \cite%
{SaffTotik1997} for further discussion of restricted range
inequalities.

\sect{Markov--Bernstein Inequalities}

Markov--Bernstein inequalities have already been mentioned above.
They are the main ingredient of converse theorems of
approximation, but also enter in many other contexts. For the
unweighted case on $\left[ -1,1\right] $, the classical Markov
inequality asserts that
\begin{equation}
\left\Vert P^{\prime }\right\Vert _{L_{\infty }\left[ -1,1\right]
}\leq n^{2}\left\Vert P\right\Vert _{L_{\infty }\left[ -1,1\right]
},
\end{equation}%
for $n\geq 1$ and all polynomials $P$ of degree $\leq n$. The
Bernstein
inequality improves this as long as we stay away from the endpoints:%
\begin{equation}
\left\vert P^{\prime }(x)\right\vert \leq \frac{n}{\sqrt{1-x^{2}%
}}\left\Vert P\right\Vert _{L_{\infty }\left[ -1,1\right] },\mbox{
}\qquad x\in \left( -1,1\right) ,
\end{equation}%
for $n\geq 1$ and all polynomials $P$ of degree $\leq n$.

The weights $W_{\alpha }(x)=\exp( -\vert
x\vert ^{\alpha }) ,\alpha >0$, already provide a lot
of insight into the weighted case on the real line. The analogue
of the Markov
inequality is%
\begin{equation}
\left\Vert P_{n}^{\prime }W_{\alpha }\right\Vert _{L_{p}\left( \mathbb{R}%
\right) }\leq C\left\Vert P_{n}W_{\alpha }\right\Vert _{L_{p}\left( \mathbb{R%
}\right) }\left\{
\begin{array}{ll}
n^{1-1/\alpha }, & \alpha >1 \\
\log (n+1) , & \alpha =1 \\
1, & \alpha <1%
\end{array}%
\right. .
\end{equation}%
This is valid for all $n\geq 1$ and polynomials $P_{n}$ of degree
$\leq n$. The constant $C$ depends on $p\in (0,\infty ]$ and
$\alpha $, but not on $n$ or $P$. It's no accident that the
factors arising are similar to those in the Jackson rates, and for
$\alpha <1$, the factor is bounded independent of
$n$. For $\alpha \geq 2$, these inequalities were proved by Freud; for $%
1<\alpha <2$, by Eli Levin and the author; and for $\alpha \leq
1$, by Paul Nevai and Vili Totik.

In order to generalize this for arbitrary $Q$, we should try to
cast this estimate in a unified form. For $\alpha >1$, we see that
\[
n^{1-1/\alpha }=C_{1}\frac{n}{a_{n}}=C_{2}Q^{\prime }\left( a_{n}\right) ,
\] %
where $Q(x)=\vert x\vert ^{\alpha }$, and recall, $%
a_{n}=C_{1}n^{1/\alpha }$. As we shall see, $Q^{\prime }\left( a_{n}\right) $ is the correct factor whenever $W=\exp(
-Q) $ and $Q$ is even and grows faster than $\left\vert
x\right\vert ^{\alpha }$ for some $\alpha
>1$. For Freud weights, where $Q$ is of polynomial growth, but still grows
faster than $\vert x\vert ^{\alpha }$, for some $\alpha
>1$, we can use $\frac{n}{a_{n}}$. However, for all polynomial
rates of growth,
including $\vert x\vert ^{\alpha },\alpha \leq 1$, we can use%
\[
\int_{1}^{Cn}\frac{ds}{Q^{\left[ -1\right] }(s) },
\] %
where if $Q$ is strictly increasing and continuous on $[0,\infty )$, $Q^{%
\left[ -1\right] }$ denotes its inverse. The most general $L_{\infty \mbox{ }%
}$result for $Q$ of polynomial growth, with this factor, is due to
Kro\'{o} and Szabados \cite[p.~48]{KrooSzabados1995}:

\begin{thh}  \textsl{Let }$W=\exp (-Q)
$\textsl{, where }$Q$\textsl{\ is even, continuous, increasing in
}$\left( 0,\infty \right) $\textsl{, and
twice differentiable for large enough }$x$\textsl{, with }%
\begin{equation}
\liminf_{x\rightarrow \infty }\frac{xQ^{\prime }(x)
}{Q(x)}>0,
\end{equation}%
\textsl{and}
\begin{equation}
\limsup_{x\rightarrow \infty }\frac{xQ^{\prime \prime }(x)}{%
Q^{\prime }(x)}<\infty .
\end{equation}%
\textsl{Then for }$n\geq 1$\textsl{\ and all polynomials
}$P$\textsl{\ of
degree }$\leq n,$\textsl{\ }%
\[
\left\Vert P^{\prime }W\right\Vert _{L_{\infty }(\mathbb{R})
}\leq \int_{1}^{Cn}\frac{ds}{Q^{\left[ -1\right] }(s) }%
\left\Vert PW\right\Vert _{L_{\infty }(\mathbb{R}) }.
\] %
\textsl{Here }$C\neq C\left( n,P\right) $.\end{thh}

The author and Eli Levin used the slightly more restrictive condition that $%
W\in \mathcal{F}^{\ast }$, so (4.18) holds. This implies, like
(7.4) and (7.5), that $Q$ grows faster than some positive power of
$x$, but slower than some other positive power of $x$. Strangely
enough, there does not seem to be a published form of this result
in $L_{p}$. Indeed, extra difficulties arise in $L_{p}$ --- we
shall discuss techniques for proving Markov--Bernstein
inequalities in the next section. However, there are results that
apply to $L_{p}$ \textit{separately} for the case where $Q$ grows
faster, or slower, than $\vert x\vert $. First, we
record the former case \cite[p.~231]{LevinLubinsky1994}:

\begin{thh} \textsl{Let }$W\in \mathcal{F}^{\ast }$. \textsl{Let }$1\leq p<\infty $%
\textsl{. Then for }$n\geq 1$\textsl{\ and all polynomials
}$P$\textsl{\ of
degree }$\leq n,$\textsl{\ }%
\[
\left\Vert P^{\prime }W\right\Vert _{L_{p}(\mathbb{R}) }\leq C%
\frac{n}{a_{n}}\left\Vert PW\right\Vert _{L_{p}(\mathbb{R}) }.
\] %
\textsl{Here }$C\neq C\left( n,P\right) $.\end{thh}

This was later extended to the case $0<p<1$ \cite[Corollary 1.16, p.~21]%
{LevinLubinsky2001}. For weights growing slower than $\left\vert
x\right\vert $, Nevai and Totik \cite[Thm.~2,
p.~122]{NevaiTotik1986} used a beautiful method to prove:

\begin{thh} \textsl{Let }$W=\exp (-Q)
$\textsl{, where }$Q$\textsl{\ is even, increasing, and concave in
}$\left( 0,\infty \right) $, \textsl{\ with}
\[
\int_{0}^{\infty }\frac{Q(x)}{1+x^{2}}dx<\infty .
\] %
\textsl{Then for }$n\geq 1$\textsl{\ and all polynomials
}$P$\textsl{\ of
degree }$\leq n,$\textsl{\ }%
\[
\left\Vert P^{\prime }W\right\Vert _{L_{\infty }(\mathbb{R}) }\leq C\left\Vert PW\right\Vert _{L_{\infty
}(\mathbb{R}) }.
\] %
\textsl{Here }$C\neq C\left( n,P\right) $.\textsl{\ If in addition, }%
\[
\lim_{x\rightarrow \infty }\frac{2Q(x)-Q(2x) }{\log x}=\infty ,
\] %
\textsl{then, given }$p>0$\textsl{, we have for }$n\geq
1$\textsl{\ and all
polynomials }$P$\textsl{\ of degree }$\leq n,$\textsl{\ }%
\[
\left\Vert P^{\prime }W\right\Vert _{L_{p}(\mathbb{R})
}\leq C\left\Vert PW\right\Vert _{L_{p}(\mathbb{R}) }.
\] %
\textsl{Here }$C\neq C\left( n,P\right) $.\end{thh}

The above are all Markov inequalities, analogues of (7.3). There
are also Bernstein inequalities, which reflect the opposite
feature of the finite interval case. The growth in $n$
\textit{decreases} towards the end of the Mhaskar--Rakhmanov--Saff
interval. To state these, we need the function
\[
\phi _{n}(x):=\max \left\{ 1-\frac{\vert x\vert }{%
a_{n}},n^{-2/3}\right\} ,\mbox{ }\qquad x\in \mathbb{R}\mbox{.}
\] %
Eli Levin and the author proved \cite[Thm.~1.1, p.~231]{LevinLubinsky1994}:%

\begin{thh}  \textsl{Let }$W\in \mathcal{F}^{\ast
}$\textsl{\ and let }$1\leq p\leq
\infty $. \textsl{Then for }$n\geq 1$\textsl{\ and all polynomials }$P$%
\textsl{\ of degree }$\leq n,$\textsl{\ }%
\[
\Vert (PW) ^{\prime }\phi _{n}^{-1/2}\Vert
_{L_{p}(\mathbb{R}) }\leq C\frac{n}{a_{n}}\left\Vert
PW\right\Vert _{L_{p}(\mathbb{R}) }.
\] %
\textsl{Here }$C\neq C\left( n,P\right) $.\end{thh}

For $p=\infty $, this was proved in \cite{LevinLubinsky1990}, and
for $1\leq p<\infty $, in \cite{LevinLubinsky1994}. If $p=\infty
$, we see that this
has the consequence%
\[
\left\vert (PW) ^{\prime }\left( a_{n}\right) \right\vert \leq C%
\frac{n^{2/3}}{a_{n}}\left\Vert PW\right\Vert _{L_{\infty }\left( \mathbb{R}%
\right) },
\] %
when $P$ has degree $\leq n$. There is in \cite[Thm.~1.2, p.~233]%
{LevinLubinsky1994} an analogue of this for weights satisfying
(4.18) in Definition 4.6 with $A>0$ only (such as exp$\left( -\vert x\vert ^{\alpha }\right) ,\alpha \leq 1$). One
fixes $\eta \in \left( 0,1\right) $, and proves that
\[
\Vert (PW) ^{\prime }\phi _{n}^{-1/2}\Vert
_{L_{p}\left( \vert x\vert \geq \eta a_{n}\right) }\leq C\frac{n}{%
a_{n}}\left\Vert PW\right\Vert _{L_{p}(\mathbb{R}) }.
\] %
Note that in the Bernstein inequalities, it is essential that we estimate $%
(PW) ^{\prime }$, and not $P^{\prime }W$. There is no
improvement for the latter, at least in general, near $\pm a_{n}$.

Analogues of Markov and Bernstein inequalities have also been
obtained for Erd\H{o}s weights, where $Q$ is of faster than
polynomial growth, as well as for exponential weights on $\left( -1,1\right) $. They are also available for non-even exponential
weights on a possibly asymmetric finite or infinite interval
\cite{LevinLubinsky2001}. For simplicity, we quote only the even
case. First, we define a suitable class of weights, which includes
both the
Freud and Erd\H{o}s weights in the real line, and exponential weights on $%
\left( -1,1\right) $:

\begin{deff}
\textsl{Let }$I=\left( -d,d\right) $\textsl{\ where }$0<d\leq \infty $%
\textsl{. Let }$Q:I\rightarrow \lbrack 0,\infty )$\textsl{\ be an
even function with the following properties:}
\begin{description}
\item[(a)] $Q^{\prime }$\textsl{\ is continuous and positive in }$I$\textsl{\ and }%
$Q(0) =0;$
\item[(b)] $Q^{\prime \prime
}$\textsl{\ exists and is positive in }$\left( 0,d\right)
;$
\item[(c)]
\[
\lim_{t\rightarrow d-}Q(t) =\infty ;
\] %
\item[(d)]\textsl{ The function }%
\[
T(t) =\frac{tQ^{\prime }(t) }{Q(t) }
\] %
\textsl{is quasi-increasing in }$\left( 0,d\right) $\textsl{, in
the sense
that }%
\[
0\leq x<y<d\quad \Longrightarrow \quad T(x)\leq
CT(y) .
\] %
\textsl{Moreover, for some }$\Lambda >1$\textsl{, }%
\[
T(t) \geq \Lambda >1\mbox{, }\qquad t\in \left( 0,d\right) .
\] %
\item[(e)]\textsl{ There exists }$C_{1}>0$\textsl{\ such that }%
\[
\frac{\left\vert Q^{\prime \prime }(x)\right\vert
}{Q^{\prime }(x)}\leq C_{1}\frac{Q^{\prime }(x) }{Q(x)},\qquad x\in \left( 0,d\right).
\] %
\end{description}
\textsl{Then we write }$W\in \mathcal{F}_{even}\left( C^{2}\right) $\textsl{%
. If also, there exists }$c\in \left( 0,d\right) $\textsl{\ such
that}
\[
\frac{\left\vert Q^{\prime \prime }(x)\right\vert
}{Q^{\prime }(x)}\geq C_{2}\frac{Q^{\prime }(x) }{Q(x)},\qquad x\in \left( c,d\right) ,
\] %
\textsl{then we write }$W\in \mathcal{F}_{even}\left( C^{2}+\right) $.\end{deff}

Examples of weights satisfying this are $W_{\alpha },\alpha >1$, as well as $%
W=\exp (-Q) $, where
\begin{equation}
Q(x)=\exp _{\ell }\left( \vert x\vert
^{\alpha }\right) -\exp _{\ell }(0) ,\mbox{ }\qquad
x\in \left( -\infty ,\infty \right),
\end{equation}%
with $\ell \geq 1,\alpha >1,$ and%
\[
\exp _{\ell }:=\exp(( \cdots\exp( {})
))
\] %
denoting the $\ell $th iterated exponential. We set $\exp
_{0}(x) :=x$. Others include%
\begin{equation}
Q(x)=\exp _{\ell }\left( \left( 1-x^{2}\right)
^{-\alpha }\right) -\exp _{\ell }(1) ,\mbox{ }\qquad
x\in \left( -1,1\right) ,
\end{equation}%
where now $\ell \geq 0$ and $\alpha >0$.

\begin{thh}
\textsl{Let }$0<p\leq \infty $\textsl{\ and }$W=e^{-Q}\in \mathcal{F}%
_{even}\left( C^{2}\right) $.\textsl{\ For }$n\geq 1$\textsl{, let }%
\begin{equation}
\varphi _{n}(x)=\frac{a_{n}}{n}\frac{\left\vert 1-\frac{%
\vert x\vert }{a_{2n}}\right\vert }{\sqrt{\left\vert 1-\frac{%
\vert x\vert }{a_{n}}\right\vert +\left( nT\left( a_{n}\right) \right) ^{-2/3}}}.
\end{equation}%
\textsl{Then for }$n\geq 1$\textsl{\ and polynomials
}$P_{n}$\textsl{\ of
degree }$\leq n,$%
\[
\left\Vert (PW) ^{\prime }\varphi _{n}\right\Vert
_{L_{p}(\mathbb{R}) }\leq C\left\Vert PW\right\Vert _{L_{p}\left( \mathbb{R}%
\right) }.
\] %
\textsl{Moreover,}
\[
\left\Vert P^{\prime }W\right\Vert _{L_{p}(\mathbb{R}) }\leq C%
\frac{n}{a_{n}}T( a_{n}) ^{1/2}\left\Vert PW\right\Vert
_{L_{p}(\mathbb{R}) }.
\]\end{thh}

This result is a special case of Theorem 1.15 and Corollary 1.16
in \cite[p.~21]{LevinLubinsky2001}. There the restrictions on $Q$
are weaker, but the definition of the classes is more implicit, so
we restrict ourselves to the
smallest even case considered there. The main feature is the extra factor $%
T( a_{n}) ^{1/2}$, which really is there --- weighted
extremal polynomials attain it, see \cite{LevinLubinsky2001}.

As examples, $Q$ of (7.6) has Markov factor \cite[p.~30]{LevinLubinsky2001}%
\[
\frac{n}{a_{n}}T( a_{n}) ^{1/2}\sim \frac{n}{\left( \log _{\ell }n\right) ^{1/\alpha }}\left( \prod\limits_{j=1}^{\ell
}\log _{j}n\right) ^{1/2}.
\] %
Here
\[
\log _{j}:=\log \left( \log \left( \cdots\log \left( {}\right)
\right) \right)
\] %
denotes the $j$th iterated logarithm. For $Q$ of (7.7) and $\ell
=0$, the Markov factor is \cite[p.~32]{LevinLubinsky2001}
\[
\frac{n}{a_{n}}T( a_{n}) ^{1/2}\sim n^{\frac{2\alpha
+2}{2\alpha +1}}
\] %
while if $\ell \geq 1$, it is \cite[p.~34]{LevinLubinsky2001}
\[
\frac{n}{a_{n}}T( a_{n}) ^{1/2}\sim n\left( \log _{\ell
}n\right) ^{\frac{1}{2}\left( 1+\frac{1}{\alpha }\right) }\left( \prod\limits_{j=1}^{\ell -1}\log _{j}n\right) ^{1/2}.
\]

\sect{Methods to Prove Markov--Bernstein Inequalities}

In this section, we outline some of the methods that have been
used to prove Markov--Bernstein inequalities for exponential
weights. An extensive treatment of Markov--Bernstein inequalities
for both weighted and unweighted cases is given in
\cite{BorweinErdelyi1995}, \cite{Milovanovicetal1994}. We begin
with

\subsection{Freud's Method via de la Vall\'ee Poussin Means}

Recall that\ in Lemma 4.5 we proved
\[
\left\Vert V_{n}\left[ f\right] W\right\Vert _{L_{\infty }\left( \mathbb{R}%
\right) }\leq C\left\Vert fW\right\Vert _{L_{\infty }\left( \mathbb{R}%
\right) },
\] %
where $C\neq C\left( n,f\right) $. Thereafter, at (4.17), we
showed how
duality can be used to prove%
\[
\left\Vert V_{n}\left[ f\right] W\right\Vert _{L_{1}(\mathbb{R}) }\leq C\left\Vert fW\right\Vert _{L_{1}\left(
\mathbb{R}\right) },
\] %
where $C\neq C\left( n,f\right) $. In much the same way, one can
prove that when $W\in \mathcal{F}^{\ast },$
\begin{equation}
\Vert V_{n}^{\prime }\left[ f\right] W\Vert _{L_{\infty
}(\mathbb{R}) }\leq C\frac{n}{a_{n}}\left\Vert
fW\right\Vert _{L_{\infty }(\mathbb{R}) }.
\end{equation}%
The main new technical ingredient required is the estimate%
\begin{equation}
1\Bigl/\;\sum_{k=0}^{n-1}p_{k}^{\prime }(\xi) ^{2}\geq C\left( \frac{%
a_{n}}{n}\right) ^{3}W^{2}(\xi) ,\qquad \xi \in
\mathbb{R}.
\end{equation}%
This is a cousin of the Christoffel function estimate (4.13) and
can be proved using much the same ideas. Freud did this in
\cite{Freud1977}. A proof is also given in
\cite[p.~64]{Mhaskar1996}. Now if $P$ is a polynomial of degree
$\leq n$, then
\[
V_{n}\left[ P\right] =P.
\] %
So (8.1) immediately gives the $L_{\infty }$ Markov--Bernstein inequality%
\[
\Vert P^{\prime }W\Vert _{L_{\infty }(\mathbb{R})
}\leq C\frac{n}{a_{n}}\left\Vert PW\right\Vert _{L_{\infty }\left( \mathbb{R}%
\right) }.
\] %
To extend this to $L_{1}$, we use duality on $V_{n}^{\prime }$:%
\begin{equation}
\left\Vert V_{n}^{\prime }\left[ f\right] W\right\Vert _{L_{1}\left( \mathbb{%
R}\right) }=\sup \int_{-\infty }^{\infty }V_{n}^{\prime }\left[
f\right] gW^{2},
\end{equation}%
where the $\sup $ is taken over all measurable functions $g$ with $%
\left\Vert gW\right\Vert _{L_{\infty }(\mathbb{R})
}\leq 1$. Next, we use the fact that if $P$ is a polynomial of
degree $\leq 2n$, then
\[
\int_{-\infty }^{\infty }\left( g-V_{4n}\left[ g\right] \right)
PW^{2}=0.
\] %
Indeed, we can prove this by considering the special case
$P=p_{j}$, $0\leq j\leq 2n$, and recalling that the first $2n$
Fourier series coefficients (with respect to $\left\{
p_{j}\right\} $) of $V_{4n}\left[ g\right]$ are the same as those
of $g$. The above considerations and an integration by parts give
\begin{eqnarray}
\int_{-\infty }^{\infty }V_{n}^{\prime }\left[ f\right] gW^{2}
&=&\int_{-\infty }^{\infty }V_{n}^{\prime }\left[ f\right] V_{4n}\left[ g%
\right] W^{2}  \nonumber \\[10pt]
&=&-\int_{-\infty }^{\infty }V_{n}\left[ f\right] V_{4n}^{\prime }\left[ g%
\right] W^{2}+2\int_{-\infty }^{\infty }V_{n}\left[ f\right] V_{4n}\left[ g%
\right] Q^{\prime }W^{2}.
\end{eqnarray}%
Here%
\begin{eqnarray}
\left\vert \int_{-\infty }^{\infty }V_{n}\left[ f\right] V_{4n}^{\prime }%
\left[ g\right] W^{2}\right\vert &\leq& \left\Vert V_{4n}^{\prime
}\left[ g\right] W\right\Vert _{L_{\infty
}(\mathbb{R}) } \int_{-\infty }^{\infty }\left\vert V_{n}\left[ f%
\right]\right\vert W  \nonumber \\[10pt]
&\leq &C\frac{n}{a_{n}}\left\Vert gW\right\Vert _{L_{\infty }\left( \mathbb{R%
}\right) }\left\Vert fW\right\Vert _{L_{1}(\mathbb{R})
}
\leq C\frac{n}{a_{n}}\left\Vert fW\right\Vert _{L_{1}\left( \mathbb{R}%
\right) }.
\end{eqnarray}%
In the second last line, we used our $L_{1}$ bound (4.17) for
$V_{n}$, and the $L_{\infty }$ bound (8.1) for $V_{n}^{\prime }$,
and in the last line we used our bound on the sup norm of $gW$. To
bound the second term in (8.4), one needs

\begin{lem} \textsl{Let }$h$\textsl{\ be absolutely
continuous with }$h(0) =0 $\textsl{. Then}
\begin{equation}
\left\Vert Q^{\prime }hW\right\Vert _{L_{\infty }(\mathbb{R})
}\leq C\left\Vert h^{\prime }W\right\Vert _{L_{\infty }\left( \mathbb{R}%
\right) }.
\end{equation}\end{lem}

\proof.
Observe that if $x>0,$%
\begin{eqnarray}
\left\vert Q^{\prime }(x)h(x)W(x) \right\vert &=&\left\vert Q^{\prime }(x)
W(x)
\int_{0}^{x}h^{\prime }(t) dt\right\vert  \nonumber \\
&\leq &\left\Vert h^{\prime }W\right\Vert _{L_{\infty }\left( \mathbb{R}%
\right) }\left\vert Q^{\prime }(x)W(x)
\int_{0}^{x}W^{-1}(t) dt\right\vert .
\end{eqnarray}%
We now assume that $Q^{\prime }$ is increasing, and
\[
\lim_{t\rightarrow \infty }\frac{Q^{\prime \prime }(t) }{%
Q^{\prime }(t) ^{2}}=0.
\] %
This latter condition is true for regularly behaved Freud and
Erd\H{o}s weights such as $\mathcal{F}^{\ast },\mathcal{E}$.
Choose $A>0$ such that
\[
\frac{Q^{\prime \prime }(t) }{Q^{\prime }(t) ^{2}}%
\leq \frac{1}{2},\mbox{ }\quad t\geq A.
\] %
If $x\geq A,$ an integration by parts gives%
\begin{eqnarray*}
\int_{A}^{x}W^{-1}(t) dt &=&\int_{A}^{x}Q^{\prime
}(t) ^{-1}Q^{\prime }(t) W^{-1}(t) dt \\[10pt]
&=&Q^{\prime }(t) ^{-1}W^{-1}(t)
\Bigl|_{t=A}^{t=x}\;+\int_{A}^{x}\frac{Q^{\prime \prime }(t) }{%
Q^{\prime }(t) ^{2}}W^{-1}(t) dt \\[10pt]
&\leq &Q^{\prime }(x)^{-1}W^{-1}(x)+\frac{1}{2}%
\int_{A}^{x}W^{-1}(t) dt.
\end{eqnarray*}%
Then%
\[
\int_{A}^{x}W^{-1}(t) \mbox{ }dt\leq 2Q^{\prime
}(x)^{-1}W^{-1}(x)
\] %
and%
\[
x\geq A\;\Rightarrow\; Q^{\prime }(x)W(x)
\int_{A}^{x}W^{-1}(t) \mbox{ }dt\leq 2.
\] %
As $Q^{\prime }W$ is bounded in $\left( 0,\infty \right) $, we obtain%
\[
x\geq 0\;\Rightarrow\; Q^{\prime }(x)W(x)
\int_{0}^{x}W^{-1}(t) \mbox{ }dt\leq C.
\] %
A similar bound holds for $x<0$, and then we obtain (8.6) from (8.7).
 \endproof

Next, we apply (8.6) to the second term in the right-hand side of
(8.4), with $h=V_{4n}\left[ g\right] -V_{4n}\left[ g\right] (0) $. We
obtain%
\begin{eqnarray*}
\left\vert V_{4n}\left[ g\right] Q^{\prime }W\right\vert (x)
&\leq &\left\vert V_{4n}\left[ g\right] (x)-V_{4n}\left[ g%
\right] (0) \right\vert \left\vert Q^{\prime }(x) \right\vert W(x)+V_{4n}\left[ g\right]
(0)
\left\vert Q^{\prime }(x)\right\vert W(x)\\[10pt]
&\leq &C\Vert V_{4n}^{\prime }\left[ g\right] W\Vert
_{L_{\infty }(\mathbb{R}) }+C\left\vert V_{4n}\left[
g\right] (0) \right\vert
\leq C\frac{n}{a_{n}}\left\Vert gW\right\Vert _{L_{\infty }\left( \mathbb{R%
}\right) }\leq C\frac{n}{a_{n}},
\end{eqnarray*}%
so%
\begin{eqnarray}
\left\vert \int_{-\infty }^{\infty }V_{n}\left[ f\right] V_{4n}\left[ g%
\right] Q^{\prime }W^{2}\right\vert
\leq C\frac{n}{a_{n}}\left\vert \int_{-\infty }^{\infty }V_{n}\left[ f%
\right] W\right\vert
\leq C\frac{n}{a_{n}}\left\Vert fW\right\Vert _{L_{1}\left( \mathbb{R}%
\right) }.
\end{eqnarray}%
Combining (8.4), (8.5), and (8.8) gives%
\[
\Vert V_{n}^{\prime }\left[ f\right] W\Vert _{L_{1}\left( \mathbb{%
R}\right) }\leq C\frac{n}{a_{n}}\left\Vert fW\right\Vert
_{L_{1}(\mathbb{R}) }.
\] %
Thus, recalling (8.1), we have that for both $p=1$ and $p=\infty
$,
\[
\Vert V_{n}^{\prime }\left[ f\right] W\Vert _{L_{p}( \mathbb{%
R}) }\leq C\frac{n}{a_{n}}\Vert fW\Vert
_{L_{p}(\mathbb{R}) }.
\] %
By interpolation, we obtain%
\[
\Vert V_{n}^{\prime }\left[ f\right] W\Vert _{L_{p}( \mathbb{%
R}) }\leq C\frac{n}{a_{n}}\left\Vert fW\right\Vert
_{L_{p}(\mathbb{R}) },
\] %
for all $1\leq p\leq \infty $, where $C\neq C\left( n,P\right) $.
Applying
this with $f=P$, a polynomial of degree $\leq n$, gives%
\[
\Vert P^{\prime }W\Vert _{L_{p}(\mathbb{R}) }\leq C%
\frac{n}{a_{n}}\Vert PW\Vert _{L_{p}(\mathbb{R}) }.
\] %
This method is elegant, and yields more than just
Markov--Bernstein inequalities. As we have seen, we also obtain a
proof of a Jackson--Favard type inequality, bounds on Ces\`aro
means of orthogonal expansions, and other useful information. Our
main technical ingredients were bounds on Christoffel functions
and their derivative analogue (8.2), and some technical estimates
involving $Q^{\prime }.$

\subsection{Replacing the Weight by a Polynomial}

This method is very simple, but applies only to a limited class of
weights.
Suppose that for some fixed $K>0,$ we have polynomials $S_{n}$ of degree $%
\leq Kn$, with
\begin{equation}
C_{1}\leq S_{n}/W\leq C_{2}\mbox{\ \ in\ \ }[
-2a_{n},2a_{n}],
\end{equation}%
and
\begin{equation}
\left\vert S_{n}^{\prime }\right\vert /W\leq C_{3}\frac{n}{a_{n}}\mbox{\ \ in\ \ }%
[ -a_{n},a_{n}] \mbox{.}
\end{equation}%
These polynomials enable us to reduce weighted Bernstein
inequalities to classical unweighted Bernstein inequalities. If
$P$ is a polynomial of
degree $\leq n$, then our restricted range inequalities give%
\begin{eqnarray*}
\left\Vert P^{\prime }W\right\Vert _{L_{p}(\mathbb{R})
} &\leq &C\left\Vert P^{\prime }W\right\Vert _{L_{p}\left[
-a_{n},a_{n}\right] } \leq CC_{1}^{-1}\left\Vert P^{\prime
}S_{n}\right\Vert _{L_{p}\left[
-a_{n},a_{n}\right] }\\[10pt]
&\leq& CC_{1}^{-1}\left[ \left\Vert \left( PS_{n}\right) ^{\prime
}\right\Vert _{L_{p}\left[ -a_{n},a_{n}\right] }+\left\Vert
PS_{n}^{\prime
}\right\Vert _{L_{p}\left[ -a_{n},a_{n}\right] }\right] \\[10pt]
&\leq &CC_{1}^{-1}\left[ \frac{2n}{a_{n}}\left\Vert PS_{n}\right\Vert _{L_{p}%
\left[ -2a_{n},2a_{n}\right] }+\frac{n}{a_{n}}\left\Vert
PW\right\Vert _{L_{p}\left[ -a_{n},a_{n}\right] }\right] ,
\end{eqnarray*}%
by our hypotheses on $S_{n}^{\prime }$, and the classical
Bernstein
inequality, scaled from $\left[ -1,1\right] $ to $\left[ -2a_{n},2a_{n}%
\right] $. Using (8.9) again, we obtain%
\[
\Vert P^{\prime }W\Vert _{L_{p}(\mathbb{R})
}\leq C\left\Vert PW\right\Vert _{L_{p}(\mathbb{R}) }.
\] %
This works in any $L_{p}$, $0<p\leq \infty $. If $p=\infty $, we
can weaken the requirements on $S_{n}$, which can be made
different for each $x$. We only need the upper bound on
$S_{n}^{\prime }$ at a given $x$, and the lower bound on $S_{n}$
at that $x$ (but we still need the upper bound on $S_{n}$
throughout $\left[ -2a_{n},2a_{n}\right] $).

Freud and Nevai used this for weights like $W_{2m}(x)
=\exp( -x^{2m}) $, where $m\geq 1$ is a positive
integer. The partial
sums of these entire weights can be used for $\left\{ S_{n}\right\} $. For $%
W_{\alpha },\alpha >1$, Eli Levin and the author used canonical
products of Weierstrass primary factors, such as
\[
\prod\limits_{n=1}^{\infty }E\left( -x/n^{1/\alpha };\ell \right)
,
\] %
to generate these polynomials \cite{LevinLubinsky1987A}, \cite%
{LevinLubinsky1987B}. This method can provide quick easy proofs in
special cases, which would be useful for teaching a course on
weighted approximation.

\subsection{Nevai--Totik's Method}

This involves fast decreasing polynomials, a topic initiated by
Kamen Ivanov, Paul Nevai and Vili Totik, and works well for slowly
decreasing weights such as $W_{\alpha },\alpha \leq 1$. Let us
suppose that for $n\geq
1 $, and some $K>0$, we have polynomials $S_{n}^{\#}$ that satisfy%
\begin{equation}
S_{n}^{\#}(0) =1, \qquad S_{n}^{\#\prime }(0) =0,
\end{equation}%
and
\begin{equation}
\left\vert S_{n}^{\#}(x)\right\vert \leq Ke^{-Q(x) },\qquad x\in \left[ -a_{n},a_{n}\right] .
\end{equation}%
Typically, $S_{n}^{\#}$ has a unique maximum in $\left[ -1,1\right]
$ at $0$,
and decreases rapidly in $\left( 0,1\right) $. Nevai and Totik \cite%
{NevaiTotik1986} used a construction of Marchenko to generate such
polynomials. One starts with an entire function of the form%
\[
B(z) =\prod\limits_{k\geq 1}\frac{\sin ^{2}\left( \frac{\pi }{2}%
\sqrt{1+\left( \frac{z}{t_{k}}\right) ^{2}}\right) }{1+\left( \frac{z}{t_{k}}%
\right) ^{2}}.
\] %
The product may be finite or infinite, and the $\left\{
t_{k}\right\} $ are positive numbers with
\[
T=\sum_{k}\frac{1}{t_{k}}<\infty .
\] %
Assuming that $Q$ is even, increasing on $\left( 0,\infty \right)
$, and
\begin{equation}
\int_{0}^{\infty }\frac{Q(x)}{1+x^{2}}dx<\infty ,
\end{equation}%
one can choose $\left\{ t_{k}\right\} $ such that $T\leq 1/\pi $
and
\[
\left\vert B(x)\right\vert \leq K\exp( -Q(x)) ,\qquad \mbox{ }x\in \mathbb{R}.
\] %
In this case $B$ is an entire function of exponential type $\leq
1$. Assuming (8.13), Nevai and Totik used the partial sums of $B$
to construct polynomials $P_n$ of degree $\leq n$, with a local
maximum at $0$, and
\begin{equation}
P_{n}(0) =1\mbox{ and }\left\vert P_{n}(x)
\right\vert \leq K\exp( -Q(nx)) \mbox{,
}\qquad x\in \left[ -1,1\right] .
\end{equation}%
Then
\[
S_{n}^{\#}(x)=P_{\left[ a_{n}\right] }\left( \frac{x}{\left[ a_{n}\right] }\right)
\] %
satisfies (8.11) and (8.12), and is of degree $\leq \left[
a_{n}\right] $. It is easy to derive the Markov inequality at $0$:
\begin{eqnarray*}
\left\vert (P^{\prime }W)(0) \right\vert &=&\left\vert
P^{\prime }(0) \right\vert =\left\vert( PS_{n}^{\#}) ^{\prime }(0) \right\vert \leq
\frac{n+\floor{ a_{n}} }{a_{n}}\Vert
PS_{n}^{\#}\Vert _{L_{\infty }[ -a_{n},a_{n}] },
\end{eqnarray*}%
by the usual unweighted Bernstein inequality. We continue this as
\begin{equation}
\left\vert (P^{\prime }W)(0) \right\vert \leq
C\left\Vert PW\right\Vert _{L_{\infty }\left[ -a_{n},a_{n}\right]
},
\end{equation}%
using the property (8.12), and the fact that the convergence (8.13) implies $%
n=O( a_{n}) $, which we shall not prove. To extend this
to general $x\geq 0$, we use the evenness and concavity of $Q$,
which implies that for $x,y\in \mathbb{R},$
\[
W(x)W(y) \leq W(x+y).
\] %
Now apply (8.15) to the polynomial $R(y) =P(x+y)
W(x)$, for fixed $x$:%
\begin{eqnarray*}
\left\vert ( P^{\prime }W) (x)\right\vert
&=&\left\vert R^{\prime }(0) W(0)
\right\vert \leq C\left\Vert RW\right\Vert _{L_{\infty }(\mathbb{R}) } =C\sup_{y\in \mathbb{R}}\left\vert P(
x+y) W(x)
W(y) \right\vert \\
&\leq &C\sup_{y\in \mathbb{R}}\left\vert P(x+y)
W(x+y) \right\vert =C\left\Vert PW\right\Vert
_{L_{\infty }(\mathbb{R}) }.
\end{eqnarray*}%
This method can be modified to give the correct Markov inequality
for $W_{1}$. Nevai and Totik were the first to do so.

\subsection{Dzrbasjan/Kro\'{o}--Szabados' Method}

Dzrbasjan was apparently the first researcher to investigate the
degree of
approximation for general exponential weights, in his 1955 paper \cite%
{Dzrbasjan1955}. Although his approximation estimates worked only
on a finite interval, he nevertheless came up with a great many
ideas. Kro\'{o} and Szabados subsequently used Dzrbasjan's method,
and restricted range inequalities to establish Theorem 7.1.

We start with Cauchy's integral formula for derivatives (or, if
you prefer,
Cauchy's estimates):%
\begin{eqnarray*}
\left\vert P^{\prime }(x)\right\vert =\left\vert \frac{1}{%
2\pi i}\int_{\left\vert t-x\right\vert =\varepsilon }\frac{P(t) }{(t-x) ^{2}}dt\right\vert \leq
\frac{1}{\varepsilon }\sup \left\{ \left\vert P(t)
\right\vert :\left\vert t-x\right\vert =\varepsilon \right\} .
\end{eqnarray*}%
The number $1/\varepsilon $ is invariably chosen as the size of
the Markov--Bernstein factor. To estimate $P(t) $ in
terms of the values of $PW$ on the real line, we write $t=u+iv$,
and use an inequality that often arises in the theory of functions
analytic in the upper half-plane. We already used one form of this in the proof of Lemma 2.2.%
\begin{eqnarray}
\log \left\vert P(u+iv) \right\vert
&\leq &\frac{\vert v\vert }{\pi }\int_{-\infty }^{\infty }\frac{%
\log \left\vert P(s) \right\vert }{(s-u) ^{2}+v^{2}}%
ds  \nonumber
\leq \frac{\vert v\vert }{\pi }\int_{-\infty }^{\infty }\frac{%
\log M+Q(s) }{(s-u) ^{2}+v^{2}}ds  \nonumber \\[10pt]
&=&\log M+\frac{\vert v\vert }{\pi }\int_{-\infty }^{\infty }%
\frac{Q(s) }{(s-u) ^{2}+v^{2}}ds,
\end{eqnarray}%
where
\[
M:=\left\Vert PW\right\Vert _{L_{\infty }(\mathbb{R})
}.
\] %
Of course, we need to assume the integral converges, which is very
restrictive. We shall assume more precisely that
\[
\int_{0}^{\infty }\frac{Q(t) }{1+t^{2}}dt<\infty
\mbox{,}
\] %
and that for $s\geq 1$,%
\[
Q(s) \sim sQ^{\prime }(s) \sim
s^{2}Q^{\prime \prime }(s) .
\] %
Then the above considerations show that
\begin{eqnarray}
\log \left\vert P^{\prime }W\right\vert (x)\leq \log \frac{1}{%
\varepsilon }+\log M+  \sup_v
\frac{\vert v\vert }{\pi }\int_{-\infty }^{\infty }\frac{%
Q(s) -Q(x)}{(s-u)
^{2}+v^{2}}ds.
\end{eqnarray}%
The technical challenge is to estimate the integral in the last
right-hand side. We break it into several pieces. If we assume
$x\geq 1,0\leq u\leq 2x$ and $Q(x)\geq 0$, we see
that
\[
I_{1}:=\frac{\vert v\vert }{\pi }\int_{3u/2}^{\infty }\frac{%
Q(s) -Q(x)}{(s-u) ^{2}+v^{2}}ds\leq 9%
\frac{\varepsilon }{\pi }\int_{3u/2}^{\infty }\frac{Q(s) }{s^{2}}%
ds,
\] %
and%
\[
I_{2}:=\frac{\vert v\vert }{\pi }\int_{-\infty
}^{0}\frac{Q(s) -Q(x)}{(s-u) ^{2}+v^{2}}ds\leq 9\frac{%
\varepsilon }{\pi }\int_{-\infty }^{0}\frac{Q(s)
}{s^{2}}ds.
\] %
(Recall that $\vert v\vert \leq \left\vert
t-x\right\vert \leq \varepsilon $.) In $\left[
0,\frac{u}{2}\right] $, we see that the integrand is non-positive,
as $s\leq u/2\leq x$, so
\[
I_{3}:=\frac{\vert v\vert }{\pi
}\int_{0}^{u/2}\frac{Q(s) -Q(x)}{(s-u) ^{2}+v^{2}}ds\leq 0.
\] %
Finally, we handle the difficult central integral%
\begin{eqnarray*}
I_{4} &:= &\frac{\vert v\vert }{\pi }\int_{u/2}^{3u/2}\frac{%
Q(s) -Q(x)}{(s-u) ^{2}+v^{2}}ds \\[10pt]
&\leq &\frac{\vert v\vert }{\pi
}\int_{u/2}^{3u/2}\frac{Q(s) -Q(u)
}{(s-u) ^{2}+v^{2}}ds+\left\vert
Q(u) -Q(x)\right\vert \\[10pt]
&=:&I_{41}+I_{42}.
\end{eqnarray*}%
Here%
\begin{eqnarray*}
I_{41} =\frac{\vert v\vert }{\pi }\int_{u/2}^{3u/2}\frac{%
Q(s) -Q(u) }{(s-u)
^{2}+v^{2}}ds =\frac{\vert v\vert }{\pi
}\int_{u/2}^{u}\frac{Q(r) -2Q(u) +Q(2u-r) }{(r-u) ^{2}+v^{2}}dr.
\end{eqnarray*}%
For some $\xi $ between $r$ and $u$,%
\[
Q(r) -2Q(u) +Q(2u-r)
=Q^{\prime \prime }(\xi) (r-u) ^{2}
\] %
and one can show that \cite[p.~53]{KrooSzabados1995}
\[
Q^{\prime \prime }(\xi) \leq C\frac{Q(r)
}{r^{2}},
\] %
so%
\[
I_{41}\leq C\frac{\varepsilon }{\pi }\int_{u/2}^{2u}\frac{Q(r) }{%
r^{2}}dr.
\] %
Also,%
\begin{eqnarray*}
I_{42} =\left\vert Q(u) -Q(x)\right\vert
=Q^{\prime }(\xi) \left\vert u-x\right\vert \leq
CQ^{\prime }(x)\varepsilon \leq C\varepsilon
\int_{x}^{2x}\frac{Q(s) }{s^{2}}ds.
\end{eqnarray*}%
Combining the above estimates, we have shown that for $x\geq 1,$%
\[
\log \left\vert P^{\prime }W\right\vert (x)\leq \log \frac{1}{%
\varepsilon }+\log M+C\varepsilon \int_{1}^{\infty }\frac{Q(y) }{%
y^{2}+1}dy.
\] %
We simply choose $\varepsilon =1$. Similar estimates hold over
$(-\infty ,-1] $, and the range $\left[ -1,1\right] $ is easy to
handle. This leads to
the Markov inequality%
\[
\Vert P^{\prime }W\Vert _{L_{\infty }(\mathbb{R}) }\leq C\left\Vert PW\right\Vert _{L_{\infty
}(\mathbb{R}) },
\] %
valid for all $n\geq 1$ and polynomials of degree $\leq n$.

When dealing with weights $W=\exp (-Q) $, where $Q$
grows faster than $\vert x\vert ^{\alpha }$, some
$\alpha >1$, one typically chooses
\[
\varepsilon =\frac{a_{n}}{n},
\] %
uses restricted range inequalities, and needs more work to
estimate the various integrals.

\subsection{Levin--Lubinsky's Method}

Like the previous method, we use Cauchy's integral formula for
derivatives, and then go back to the real line. However, instead
of using (8.16), we use the potential theoretic functions
associated with exponential weights. These give analogues of the
Bernstein--Walsh inequality for growth of polynomials in the
complex plane. For $p<\infty $, going back from the plane to the
real line is quite complicated, and requires Carleson measures. We
shall outline the method for $1\leq p<\infty $.

Fix $x\geq 0$, $\varepsilon >0$, and define an entire function $F_{x}$ by%
\[
F_{x}(z) :=\exp( -Q(x)-Q^{\prime
}(x)(z-x)) .
\] %
Observe that
\[
F_{x}^{(j) }(x)=W^{(j)
}(x)\mbox{ for }j=0,1.
\] %
We have%
\begin{eqnarray}
\left\vert (PW) ^{\prime }(x)\right\vert
&=&\left\vert( PF_{x}) ^{\prime }(x)
\right\vert
=\left\vert \frac{1}{2\pi i}\int_{\left\vert t-x\right\vert =\varepsilon }%
\frac{( PF_{x})(t) }{(t-x) ^{2}}%
dt\right\vert \nonumber\\[10pt]
&\leq &\left( \frac{1}{2\pi }\int_{\left\vert t-x\right\vert
=\varepsilon }\vert PF_{x}\vert (t)
^{p}\left\vert dt\right\vert \right) ^{1/p}\left( \frac{1}{2\pi
}\int_{\vert t-x\vert =\varepsilon }\vert
t-x\vert ^{-2q}\vert dt\vert
\right) ^{1/q} \nonumber\\[10pt]
&=&\frac{1}{\varepsilon }\left( \frac{1}{2\pi }\int_{-\pi }^{\pi
}\vert PF_{x}\vert ^{p}( x+\varepsilon e^{i\theta
}) d\theta \right) ^{1/p}.
\end{eqnarray}%
In the second last line, we used H\"{o}lder's inequality, and $q=\frac{p}{p-1%
}$ there. We shall choose
\[
\varepsilon =\varepsilon \left( n,x\right) =\varphi _{n}(x) ,
\] %
where $\varphi _{n}(x)$ is given by (7.8). This guarantees \cite%
[Lemma 10.6, p.~301]{LevinLubinsky2001} that
\begin{equation}
\left\vert F_{x}(z) \right\vert \leq CW\left( \vert z\vert \right) \quad\mbox{ for }\quad\left\vert
z-x\right\vert \leq \varphi _{n}(x).
\end{equation}%
The proof of this involves a Taylor series expansion of $Q( \func{Re}%
z) $ about $Q(x)$. Now comes the potential
theory bit. There is \cite[Lemma 10.7, p.~303]{LevinLubinsky2001}
a function $G_{n}$ analytic in $\mathbb{C}\backslash \left[
-a_{n},a_{n}\right] $, with boundary values $G_{n}(x)
$ on $\left[ -a_{n},a_{n}\right] $ from the upper half-plane
satisfying
\[
\left\vert G_{n}(x)\right\vert =W(x),\qquad x\in %
\left[ -a_{n},a_{n}\right] .
\] %
Moreover, $z^{n}G_{n}(z) $ has a finite limit at
$\infty $, and uniformly for $x\in \left[ -a_{2n},a_{2n}\right] ,$
\begin{equation}
W( \vert z\vert) \leq C\left\vert
G_{n}(z) \right\vert \mbox{, }\qquad \left\vert
z-x\right\vert \leq \varphi _{n}(x).
\end{equation}%
One representation for $G_{n}$ is
\[
G_{n}(z) =\exp \left( -\int_{-a_{n}}^{a_{n}}\log
(z-t) d\mu _{n}(t) -c_{n}\right) ,
\] %
where $\mu _{n}$ is the equilibrium measure of mass $n$ for the
weight $W$, and $c_{n}$ is an equilibrium constant. Recall that we
used something similar in deriving restricted range inequalities,
in Section 6. Then we obtain from (8.18) to (8.20),
\[
\varphi _{n}^{p}(x)\left\vert (PW)
^{\prime }(x)\right\vert ^{p}\leq C\int_{-\pi }^{\pi
}\vert PG_{n}\vert ^{p}( x+\varphi _{n}(x) e^{i\theta }) d\theta .
\] %
Integrating gives%
\begin{eqnarray}
\int_{-a_{2n}}^{a_{2n}}\varphi _{n}^{p}\vert (PW)
^{\prime }\vert ^{p}  \nonumber &\leq
&C\int_{-a_{2n}}^{a_{2n}}\left[ \int_{-\pi }^{\pi }\vert
PG_{n}\vert ^{p}( x+\varphi _{n}(x)
e^{i\theta
}) d\theta \right] dx  \nonumber \\[10pt]
&=&C\int \left\vert PG_{n}\right\vert ^{p}d\left[ \nu _{n}^{+}+\nu _{n}^{-}%
\right] \mbox{,}
\end{eqnarray}%
where $\nu _{n}^{+}$ is a measure on the upper half-plane, and
$\nu _{n}^{-}$ is a measure on the lower half-plane. For Borel
measurable sets $S$, with
characteristic function $\chi _{S}$, we have%
\begin{eqnarray*}
\nu _{n}^{+}(S) &=&\int_{-a_{2n}}^{a_{2n}}\left[
\int_{0}^{\pi
}\chi _{S}( x+\varphi _{n}(x)e^{i\theta }) d\theta %
\right] dx; \\[10pt]
\nu _{n}^{-}(S) &=&\int_{-a_{2n}}^{a_{2n}}\left[
\int_{-\pi }^{0}\chi _{S}( x+\varphi _{n}(x)
e^{i\theta }) d\theta \right] dx.
\end{eqnarray*}

Next, we use a famous inequality of Carleson. We say a measure
$\nu $ on the upper half-plane is a \dword{Carleson measure} if
there exists $A>0$ such that for all squares $K$ in the upper
half-plane, with base on the real axis, and side $h>0,$
\[
\nu (K) \leq Ah.
\] %
The smallest such number $A$ is called $N\left( \nu \right) $, the
\dword{Carleson norm} of $\nu $. Let $H^{p}$ denote the Hardy space in the
upper half-plane, consisting of all functions analytic there,
whose boundary values on the
real axis lie in $L_{p}(\mathbb{R}) $. For $f\in H^{p},$%
\[
\int \vert f\vert ^{p}d\nu \leq CN\left( \nu \right)
\int_{-\infty }^{\infty }\vert f\vert ^{p}.
\] %
(For $p<1$, there is an analogous statement.)

Applying this to (8.21) gives%
\begin{eqnarray*}
\int_{-a_{2n}}^{a_{2n}}\varphi _{n}^{p}\left\vert (PW)
^{\prime }\right\vert ^{p} \leq C\left( N\left[ \nu
_{n}^{+}\right] +N\left[ \nu _{n}^{-}\right] \right) \int_{-\infty
}^{\infty }\left\vert PG_{n}\right\vert ^{p}.
\end{eqnarray*}%
Here although $\nu _{n}^{-}$ is a measure on the lower half-plane,
it is
obvious what is meant by its Carleson norm. Next, one shows that%
\[
\sup_{n}\left( N\left[ \nu _{n}^{+}\right] +N\left[ \nu
_{n}^{-}\right] \right) <\infty .
\] %
Finally,
\[
\int_{-a_{n}}^{a_{n}}\left\vert PG_{n}\right\vert
^{p}=\int_{-a_{n}}^{a_{n}}\left\vert PW\right\vert ^{p},
\] %
and one can show using Hilbert transform estimates that%
\[
\int_{\vert x\vert \geq a_{n}}\left\vert
PG_{n}\right\vert ^{p}\leq C\int_{-a_{n}}^{a_{n}}\left\vert
PW\right\vert ^{p},
\] %
giving%
\[
\int_{-a_{2n}}^{a_{2n}}\varphi _{n}\left\vert (PW)
^{\prime }\right\vert ^{p}\leq C\int_{-a_{n}}^{a_{n}}\left\vert
PW\right\vert ^{p}.
\] %
The integral over $\mathbb{R}\backslash \left[
-a_{2n},a_{2n}\right] $ may be handled using restricted-range
inequalities.

This method is the deepest of those we presented --- and it is the
only one that gives the Bernstein inequalities in Theorem 7.4 and
7.6. For $p=\infty $,
the proof is easier, as one can avoid Carleson measures \cite%
{LevinLubinsky1990}, \cite{Mhaskar1990}.

\subsection{Sieved Markov--Bernstein Inequalities}

The methods of this section have been extensively developed by
Paul Nevai and others \cite{Lubinskyetal1987},
\cite{LubinskyNevai1987}, and illustrate the power of Jensen's
inequality. On finite intervals, they go back at least to Zygmund
\cite{Zygmund2002}. The idea is to start with the Christoffel type
estimates
\begin{eqnarray}
(PW) ^{2}(\xi) &\leq
&C\frac{n}{a_{n}}\int_{-\infty
}^{\infty }(PW) ^{2},
\end{eqnarray}
and
\begin{eqnarray}
( P^{\prime }W) ^{2}(\xi) &\leq &C\left( \frac{n}{%
a_{n}}\right) ^{3}\int_{-\infty }^{\infty }(PW) ^{2},
\end{eqnarray}%
valid for Freud weights $W$, $n\geq 1$, polynomials $P$ of degree
$\leq n$, and $\xi \in \mathbb{R}$. We used these already in
Section 8.1.

We now extend these to $L_{p}$ type Christoffel function estimates. Let $%
0<p<2$. From (8.22), we derive
\[
\left\Vert PW\right\Vert _{L_{\infty }(\mathbb{R}) }^{2}\leq C%
\frac{n}{a_{n}}\left\Vert PW\right\Vert _{L_{\infty }\left( \mathbb{R}%
\right) }^{2-p}\int_{-\infty }^{\infty }\left\vert PW\right\vert
^{p},
\] %
and hence%
\begin{equation}
\left\Vert PW\right\Vert _{L_{\infty }(\mathbb{R}) }^{p}\leq C%
\frac{n}{a_{n}}\int_{-\infty }^{\infty }\left\vert PW\right\vert
^{p}.
\end{equation}%
Next, (8.23) followed by (8.24) give%
\begin{eqnarray*}
\left\Vert P^{\prime }W\right\Vert _{L_{\infty }(\mathbb{R}) }^{2} &\leq &C\left( \frac{n}{a_{n}}\right)
^{3}\left\Vert PW\right\Vert _{L_{\infty }(\mathbb{R})
}^{2-p}\int_{-\infty }^{\infty
}\left\vert PW\right\vert ^{p} \\[10pt]
&\leq &C\left( \frac{n}{a_{n}}\right) ^{3}\left( \frac{n}{a_{n}}%
\int_{-\infty }^{\infty }\left\vert PW\right\vert ^{p}\right) ^{\frac{2-p}{p}%
}\int_{-\infty }^{\infty }\left\vert PW\right\vert ^{p}.
\end{eqnarray*}%
Rearranging gives%
\begin{equation}
\left\Vert P^{\prime }W\right\Vert _{L_{\infty }(\mathbb{R}) }^{p}\leq C\left( \frac{n}{a_{n}}\right)
^{p+1}\int_{-\infty }^{\infty }\left\vert PW\right\vert ^{p}.
\end{equation}%
Thus far we have (8.24) and (8.25) for $0<p\leq 2$. Our weight is
$W=\exp (-Q) $. We apply these inequalities instead to
the weight $\exp (-rQ) $ for fixed $r>0$. Its
Mhaskar--Rakhmanov--Saff number is a multiple of that for $W=\exp
(-Q) $. Since for any fixed $s>0,$
\[
a_{sn}\sim a_{n}\mbox{ uniformly in }n\mbox{,}
\] %
we obtain%
\begin{equation}
\left\Vert PW^{r}\right\Vert _{L_{\infty }(\mathbb{R})
}^{p}\leq C\frac{n}{a_{n}}\int_{-\infty }^{\infty }\left\vert
PW^{r}\right\vert ^{p},
\end{equation}%
and%
\begin{equation}
\left\Vert P^{\prime }W^{r}\right\Vert _{L_{\infty }(\mathbb{R}) }^{p}\leq C\left( \frac{n}{a_{n}}\right)
^{p+1}\int_{-\infty }^{\infty }\left\vert PW^{r}\right\vert ^{p},
\end{equation}%
for all $0<p\leq 2,$ $r>0$, $n\geq 1$, and $P$ of degree $\leq n$.

Now comes the sieving idea. Let $L,M$ be positive integers and fix
$\xi .$
We apply these inequalities to the polynomial in $t,$%
\begin{equation}
P(t) =S(t) \left( K_{Mn}( \xi
,t) \right) ^{L},
\end{equation}%
where $S$ has degree $\leq n$, and
\[
K_{n}( \xi ,t) =\sum_{j=0}^{n-1}p_{j}(\xi)
p_{j}(t)
\] %
is the $n$th reproducing kernel for the weight $W^{2}$. As
\[
\frac{(LM+1) n}{a_{(LM+1) n}}\sim
\frac{n}{a_{n}},
\] %
uniformly in $n$, we obtain from (8.26),
\[
\left\{ \left\vert S(\xi) K_{Mn}^{L}( \xi ,\xi
)
\right\vert W^{r}(\xi) \right\} ^{p}\leq C\frac{n}{a_{n}}%
\int_{-\infty }^{\infty }\left\vert( SW^{r}) (t) K_{Mn}( \xi ,t) ^{L}\right\vert ^{p}dt.
\] %
Here, for Freud weights $W\in \mathcal{F}^{\ast }$
\cite{LevinLubinsky1992},
\begin{equation}
K_{n}( \xi ,\xi) =1/\lambda _{n}( W^{2},\xi
) \sim \frac{n}{a_{n}}W^{-2}(\xi), \qquad
\left\vert \xi \right\vert \leq \frac{1}{2}a_{n},
\end{equation}%
while%
\begin{equation}
K_{n}( \xi ,\xi) =1/\lambda _{n}( W^{2},\xi) \leq C%
\frac{n}{a_{n}}W^{-2}(\xi), \qquad \xi \in
\mathbb{R},
\end{equation}%
so if $M$ is so large that
\[
\frac{1}{2}a_{Mn}\geq 2a_{n}\mbox{,}
\] %
we obtain for all polynomials $S$ of degree $\leq n$, and all
$\left\vert
\xi \right\vert \leq 2a_{n},$%
\[
\left\vert S(\xi) \right\vert ^{p}W^{(r-2L) p}(\xi) \leq C\left(
\frac{n}{a_{n}}\right) ^{1-Lp}\int_{-\infty }^{\infty }\left\vert
( SW^{r}) (t) K_{Mn}( \xi ,t)
^{L}\right\vert ^{p}dt.
\] %
Now assume that $L$ is chosen so large that $Lp>2$. By
Cauchy--Schwarz, for
all $n$ and $t,\xi \in \mathbb{R},$%
\begin{eqnarray}
\left\vert K_{Mn}( \xi ,t) \right\vert ^{Lp-2} &\leq
&\left( K_{Mn}( \xi ,\xi) ^{1/2}K_{Mn}( t,t) ^{1/2}\right)
^{Lp-2} \nonumber\\
&\leq &C\left( \frac{n}{a_{n}}\right) ^{Lp-2}(W^{-1}( \xi
) W^{-1}(t) )^{Lp-2}.
\end{eqnarray}%
So%
\begin{equation}
\left\vert S(\xi) \right\vert ^{p}W^{(r-L) p-2}(\xi) \leq
C\frac{a_{n}}{n}\int_{-\infty }^{\infty }\vert S(t)\vert^p K_{Mn}(
\xi ,t) ^{2}W^{(r-L) p+2}(t) dt.
\end{equation}%
We now choose $r=L+1$, giving%
\begin{equation}
\vert (SW)(\xi)\vert ^{p} \leq C\frac{a_{n}}{n}%
\int_{-\infty }^{\infty }\vert (SW)(t)\vert^{p} \left[ W(\xi) W(t)
K_{Mn}( \xi ,t\right)] ^{2}dt,
\end{equation}%
valid for all $0<p\leq 2,$ $n\geq 1$, $S$ of degree $\leq n$, and $%
\left\vert \xi \right\vert \leq 2a_{n}$. Since (by orthonormality)
\begin{equation}
\int_{-\infty }^{\infty }K_{Mn}^{2}( \xi ,t) W(t) ^{2}dt=K_{Mn}( \xi ,\xi) \sim
\frac{n}{a_{n}}W^{-2}(\xi) ,
\end{equation}%
we can also write this as:

\begin{lem} \textsl{Let }$0<p\leq 2,$\textsl{\ }$n\geq
1$\textsl{, }$S$\textsl{\ of
degree }$\leq n$\textsl{, and }$\left\vert \xi \right\vert \leq 2a_{n}$%
\textsl{. Then}%
\begin{equation}
\vert (SW)(\xi)\vert ^{p} \leq C\frac{\int_{-\infty }^{\infty
}\vert (SW)(t)\vert ^{p} K_{Mn}( \xi ,t) ^{2}W^{2}(t)
dt}{\int_{-\infty }^{\infty }K_{Mn}( \xi ,t) ^{2}W^{2}(t) dt}.
\end{equation}\end{lem}

From this follows:

\begin{lem}  {\bf (Fundamental lemma of sieving)}
\textsl{Let }$\psi :[0,\infty )\rightarrow \lbrack 0,\infty
)$\textsl{\ be a
convex increasing function of }$x$\textsl{, with }$\psi (0) =0.$%
\textsl{\ Let }$p>0,$\textsl{\ }$n\geq 1$\textsl{, }$S$\textsl{\ of degree }$%
\leq n$\textsl{, and }$\left\vert \xi \right\vert \leq 2a_{n}$\textsl{. Then}%
\begin{equation}
\psi(\vert SW\vert ^{p}) ( \xi
) \leq \frac{\int_{-\infty }^{\infty }\psi ( C\vert SW\vert ^{p}) (t) K_{Mn}(
\xi ,t) ^{2}W^{2}(t) dt}{\int_{-\infty }^{\infty
}K_{Mn}( \xi ,t) ^{2}W^{2}(t) dt}.
\end{equation}%
\textsl{If in addition, for some }$A>0,$\textsl{\ }%
\[
\psi (2t) \leq A\psi (t) ,\qquad t\in
\lbrack 0,\infty ),
\] %
\textsl{then}
\begin{equation}
\psi(\vert SW\vert ^{p}) (\xi) \leq C%
\frac{\int_{-\infty }^{\infty }\psi(\vert
SW\vert ^{p}) (t) K_{Mn}( \xi
,t) ^{2}W^{2}(t) dt}{\int_{-\infty }^{\infty
}K_{Mn}( \xi ,t) ^{2}W^{2}(t) dt}.
\end{equation}\end{lem}

\proof. For $p\leq 2$, this follows from (8.35)
by a single application of Jensen's inequality. For general $p$,
one instead applies Jensen's inequality with the convex function
$t\mapsto \psi \left( t^{a}\right) $, with large enough $a$. \endproof

Next, we extend this to derivatives. We again use the polynomial
$P$ of
(8.28). We let%
\[
K_{n}^{\prime }( x,t) =\frac{\partial }{\partial
x}K_{n}( x,t) .
\] %
By (8.25),
\begin{eqnarray*}
&&\hskip-4cm\left\vert S^{\prime }(\xi)
K_{Mn}^{L}( \xi ,\xi) +L\left( K_{Mn}( \xi ,\xi
) \right) ^{L-1}K_{Mn}^{\prime }( \xi ,\xi)
S(\xi) \right\vert ^{p}W^{rp}(\xi)
\\
&\leq &C\left( \frac{n}{a_{n}}\right) ^{p+1}\int_{-\infty
}^{\infty }\left\vert( SW^{r}) (t)
K_{Mn}( \xi ,t) ^{L}\right\vert ^{p}dt.
\end{eqnarray*}%
Here we recall (8.29) and (8.30), while Cauchy--Schwarz and (8.2)
give for all $x,t\in
\mathbb{R}$,%
\begin{eqnarray*}
\left\vert K_{Mn}^{\prime }( \xi ,t) \right\vert \leq
\left( \sum_{k=0}^{Mn-1}p_{k}^{\prime }(\xi)
^{2}\right) ^{1/2}\left( \sum_{k=0}^{Mn-1}p_{k}(t)
^{2}\right) ^{1/2} \leq C\left( \frac{n}{a_{n}}\right)
^{2}W^{-1}(\xi) W^{-1}(t) .
\end{eqnarray*}%
Then%
\begin{eqnarray*}
\left\vert S^{\prime }(\xi) \right\vert ^{p}W(\xi) ^{(r-L)p-2} &\leq &C\left( \frac{n}{a_{n}}\right)
^{p}\left\vert S( \xi
) \right\vert ^{p}W^{(r-L)p-2}(\xi) \\
&&+\ C\left( \frac{n}{a_{n}}\right) ^{p-1}\int_{-\infty }^{\infty
}\left\vert S(t) \right\vert ^{p}K_{Mn}( \xi
,t) ^{2}W^{(r-L) p+2}(t) dt \\
&\leq &C\left( \frac{n}{a_{n}}\right) ^{p-1}\int_{-\infty
}^{\infty }\left\vert S(t) \right\vert ^{p} K_{Mn}( \xi ,t)
^{2}W^{(r-L) p+2}(t) dt,
\end{eqnarray*}%
by (8.32). We now choose $r=L+1$, giving%
\begin{equation}
\vert S^{\prime }W\vert ^{p}(\xi) \leq C\left( \frac{%
n}{a_{n}}\right) ^{p-1}\int_{-\infty }^{\infty }\vert
SW\vert ^{p}(t) \left[ W(\xi)
W(t) K_{Mn}( \xi ,t) \right] ^{2}dt,
\end{equation}%
which we can reformulate as

\begin{lem}  \textsl{Let }$0<p\leq 2,$\textsl{\ }$n\geq
1$\textsl{, }$S$\textsl{\ of
degree }$\leq n$\textsl{, and }$\left\vert \xi \right\vert \leq 2a_{n}$%
\textsl{. Then}%
\begin{equation}
\vert S^{\prime }W\vert ^{p}(\xi) \leq C\left( \frac{%
n}{a_{n}}\right) ^{p}\frac{\int_{-\infty }^{\infty }\vert
SW\vert
^{p}(t) K_{Mn}( \xi ,t) ^{2}W^{2}(t) dt}{%
\int_{-\infty }^{\infty }K_{Mn}( \xi ,t)
^{2}W^{2}(t) dt}.
\end{equation}\end{lem}

From this follows:

\begin{lem}  \textsl{Let }$\psi :[0,\infty )\rightarrow
\lbrack 0,\infty )$\textsl{\ be a convex increasing function of
}$x$\textsl{, with }$\psi (0) =0$.
\textsl{Let }$p>0,$\textsl{\ }$n\geq 1$\textsl{, }$S$\textsl{\ of degree }$%
\leq n$\textsl{, and }$\left\vert \xi \right\vert \leq 2a_{n}$\textsl{. Then}%
\begin{equation}
\psi(\vert S^{\prime }W\vert ^{p})( \xi
) \leq \frac{\int_{-\infty }^{\infty }\psi( C\left\vert \frac{n}{%
a_{n}}SW\right\vert ^{p}) (t) K_{Mn}( \xi
,t) ^{2}W^{2}(t) dt}{\int_{-\infty }^{\infty
}K_{Mn}( \xi ,t) ^{2}W^{2}(t) dt}.
\end{equation}%
\textsl{If in addition, for some }$A>0,$\textsl{\ }%
\begin{equation}
\psi (2t) \leq A\psi (t) ,\mbox{ }\qquad
t\in \lbrack 0,\infty ),
\end{equation}%
\textsl{then}
\begin{equation}
\psi(\vert S^{\prime }W\vert ^{p})( \xi
) \leq C\frac{\int_{-\infty }^{\infty }\psi( \left\vert\frac{n}{a_{n}}%
SW\right\vert ^{p}) (t) K_{Mn}( \xi ,t) ^{2}W^{2}(t)
dt}{\int_{-\infty }^{\infty }K_{Mn}( \xi ,t) ^{2}W^{2}(t) dt}.
\end{equation}\end{lem}

\proof. The proof uses Jensen's inequality as in
Lemma 8.3.  \endproof

The inequalities (8.40) and (8.42) are useful for more than
Markov--Bernstein inequalities. But for the moment, we deduce by
integration and restricted range inequalities
\cite{LubinskyNevai1987}:

\begin{thh}  \textsl{Let} $W\in \mathcal{F}^{\ast }$. \textsl{Let
}$\psi :[0,\infty )\rightarrow
\lbrack 0,\infty )$\textsl{\ be a convex increasing function of }$x$\textsl{%
, with }$\psi (0) =0$. \textsl{Let }$p>0,$\textsl{\ }$n\geq 1$%
\textsl{, }$S$\textsl{\ of degree }$\leq n$. \textsl{Then}
\begin{equation}
\int_{-\infty }^{\infty }\psi(\vert S^{\prime
}W\vert ^{p}) (\xi) d\xi \leq
\int_{-\infty }^{\infty }\psi( C\frac{n}{a_{n}}\vert
SW\vert ^{p}) (t) dt.
\end{equation}%
\textsl{If in addition, for some }$A>0,$\textsl{\ (8.41) holds, then }%
\begin{equation}
\int_{-\infty }^{\infty }\psi(\vert S^{\prime
}W\vert ^{p}) (\xi) d\xi \leq
C\int_{-\infty }^{\infty }\psi( \frac{n}{a_{n}}\vert
SW\vert ^{p}) (t) dt.
\end{equation}\end{thh}

\sect{Nikolskii Inequalities}

We already proved an inequality of this type in the last section:%
\[
\left\Vert PW\right\Vert _{L_{\infty }(\mathbb{R})
}\leq C\left( \frac{n}{a_{n}}\right) ^{1/p}\left\Vert PW\right\Vert _{L_{p}\left( \mathbb{R%
}\right) },
\] %
for $n\geq 1$, and polynomials $P$ of degree $\leq n$ --- recall
(8.24) and (8.35). Thus we compared the weighted sup norm and
weighted $L_{p}$ norm of a polynomial. More generally,
inequalities that compare the norms of polynomials of degree $\leq
n$ in different spaces are called \dword{Nikolskii inequalities}.
They are not difficult to prove, here is a sample:

\begin{thh}  \textsl{Let }$W\in \mathcal{F}^{\ast }$\textsl{. Let }$0<p,r\leq \infty $%
\textsl{. Then for }$n\geq 1$\textsl{\ and polynomials
}$P$\textsl{\ of
degree }$\leq n,$%
\[
\left\Vert PW\right\Vert _{L_{p}(\mathbb{R}) }\leq
CN_{n}( p,r) \left\Vert PW\right\Vert _{L_{r}(\mathbb{R}) },
\] %
\textsl{where}%
\[
N_{n}( p,r) =\left\{
\begin{array}{ll}
a_{n}^{\frac{1}{p}-\frac{1}{r}}, & r>p \\[10pt]
\left( \frac{n}{a_{n}}\right) ^{\frac{1}{r}-\frac{1}{p}}, & r<p%
\end{array}%
\right. .
\] \end{thh}

\proof.
If first $\infty >r>p$, we can use restricted-range inequalities, and then H\"{o}lder's inequality:%
\begin{eqnarray*}
\left\Vert PW\right\Vert _{L_{p}(\mathbb{R}) }^{p}
\leq C\int_{-a_{n}}^{a_{n}}\left\vert PW\right\vert ^{p}
\leq C\left( \int_{-a_{n}}^{a_{n}}\left\vert PW\right\vert ^{p\frac{r}{p}%
}\right) ^{\frac{p}{r}}\left( \int_{-a_{n}}^{a_{n}}1\right)
^{1-\frac{p}{r}} \leq C\left\Vert PW\right\Vert _{L_{r}(\mathbb{R}) }^{p}a_{n}^{1-\frac{p}{r}},
\end{eqnarray*}%
that is%
\[
\left\Vert PW\right\Vert _{L_{p}(\mathbb{R}) }\leq Ca_{n}^{\frac{%
1}{p}-\frac{1}{r}}\left\Vert PW\right\Vert _{L_{r}(\mathbb{R}) }.
\] %
If $p=\infty $, the proof is easier. Next, if $r<p<\infty $, we
use
\begin{eqnarray*}
\left\Vert PW\right\Vert _{L_{p}(\mathbb{R}) }^{p}
=\int_{-\infty }^{\infty }\left\vert PW\right\vert ^{p} \leq
\left\Vert PW\right\Vert _{L_{\infty }(\mathbb{R})
}^{p-r}\int_{-\infty }^{\infty }\left\vert PW\right\vert ^{r}.
\end{eqnarray*}%
Letting%
\[
\Lambda _{n,p}:=\sup_{\deg (P) \leq n}\left( \frac{\left\Vert PW\right\Vert _{L_{\infty }\left(
\mathbb{R}\right) }}{\left\Vert PW\right\Vert _{L_{p}(\mathbb{R}) }}\right) ^{p},
\] %
we obtain%
\begin{equation}
\left\Vert PW\right\Vert _{L_{p}(\mathbb{R}) }\leq
\Lambda _{n,p}^{\frac{1}{r}-\frac{1}{p}}\left\Vert PW\right\Vert
_{L_{r}(\mathbb{R}) }^{p}.
\end{equation}%
In the case of Freud weights $W\in \mathcal{F}^{\ast }$ that grow
at least as fast as $\vert x\vert ^{\alpha }$, some
$\alpha >1$, we know from (8.24) that
\[
\Lambda _{n,p}\leq C\frac{n}{a_{n}}. \qquad \meop
\]

For the canonical weights $W_{\alpha }$, we have
\[
\Lambda _{n,p}\sim \left\{
\begin{array}{ll}
n^{1-1/\alpha }, & \alpha >1 \\
\log (n+1) , & \alpha =1 \\
1, & \alpha <1%
\end{array}%
\right. .
\] %
The sharpness of these was proved by Nevai and Totik \cite%
{Milovanovicetal1994}, \cite{NevaiTotik1987}. Observe that
$\Lambda _{n,p}$ grows independently of $p$. In general, it seems
that
\[
\Lambda _{n,p}\sim \sup_{x\in \mathbb{R}}\lambda _{n}^{-1}( W^{2},x) W^{2}(x).
\] %
For general convex exponential weights, it is known \cite[p.~295]%
{LevinLubinsky2001} that:

\begin{thh} \textsl{ Let }$W\in \mathcal{F}_{even}\left( C^{2}\right) $\textsl{\ be as in
Definition 7.5. Let }$0<p,r\leq \infty $\textsl{. Then for }$n\geq 1$\textsl{%
\ and polynomials }$P$\textsl{\ of degree }$\leq n,$%
\[
\left\Vert PW\right\Vert _{L_{p}(\mathbb{R}) }\leq
CN_{n}( p,r) \left\Vert PW\right\Vert _{L_{r}(\mathbb{R}) },
\] %
\textsl{where}%
\[
N_{n}( p,r) =\left\{
\begin{array}{ll}
a_{n}^{\frac{1}{p}-\frac{1}{r}}, & r>p \\[10pt]
\left( \frac{nT( a_{n}) ^{1/2}}{a_{n}}\right) ^{\frac{1}{r}-\frac{%
1}{p}}, & r<p%
\end{array}%
\right . .
\]\end{thh}

\sect{Orthogonal Expansions}

There is an old mathematical saying that $L_{1}$, $L_{2}$ and
$L_{\infty }$ were invented by the Almighty, and man invented all
else. The author heard this many years ago, but was interested to
see it used as the opening quote in a chapter of Simon's treatise
\cite{Simon2005}. That orthogonal expansions naturally live in
$L_{2}$ is obvious. Among the many manifestations of this, is the
best approximation property
\[
\left\Vert \left( f-S_{n}\left[ f\right] \right) W\right\Vert
_{L_{2}(\mathbb{R}) }=E_{n}\left[ f;W\right]
_{2}=\inf_{\deg (P)
\leq n}\left\Vert (f-P) W\right\Vert _{L_{2}\left( \mathbb{R}%
\right) }.
\] %
This ensures that when the polynomials are dense (in an obvious
sense)
\[
\lim_{n\rightarrow \infty }\left\Vert \left( f-S_{n}\left[
f\right] \right) W\right\Vert _{L_{2}(\mathbb{R}) }=0
\] %
for all functions $f$ for which $fW\in L_{2}(\mathbb{R}) $.

But it is part of the mathematician's spirit to take important
tools out of their natural domain, so it is not surprising that
much effort has been
devoted to convergence of $\left\{ S_{n}\left[ f\right] \right\} $ in $L_{p}$%
, or in a uniform sense, or at a specific point, and so on. A lot
of fundamental advances have ensued: for example, the boundedness
of the Hilbert transform in $L_{p}$, $1<p<\infty $, was
established in order to
prove that classic Fourier series converge in such $L_{p}$. The theory of $%
A_{p}$ weights started with Muckenhoupt's efforts to prove convergence in $%
L_{p}$ of Hermite expansions.

In this section, we shall discuss pointwise and mean convergence.
We begin with the latter.

\subsection{Mean Convergence}

Recall the reproducing kernel and the Christoffel--Darboux formula:%
\begin{eqnarray*}
K_{n}( x,t) =\sum_{j=0}^{n-1}p_{j}(x)
p_{j}(t) =\frac{\gamma _{n-1}}{\gamma
_{n}}\frac{p_{n}(x)p_{n-1}(t)
-p_{n-1}(x)p_{n}(t) }{x-t}.
\end{eqnarray*}%
Define the Hilbert transform%
\begin{eqnarray*}
H\left[ g\right] (x)=\int_{-\infty }^{\infty
}\frac{g(t) }{x-t}dt =\lim_{\varepsilon \rightarrow
0+}\int_{\left( -\infty ,\infty \right) \backslash \left( x-\varepsilon ,x+\varepsilon \right) }\frac{g(t)
}{x-t}dt,
\end{eqnarray*}%
whenever the limit exists. If $g\in L_{1}(\mathbb{R})
$, the transform exists for a.e.\ $x$. For a function $f$ with
$fW\in L_{2}(\mathbb{R}) $, we see that
\begin{eqnarray*}
S_{n}\left[ f\right] (x)&=&\int_{-\infty }^{\infty
}K_{n}( x,t) f(t) W^{2}(t) dt \\
&=&\frac{\gamma _{n-1}}{\gamma _{n}}\left\{ p_{n}(x)
H\left[ p_{n-1}fW^{2}\right] (x)-p_{n-1}(x) H\left[ p_{n}fW^{2}\right] (x)\right\} .
\end{eqnarray*}%
Thus if $u$ is a given function, and $1<p<\infty ,$%
\begin{eqnarray}
\left\Vert S_{n}\left[ f\right] Wu^{2}\right\Vert _{L_{p}\left( \mathbb{R}%
\right) } &\leq &\frac{\gamma _{n-1}}{\gamma _{n}}\left\{
\left\Vert p_{n}Wu\right\Vert _{L_{\infty }(\mathbb{R}) }\left\Vert H\left[ p_{n-1}fW^{2}\right]
u\right\Vert _{L_{p}(\mathbb{R})
}\right.\nonumber\\
&&\qquad\left.+\left\Vert p_{n-1}Wu\right\Vert _{L_{\infty }(\mathbb{R})
}\left\Vert H\left[ p_{n}fW^{2}\right] u\right\Vert _{L_{p}\left( \mathbb{R}%
\right) }\right\} .
\end{eqnarray}%
Next, we use the aforementioned theorem of Riesz that the Hilbert
transform is a bounded operator on $L_{p}$, provided $1<p<\infty
$: For some $C\neq
C(g) ,$%
\[
\left\Vert H\left[ g\right] \right\Vert _{L_{p}(\mathbb{R}) }\leq C\left\Vert g\right\Vert _{L_{p}\left(
\mathbb{R}\right) }.
\] %
In his investigations of Hermite expansions, Muckenhoupt \cite%
{Muckenhoupt1970} considered weighted versions: for suitable functions $u,$%
\[
\left\Vert H\left[ g\right] u\right\Vert _{L_{p}(\mathbb{R}) }\leq C\left\Vert gu\right\Vert _{L_{p}\left(
\mathbb{R}\right) }.
\] %
Such $u$ are severely restricted, in particular satisfying
Muckenhoupt's condition; examples are
\begin{equation}
u(t) =\left( 1+\vert t\vert \right) ^{a}\mbox{, }\qquad -%
\frac{1}{p}<a<1-\frac{1}{p}.
\end{equation}%
The theorem fails if $p=1$ or $p=\infty $, though in $L_{1}$,
other inequalities are available. So we continue (10.1) as
\begin{eqnarray}
\left\Vert S_{n}\left[ f\right] Wu^{2}\right\Vert _{L_{p}( \mathbb{R}%
) } &\leq &C\frac{\gamma _{n-1}}{\gamma _{n}}\left\{
\left\Vert p_{n}Wu\right\Vert _{L_{\infty }(\mathbb{R}) }\left\Vert
p_{n-1}fW^{2}u\right\Vert _{L_{p}(\mathbb{R}) }\right.\nonumber\\
&&\qquad\left.+\left\Vert p_{n-1}Wu\right\Vert _{L_{\infty }(\mathbb{R}) }\left\Vert p_{n}fW^{2}u\right\Vert
_{L_{p}(\mathbb{R}) }\right\} .
\end{eqnarray}%

For exponential weights in the real line, recall from (4.14) that
\[
\frac{\gamma _{n-1}}{\gamma _{n}}\leq Ca_{n}.
\] %
If we have the bounds
\begin{equation}
\left\Vert p_{n}Wu\right\Vert _{L_{\infty }(\mathbb{R}) }\leq Ca_{n}^{-1/2},
\end{equation}%
then for some $C\neq C\left( n,f\right) ,$%
\[
\left\Vert S_{n}\left[ f\right] Wu^{2}\right\Vert _{L_{p}( \mathbb{R}%
) }\leq C\left\Vert fW\right\Vert _{L_{p}(\mathbb{R}) }.
\] %
Once we have such a bound, the reproducing property of $S_{n}$ (that $S_{n}%
\left[ P\right] =P$ when $P$ has degree $\leq n$) and density of
polynomials
gives%
\[
\lim_{n\rightarrow \infty }\left\Vert \left( f-S_{n}(f) \right) Wu^{2}\right\Vert _{L_{p}(\mathbb{R})
}=0.
\] %
All that is required of $f$ is that $fW\in L_{p}(\mathbb{R}) .$

The problem with this procedure is that the bound (10.4) is hardly
ever
true. For Freud weights like $W_{\alpha }$, or the weights in the class $%
\mathcal{F}^{\ast }$, we typically have
\begin{equation}
\left\vert p_{n}(x)\right\vert W(x)\leq
Ca_{n}^{-1/2}\left( \left\vert 1-\frac{\vert x\vert }{a_{n}}%
\right\vert +n^{-2/3}\right) ^{-1/4},\qquad \mbox{ }x\in
\mathbb{R}
\end{equation}%
and the upper bound reflects the real growth of $p_{n}$ near $\pm
a_{n}$:
\[
\left\Vert p_{n}W\right\Vert _{L_{\infty }(\mathbb{R})
}\sim a_{n}^{-1/2}n^{1/6}.
\] %
However, at least in $\left[ -\rho a_{n},\rho a_{n}\right] $, with
$\rho \in \left( 0,1\right) $ fixed, we do have
\begin{equation}
\left\vert p_{n}W\right\vert \leq Ca_{n}^{-1/2}.
\end{equation}%
This problematic growth of the orthogonal polynomials near the
(effective) endpoints of the interval of orthogonality, already
occurs for the Legendre weight on $\left[ -1,1\right] $, and more
generally, Jacobi weights. The same Pollard that solved
Bernstein's approximation problem, came up with a fix for this. In
the context of Freud weights (a similar feature occurs for Jacobi
weights), the fix is based on the observation that $p_{n}-p_{n-2}$
has a much better bound than $p_{n}$:%
\begin{equation}
\left\vert p_{n}(x)-p_{n-2}(x)\right\vert
W(x)\leq
Ca_{n}^{-1/2}\left( \left\vert 1-\frac{\vert x\vert }{a_{n}}%
\right\vert +n^{-2/3}\right) ^{1/4},\qquad \mbox{ }x\in
\mathbb{R}\mbox{.}
\end{equation}%
In particular,%
\[
\left\vert p_{n}(x)-p_{n-2}(x)\right\vert
W(x)\leq Ca_{n}^{-1/2},\qquad \mbox{ }x\in
\mathbb{R}\mbox{. }
\] %
(Think of $\widetilde{p}_{n}\left( \cos \theta \right) =\cos
n\theta $ to see from whence this comes.)

Pollard \cite{Pollard1947}, \cite{Pollard1948}, \cite{Pollard1949}
found a clever way to rewrite the Christoffel--Darboux formula to
exploit this: let
us set%
\[
\alpha _{n}=\frac{\gamma _{n-1}}{\gamma _{n}}.
\] %
Then (see \cite{Lubinsky1995} or \cite{MhaskarXu1991} for an
accessible proof)
\[
K_{n}( x,y) =K_{n,1}(x,y) +K_{n,2}(x,y) +K_{n,3}(x,y) ,
\] %
where
\begin{eqnarray*}
K_{n,1}(x,y) &=&\frac{\alpha _{n}}{\alpha _{n}+\alpha _{n-1}}%
p_{n-1}(x)p_{n-1}(y) ; \\
K_{n,2}(x,y) &=&\frac{\alpha _{n}\alpha _{n-1}}{\alpha
_{n}+\alpha _{n-1}}p_{n-1}(y) \left[ \frac{p_{n}(x)
-p_{n-2}(x)}{x-y}\right] ; \\
K_{n,3}(x,y) &=&K_{n,2}( y,x) .
\end{eqnarray*}%
Note that there is no $x-y$ in the denominator in $K_{n,1}$, while in $%
K_{n,2}$ and $K_{n,3}$, we have the term $p_{n}-p_{n-2}$ to help.
This clever decomposition is not enough on its own; we emphasize
that in investigating mean convergence, one still has to break up
integrals into several different pieces, and work over several
different ranges.

For the Hermite weight, the simplest bound is due to Askey and Wainger \cite%
{AskeyWainger1965}:

\begin{thh}
\textsl{Let }$W(x)=\exp( -x^{2}) $\textsl{. Let }$%
\frac{4}{3}<p<4$\textsl{. There exists }$C\neq C\left( n,f\right)
$\textsl{\
such that }%
\[
\left\Vert S_{n}\left[ f\right] W\right\Vert _{L_{p}(\mathbb{R}) }\leq C\left\Vert fW\right\Vert _{L_{p}(
\mathbb{R}) }.
\]
\textsl{This inequality is not true for }$p\leq \frac{4}{3}$\textsl{\ or }$%
p\geq 4$\textsl{. }\end{thh}

Muckenhoupt found the correct extension to $1<p<\infty$. In what
follows we let $u_a(x) = (1 + |x|)^a$.

\begin{thh} \textsl{Let }$W(x)=\exp( -x^{2}) $\textsl{\ and }$%
1<p<\infty $\textsl{. Let }%
\begin{equation}
b<1-\frac{1}{p};\quad B>-\frac{1}{p};\quad b\leq B.
\end{equation}%
\textsl{Assume in addition that }%
\begin{equation}
-B+\max \left\{ b,-\frac{1}{p}\right\} +\frac{4}{3p}-1\leq 0\quad \mbox{ if }p<%
\frac{4}{3};
\end{equation}%
\textsl{and}%
\begin{equation}
b-\min \left\{ B,1-\frac{1}{p}\right\}
+\frac{1}{3}-\frac{4}{3p}\leq 0\quad \mbox{ if }p>4.
\end{equation}%
\textsl{Then \newline
}%
\begin{equation}
\left\Vert S_{n}\left[ f\right] Wu_{b}\right\Vert _{L_{p}( \mathbb{R}%
) }\leq C\left\Vert fWu_{B}\right\Vert _{L_{p}(\mathbb{R}) }
\end{equation}
\textsl{for some }$C\neq
C\left( n,f\right) $\textsl{. If }$b=B$\textsl{\
and }$p=\frac{4}{3}$\textsl{\ or }$4$\textsl{, then we insert a factor of }$%
\log \left( \vert x\vert +2\right) $\textsl{\ in the
right-hand side of (10.11). In the case of equality in (10.9) or
(10.10), we need strict inequality and replace the }$\max $
\textsl{or} $\min $ \textsl{by their second terms. All these
inequalities are also necessary.}\end{thh}

For general Freud weights, Shing Wu Jha and the author
\cite{JhaLubinsky1995} extended Muckenhoupt's result, also closing
up slight gaps between the latter's necessary and sufficient
conditions in ``boundary'' cases: we let
\[
L_{\sigma ,\tau }(n) =\left\{
\begin{array}{ll}
\left( \log (n+1) \right) ^{\left\vert \sigma
\right\vert }, &
\mbox{if }\sigma =\tau \\
1, & \mbox{otherwise.}%
\end{array}%
\right.
\] %

\begin{thh}
\textsl{Let} $W=\exp (-Q) \in \mathcal{F}^{\ast }$. \textsl{Let }%
$1<p<\infty $,\textsl{\ and\ }$b, B\in \mathbb{R}$. \newline
\begin{description}
\item[(a)]\textsl{ Then for }%
\[
\left\Vert S_{n}\left[ f\right] Wu_{b}\right\Vert _{L_{p}( \mathbb{R}%
) }\leq C\left\Vert fWu_{B}\right\Vert _{L_{p}(\mathbb{R}) }
\] %
\textsl{to hold with some }$C\neq C\left( n,f\right) $\textsl{, it
is necessary that:}
\begin{description}
\item[(I)] $(10.8)$ \textsl{holds.}
\item[(II)] \textsl{If} $p<\frac{4}{3}$\textsl{, then }%
\[
a_{n}^{\max \left\{ b,-\frac{1}{p}\right\} -B}n^{\frac{1}{6}\left( \frac{4}{p%
}-3\right) }=O\left( \frac{1}{L_{b,-\frac{1}{p}(n)
}}\right) .
\] %
\item[(III)]\textsl{ If }$p=\frac{4}{3}$\textsl{\ or }$p=4$\textsl{, then
strict inequality holds in the third inequality in $(10.8)$.}
\item[(IV)] \textsl{If }$p>4,$ \textsl{ then}%
\[
a_{n}^{b-\min \left\{ B,1-\frac{1}{p}\right\} -B}n^{\frac{1}{6}\left( 1-%
\frac{4}{p}\right) }=O\left( \frac{1}{L_{B,1-\frac{1}{p}(n) }}%
\right) .
\]
\end{description}
\item[(b)]\textsl{ Assume in addition that the orthonormal polynomials for }$W^{2}$%
\textsl{\ also satisfy $(10.7)$. Then the conditions above are also sufficient.%
}\end{description}\end{thh}

When $W=W_{\alpha }$, and $a_{n}=Cn^{1/\alpha }$, the conditions
above
reduce to essentially Muckenhoupt's:%
\[
\max \left\{ b,-\frac{1}{p}\right\} -B+\frac{\alpha }{6}\left( \frac{4}{p}%
-3\right) \left\{
\begin{array}{ll}
\leq 0, & b\neq -\frac{1}{p} \\
<0, & b=-\frac{1}{p}%
\end{array}%
.\right.
\] %
\[
b-\min \left\{ B,1-\frac{1}{p}\right\} -B+\frac{\alpha }{6}\left( 1-\frac{4}{%
p}\right) \left\{
\begin{array}{ll}
\leq 0, & B\neq 1-\frac{1}{p} \\
<0, & B=1-\frac{1}{p}%
\end{array}%
.\right.
\] %
The bound (10.7) was established for exp$\left( -x^{2m}\right) $, $%
m=1,2,3,\ldots$ in \cite{JhaLubinsky1995}. For general $\alpha
>1$, it follows
from results in that paper and work of Kriecherbauer and McLaughlin \cite%
{KriecherbauerMcLaughlin1999}. It is an interesting unsolved
problem to establish the bound (10.7) for $W\in \mathcal{F}^{\ast
}$, or to come up
with some method to prove it that works without much deeper tools. See \cite%
{KubayiMashele2003} for some work in this direction. There is also
an $L_{1}$ analogue of these \cite[p.~337]{JhaLubinsky1995}.

Further work on mean convergence of orthogonal expansions appears in \cite%
{Jungetal1999}, \cite{Ky1999}, \cite{MastroianniVertesi2006}, \cite%
{Mhaskar1984B}, \cite{Mhaskar1988}. Discrete analogues of
orthogonal
expansions have been considered by Mashele \cite{Mashele2002}, \cite%
{Mashele2002B}.

\subsection{Uniform and Pointwise Convergence}

One of the simplest and yet most elegant and general theorems on
pointwise convergence is an extension of the Dirichlet--Jordan
criterion, due to Freud
\cite{Freud1974}. It was subsequently extended by Mhaskar \cite{Mhaskar1984}%
, \cite{Mhaskar1988}. Its proof is quite simple, and follows
readily from estimates of previous sections, so we provide
it:

\begin{thh} {\bf  (Dirichlet--Jordan Criterion)}
\textsl{Let }$W\in \mathcal{F}^{\ast }$\textsl{\ and let }$f:\mathbb{R}%
\rightarrow \mathbb{R}$\textsl{\ be absolutely continuous, with
}$f^{\prime
}W\in L_{1}(\mathbb{R}) $\textsl{. Then}%
\[
\lim_{n\rightarrow \infty }\left\Vert \left( f-S_{n}\left[
f\right] \right) W\right\Vert _{L_{\infty }(\mathbb{R}) }=0.
\]\end{thh}

\proof. Note first that $fW$ is bounded and has
limit $0$ at $\pm \infty $. To see this let $A>0$, and $x\geq A$.
Then
\begin{eqnarray*}
\left\vert fW\right\vert (x)=\left\vert f(A) +\int_{A}^{x}f^{\prime }\right\vert W(x)\leq
\left\vert f(A) \right\vert W(x)
+\int_{A}^{\infty }\left\vert f^{\prime }\right\vert W,
\end{eqnarray*}%
as $W$ is decreasing. Hence%
\[
\limsup_{x\rightarrow \infty }\left\vert fW\right\vert (x) \leq \int_{A}^{\infty }\left\vert f^{\prime }\right\vert
W.
\] %
Our integrability condition on $f^{\prime }W$ ensures that this
last right-hand side can be made as small as we please. So,
indeed, $fW$ has limit $0$ at $\pm \infty .$

Next, let $m=\floor{n/2} $. Then
\begin{eqnarray}
\left\Vert \left( f-S_{n}\left[ f\right] \right) W\right\Vert
_{L_{\infty }(\mathbb{R}) } \leq\left\Vert \left( f-V_{m}\left[ f\right] \right) W\right\Vert
_{L_{\infty }(\mathbb{R}) }+\left\Vert \left( V_{m}\left[ f%
\right] -S_{n}\left[ f\right] \right) W\right\Vert _{L_{\infty
}(\mathbb{R}) }.
\end{eqnarray}%
Here by Nikolskii inequalities, such as in Theorem 9.1,
\begin{eqnarray}
\left\Vert \left( V_{m}\left[ f\right] -S_{n}\left[ f\right]
\right) W\right\Vert _{L_{\infty }(\mathbb{R}) }\
&\leq& C\sqrt{\frac{n}{a_{n}}}\left\Vert \left( V_{m}\left[ f\right] -S_{n}%
\left[ f\right] \right) W\right\Vert _{L_{2}(\mathbb{R}) } \nonumber\\[10pt]
&\leq &C\sqrt{\frac{n}{a_{n}}}\left\{ \left\Vert \left( V_{m}\left[ f\right] -f\right) W\right\Vert _{L_{2}(
\mathbb{R}) }+\left\Vert \left( f-S_{n}\left[ f\right]
\right) W\right\Vert _{L_{2}(\mathbb{R})
}\right\}  \nonumber\\[10pt]
&\leq &2C\sqrt{\frac{n}{a_{n}}}\left\Vert \left( V_{m}\left[
f\right] -f\right) W\right\Vert _{L_{2}(\mathbb{R}) }.
\end{eqnarray}%
Recall the $L_{2}$ best approximation property of the partial sums
$S_{n}$, and that $V_{m}$ has degree $\leq 2m\leq n$. Next,
\begin{eqnarray*}
\left\Vert \left( V_{m}\left[ f\right] -f\right) W\right\Vert
_{L_{2}(\mathbb{R}) }^{2} &\leq &\left\Vert \left( V_{m}\left[ f\right] -f\right) W\right\Vert _{L_{1}(
\mathbb{R}) }\left\Vert \left( V_{m}\left[ f\right]
-f\right) W\right\Vert _{L_{\infty }(\mathbb{R}) } \\[10pt]
&\leq &CE_{m}\left[ f;W\right] _{1}\left\Vert \left( V_{m}\left[
f\right]
-f\right) W\right\Vert _{L_{\infty }(\mathbb{R}) } \\[10pt]
&\leq& C\frac{a_{m}}{m}\left\Vert f^{\prime }W\right\Vert
_{L_{1}(\mathbb{R}) }\left\Vert \left( V_{m}\left[
f\right] -f\right) W\right\Vert _{L_{\infty }(\mathbb{R}) }.
\end{eqnarray*}%
We used the Jackson--Favard inequality (4.5) in the last line. Since $m\sim n$%
, so\ $\frac{a_{m}}{m}\sim \frac{a_{n}}{n}$. By substituting this
into
(10.13), we obtain%
\begin{eqnarray*}
\left\Vert \left( V_{m}\left[ f\right] -S_{n}\left[ f\right]
\right) W\right\Vert _{L_{\infty }(\mathbb{R}) }
&\leq &C\left[ \left\Vert f^{\prime }W\right\Vert _{L_{1}( \mathbb{R}%
) }\left\Vert \left( V_{m}\left[ f\right] -f\right)
W\right\Vert _{L_{\infty }(\mathbb{R}) }\right]
^{1/2},
\end{eqnarray*}%
so (10.12) becomes%
\begin{eqnarray*}
\left\Vert \left( f-S_{n}\left[ f\right] \right) W\right\Vert
_{L_{\infty }(\mathbb{R}) } &\leq &\left\Vert \left( f-V_{m}\left[ f\right] \right) W\right\Vert _{L_{\infty }(
\mathbb{R}) }+C\left[ \left\Vert f^{\prime }W\right\Vert
_{L_{1}(\mathbb{R}) }\left\Vert \left( V_{m}\left[
f\right] -f\right) W\right\Vert _{L_{\infty }(\mathbb{R}) }%
\right] ^{1/2} \\[10pt]
&\leq &C_{1}\left( E_{m}\left[ f;W\right] _{\infty }+E_{m}\left[
f;W\right] _{\infty }^{1/2}\right) ,
\end{eqnarray*}%
by Theorem 4.8 again. Since $fW$ is continuous, and has limit $0$
at
infinity, we have $E_{m}\left[ f;W\right] _{\infty }\rightarrow 0$ as $%
n\rightarrow \infty $.  \endproof

This theorem is attractive because of the simplicity of the
hypotheses. We don't need to assume anything about the ``tail'' of
$f$. There are also results which give convergence at a point,
under local assumptions. A result of Mhaskar
\cite[p.~224]{Mhaskar1996} is the archetype:

\begin{thh} \textsl{Let }$W\in \mathcal{F}^{\ast
}$\textsl{\ and in addition assume that
}$Q^{\prime \prime }$\textsl{\ is increasing on }$\left( 0,\infty \right) $%
\textsl{. Let }$f$\textsl{\ be of bounded variation, and
}$x$\textsl{\ be a
point of continuity of }$f$\textsl{. Then }%
\begin{eqnarray*}
\left\vert S_{n}\left[ f\right] (x)-f(x)
\right\vert
\leq C(x)\left[ \int_{x-\sqrt{\frac{a_{n}}{n}}}^{x+\sqrt{%
\frac{a_{n}}{n}}}W\left\vert df\right\vert +\sqrt{\frac{a_{n}}{n}}%
\int_{-\infty }^{\infty }W\left\vert df\right\vert
+\int_{\vert t\vert \geq Ca_{n}}W\left\vert
df\right\vert \right] .
\end{eqnarray*}%
\textsl{Here }$C(x)$\textsl{\ depends on
}$x$\textsl{, but not on }$n$\textsl{\ or }$f$\textsl{.}\end{thh}

Mashele \cite{Mashele2004} has established analogues of this for
general exponential weights. Unfortunately in passing to the more
general weights, some of the power of this result is lost. Ky
\cite{Ky1999} has investigated a.e.\ convergence of orthonormal
expansions for Freud weights. Trilinear
kernels and related expansions have been explored by Osilenker \cite%
{Osilenker1999}.

There is an extensive literature on expansions in Hermite
polynomials,
including transplantation theorems \cite{Thangavelu1993}, \cite%
{Thangavelu1998} that specifically exploit the structure of the
Hermite weight. There are also many results for Laguerre and
related expansions, on the interval $[0,\infty )$, see
\cite{MastroianniOccorsio2005}.

\sect{Lagrange Interpolation}

Lagrange interpolation at zeros of orthogonal polynomials is
probably the most studied of all approximation processes
associated with exponential weights. Recall that we denoted the
zeros of $p_{n}$, the $n$th orthonormal polynomial for the weight
$W^{2}$, by
\[
-\infty <x_{nn}<x_{n-1,n}<\cdots <x_{2n}<x_{1n}<\infty .
\] %
Given a function $f$, we define the $n$th Lagrange interpolation
polynomial to $f$ at $\left\{ x_{jn}\right\} _{j=1}^{n}$ to be the
unique polynomial of
degree $\leq n-1$ that satisfies%
\[
L_{n}\left[ f\right] \left( x_{jn}\right) =f\left( x_{jn}\right) \mbox{, }\qquad %
1\leq j\leq n\mbox{.}
\] %
One formula for $L_{n}\left[ f\right] $ is
\[
L_{n}\left[ f\right] =\sum_{j=1}^{n}f\left( x_{jn}\right) \ell
_{jn}
\] %
where $\left\{ \ell _{jn}\right\} _{j=1}^{n}$ are fundamental
polynomials of Lagrange interpolation, satisfying
\[
\ell _{jn}( x_{kn}) =\delta _{jk}.
\] %
The $\left\{ \ell _{jn}\right\} $ admit the representation%
\[
\ell _{jn}(x)=\prod\limits_{k=1,k\neq j}^{n}\frac{x-x_{kn}}{%
x_{jn}-x_{kn}},
\] %
but because we are using zeros of orthogonal polynomials, there are others:%
\[
\ell _{jn}(x)=\frac{K_{n}( x,x_{jn})
}{K_{n}( x_{jn},x_{jn}) }.
\] %
(To see this, look at the Christoffel--Darboux formula.)

Another advantage of interpolating at zeros of orthogonal
polynomials is the ease of analysis, largely because we have the
Gauss quadrature formula. For
polynomials $P$ of degree $\leq 2n-1,$%
\[
\int_{-\infty }^{\infty }PW^{2}=\sum_{j=1}^{n}\lambda _{jn}P( x_{jn}) ,
\] %
where
\[
\lambda _{jn}=\lambda _{n}( W^{2},x_{jn})
=1/K_{n}( x_{jn},x_{jn})
\] %
are called the \dword{Christoffel number}s. Since $L_{n}\left[ f\right]
^{2}$ has degree $\leq 2n-2$, and agrees with $f$ at $\left\{
x_{jn}\right\} $, we see that
\begin{equation}
\int_{-\infty }^{\infty }L_{n}\left[ f\right] ^{2}W^{2}=\sum_{j=1}^{n}%
\lambda _{jn}f^{2}( x_{jn}) .
\end{equation}%
Thus boundedness and convergence of $\left\{ L_{n}\left[ f\right]
\right\} $ is closely associated with convergence of Gauss
quadrature. There are at
least two recent monographs devoted to Lagrange interpolation \cite%
{MastroianniMilovanovic2006}, \cite{SzabadosVertesi1990}.

\subsection{Mean Convergence}

Just as $\left\{ S_{n}\left[ f\right] \right\} $ converge in
$L_{2}$, so do
the $\left\{ L_{n}\left[ f\right] \right\} $, a famous theorem of Erd\H{o}%
s-Tur\'{a}n. Indeed, we may think of $L_{n}\left[ f\right] $ as a
discretization of $S_{n}\left[ f\right]:$%
\begin{eqnarray*}
L_{n}\left[ f\right] (x)=\sum_{j=1}^{n}\lambda
_{jn}f\left( x_{jn}\right) K_{n}( x,x_{jn}) \approx
\int_{-\infty }^{\infty }W^{2}(t) f(t)
K_{n}( x,t) dt=S_{n}\left[ f\right] (x).
\end{eqnarray*}

At this stage, may I share two Erd\H{o}s stories with the reader.
I was doing a postdoc at the Technion in Haifa in 1982. Erd\H{o}s
often visited the Technion, and I was one of several who went to
lunch with him. Later that day, I carried his suit to the
drycleaners, and the reward was some time with him. He asked me
what I was doing, and I mentioned numerical quadrature. He replied
that I should go and read the book of Freud on orthogonal
polynomials, which I did, and it changed my career.

Here is another story, which illustrates how little we communicate
across fields. In early 2006, I had the privilege of visiting
Texas A\&M in College Station, to give some Frontiers Lectures. I
mentioned the famous theorem of Erd\H{o}s--Tur\'{a}n on
distribution of zeros of polynomials. At the end of the talk, one
of the very first comments was ``I thought Erd\H{o}s was a graph
theorist". Hmmph. He spent a large amount of his life on Lagrange
interpolation, orthogonal polynomials and approximation theory!

Since Lagrange interpolation samples a function only at a discrete
sets of points, and sets of measure zero make no difference to
Lebesgue integrals, it is clear that we shall need to use
something other than Lebesgue integrability. As Gauss quadrature
sums are Riemann sums or Riemann--Stieltjes sums
\cite[p.~50]{Szego1939}, the Riemann integral is the
appropriate tool. The following result was actually proved by Shohat \cite%
{Shohat1937}; Erd\H{o}s--Tur\'{a}n considered only finite
intervals. The
growth of $f$ at $\infty $ poses extra problems for the infinite interval.%

\begin{thh} {\bf (Shohat's version of the Erd\H{o}s--Tur\'{a}n Theorem)} %
\textsl{Let }$f$\textsl{\ be bounded and Riemann integrable on
each finite
interval. Assume, moreover, that there is an even entire function }$G$\textsl{%
\ with all non-negative Maclaurin series coefficients such that }$GW^{2}$%
\textsl{\ is integrable over the real line, and }%
\begin{equation}
\lim_{\vert x\vert \rightarrow \infty }f^{2}(x) /G(x)=0.
\end{equation}%
\textsl{Then}
\begin{equation}
\lim_{n\rightarrow \infty }\int_{-\infty }^{\infty }\left( f-L_{n}\left[ f%
\right] \right) ^{2}W^{2}=0.
\end{equation}\end{thh}

\proof. Let $P$ be a fixed polynomial. For large
enough $n$, $L_{n}\left[ P\right] =P $, so
\begin{eqnarray}
\int_{-\infty }^{\infty }\left( f-L_{n}\left[ f\right] \right)
^{2}W^{2} \nonumber &=&\int_{-\infty }^{\infty }\left( f-P-L_{n}\left[ f-P\right] \right)
^{2}W^{2}  \nonumber \\[10pt]
&\leq &2\int_{-\infty }^{\infty }(f-P)
^{2}W^{2}+2\int_{-\infty
}^{\infty }\left( L_{n}\left[ f-P\right] \right) ^{2}W^{2}  \nonumber \\[10pt]
&=&2\int_{-\infty }^{\infty }(f-P)
^{2}W^{2}+\sum_{j=1}^{n}\lambda _{jn}(f-P) ^{2}( x_{jn}) .
\end{eqnarray}%
In the second last line, we used the inequality
\[
(x+y) ^{2}\leq 2\left( x^{2}+y^{2}\right)
\] %
and in the last line we used (11.1). Now if $G$ is not a
polynomial (which
we may assume, just add a very slowly growing entire function if it is),%
\[
\lim_{\vert x\vert \rightarrow \infty }(f-P) ^{2}(x)/G(x)=0.
\] %
Let $\varepsilon >0$. Choose $A>0$ such that
\[
\vert x\vert \geq A\;\Rightarrow\; (f-P)
^{2}(x)/G(x)\leq \varepsilon .
\] %
Then%
\begin{eqnarray*}
\sum_{j:\left\vert x_{jn}\right\vert \geq A}\lambda _{jn}(f-P) ^{2}( x_{jn}) \leq \varepsilon
\sum_{j:\left\vert x_{jn}\right\vert \geq A}\lambda _{jn}G\left( x_{jn}\right) \leq \varepsilon \sum_{j=1}^{n}\lambda _{jn}G\left(
x_{jn}\right) \leq \varepsilon \int_{-\infty }^{\infty }GW^{2},
\end{eqnarray*}%
by a special case of the Posse--Markov--Stieltjes inequalities \cite[p.~92]%
{Freud1971}. These inequalities estimate quadrature sums of
absolutely monotone functions in terms of integrals. Next, since
$(f-P) ^{2}\chi _{\left[ -A,A\right] }$ is Riemann
integrable, while Gauss
quadrature sums are Riemann sums,%
\[
\lim_{n\rightarrow \infty }\sum_{j:\left\vert x_{jn}\right\vert
\leq A}\lambda _{jn}(f-P) ^{2}( x_{jn})
=\int_{-A}^{A}(f-P) ^{2}W^{2}.
\] %
Substituting all these in (11.4) gives%
\begin{eqnarray*}
\limsup_{n\rightarrow \infty }\int_{-\infty }^{\infty }\left( f-L_{n}\left[ f\right] \right) ^{2}W^{2} \leq 4\int_{-\infty
}^{\infty }(f-P) ^{2}W^{2}+\varepsilon \int_{-\infty
}^{\infty }GW^{2}.
\end{eqnarray*}%
Here $P$ is independent of $\varepsilon $, and by $L_{2}$
analogues of the Bernstein problem, may be chosen so that
$\int_{-\infty }^{\infty }(f-P) ^{2}W^{2}$ is as small
as we please.  \endproof

Note that if $f$ is of polynomial growth at $\infty $, we can always choose $%
G$ to satisfy (11.2). For large classes of exponential weights,
the author \cite{Lubinsky1986}, \cite{Lubinsky1989} constructed
entire functions $G$ satisfying the hypotheses above, and with
\[
G(x)W^{2}(x)\sim \left( 1+\left\vert
x\right\vert \right) ^{-a},
\] %
any $a>1$. So typically, the growth restriction on $f$ is hardly
more than that required for $(fW) ^{2}$ to be
integrable.

Of course this is a really beautiful result, and its proof is
elegant! The extensions to $L_{p}$ require a whole lot more work.
But like their cousins in orthogonal expansions, a lot of great
mathematics has come out of these efforts. It was Paul Nevai
\cite{Nevai1980} who obtained the first major result in $L_{p}$,
for the Hermite weight.

\begin{thh}  \textsl{Let }$f$\textsl{\ be continuous on
the real line. Let }$W(x)=\exp \left( -\frac{1}{2}x^{2}\right) $\textsl{\ be the Hermite
weight, and assume that }%
\[
\lim_{\vert x\vert \rightarrow \infty }f(x)
\left( 1+\vert x\vert \right) W(x)=0.
\] %
\textsl{Let }$\left\{ L_{n}\left[ f\right] \right\} $\textsl{\
denote the Lagrange interpolation polynomials to }$f$\textsl{\ at
the zeros of Hermite polynomials. Then for }$p>1,$
\[
\lim_{n\rightarrow \infty }\left\Vert \left( f-L_{n}\left[
f\right] \right) W\right\Vert _{L_{p}(\mathbb{R}) }=0.
\]\end{thh}

Note that the $p$th power of the weight $W$ appears in the
integral. Nevai showed that the presence of such a factor is
necessary --- see his still relevant survey for more perspectives
on this \cite{Nevai1976}. The author and D.~Matjila
\cite{LubinskyMatjila1995} subsequently extended Nevai's result to
general Freud weights:

\begin{thh}
\textsl{Let }$f$\textsl{\ be continuous on the real line. Let }$\beta >1$, $%
p>1,$\textsl{\ and }$W(x)=\exp \left( -\frac{1}{2}\left\vert
x\right\vert ^{\beta }\right) $\textsl{. Let }$\Delta \in \mathbb{R}$\textsl{%
, }$\alpha >0$\textsl{, and }$\hat{\alpha}:=\min \left\{ 1,\alpha \right\} $%
\textsl{. Moreover, let}
\begin{equation}
\tau :=\tau (p) :=\frac{1}{p}-\hat{\alpha}+\frac{\beta
}{6}\max \left\{ 0,1-\frac{4}{p}\right\} .
\end{equation}%
\textsl{Let }$\left\{ L_{n}\left[ f\right] \right\} $\textsl{\
denote the Lagrange interpolation polynomials to }$f$ \textsl{\ at
the zeros of} $\left\{
p_{n}( W^{2},x) \right\}$. \textsl{For }%
\begin{equation}
\lim_{n\rightarrow \infty }\left\Vert \left( f-L_{n}\left[
f\right] \right) (x)W(x)\left( 1+\vert x\vert \right) ^{-\Delta }\right\Vert
_{L_{p}(\mathbb{R}) }^{p}=0
\end{equation}%
\textsl{to hold for every continuous function satisfying }%
\begin{equation}
\lim_{\vert x\vert \rightarrow \infty }f(x)
\left( 1+\vert x\vert \right) ^{\alpha }W(x) =0
\end{equation}%
\textsl{it is necessary and sufficient that}
\[
\begin{array}{cl}
\Delta >\tau & \mbox{if }1<p\leq 4; \\[10pt]
\Delta >\tau & \mbox{if }p>4\mbox{ and }\alpha =1; \\[10pt]
\Delta \geq \tau & \mbox{if }p>4\mbox{ and }\alpha \neq 1.%
\end{array}%
\] %
\textsl{Moreover, for $(11.6)$ to hold for every }$1<p<\infty
$\textsl{\ and every continuous }$f$\textsl{\ satisfying (11.7),
it is necessary and
sufficient that }$\Delta \geq -\hat{\alpha}+\max \left\{ 1,\frac{\beta }{6}%
\right\} $.\end{thh}

This in turn, is a special case of a result for general Freud weights:%

\begin{thh}
\textsl{Let }$W=\exp (-Q) \in \mathcal{F}^{\ast }$. \textsl{Let }%
$p>1,$ $\Delta \in \mathbb{R}$, $\alpha >0$, \textsl{and} $\hat{\alpha}%
:=\min \left\{ 1,\alpha \right\} $. \textsl{Let }$\left\{ L_{n}\left[ f%
\right] \right\} $\textsl{\ denote the Lagrange interpolation
polynomials to }$f$ \textsl{\ at the zeros of} $\left\{
p_{n}( W^{2},x) \right\}$. \textsl{For $(11.6)$ to hold
for every continuous function satisfying $(11.7)$, it is necessary
and sufficient that}
\begin{equation}
\begin{array}{ll}
\Delta >-\hat{\alpha}+\frac{1}{p} & \mbox{if }1<p\leq 4; \\[10pt]
a_{n}^{\frac{1}{p}-\left( \hat{\alpha}+\Delta \right)
}n^{\frac{1}{6}\left( 1-\frac{4}{p}\right) }=O\left( \frac{1}{\log
n}\right) & \mbox{if }p>4\mbox{
and }\alpha =1; \\[10pt]
a_{n}^{\frac{1}{p}-\left( \hat{\alpha}+\Delta \right)
}n^{\frac{1}{6}\left( 1-\frac{4}{p}\right) }=O(1) &
\mbox{if }p>4\mbox{ and }\alpha
\neq 1.%
\end{array}%
\end{equation}\end{thh}

The main ideas of proof of these results already appear in Nevai's
1980 paper, although, obviously, technical details are more
complicated for general weights. In analyzing $L_{n}\left[
f\right] $, one fixes a
polynomial $P$ and splits the function $f-P$ into two pieces, for example, %
\[
L_{n}\left[ f-P\right] =L_{n}\left[ (f-P) \chi _{\left[ -\frac{1%
}{4}a_{n},\frac{1}{4}a_{n}\right] }\right] +L_{n}\left[ (f-P) \chi _{\mathbb{R}\backslash \left[
-\frac{1}{4}a_{n},\frac{1}{4}a_{n}\right] }\right] ,
\] %
and handles them separately, using different techniques. One of
the reasons for the split is the bound
\[
\left\vert p_{n}(x)\right\vert W(x)\leq
Ca_{n}^{-1/2},\qquad \vert x\vert \leq \left( 1-\varepsilon \right) a_{n},
\] %
any fixed $\varepsilon >0$, which does not extend all the way to $\pm a_{n}$;
recall the discussion at (10.5) and (10.6). Indeed, it is the fact that $%
p_{n}(x)W(x)$ behaves roughly like $\left( 1-\frac{%
\vert x\vert }{a_{n}}\right) ^{-1/4}$ near $\pm a_{n}$
that accounts for the change in phenomena at $p=4$, as well as for
the complicated conditions.

The weighted $L_{p}$ norms of the two ``pieces'' of $f$ are
themselves split into a number of pieces, with different
techniques used to estimate over different ranges. Tools include
estimates of quadrature sums in terms of integrals (see Section
12), bounds on Hilbert transforms, and estimates of special
quadrature sums. Here we shall establish just one of the main
estimates, using ideas that go back to Richard Askey
\cite{Askey1972}.

\begin{lem} \textsl{Let}%
\begin{equation}
\psi _{n}=(f-P) \chi _{\left[ -\frac{1}{4}a_{n},\frac{1}{4}a_{n}%
\right] }
\end{equation}%
\textsl{and}%
\begin{equation}
\Lambda =\left\Vert (f-P)(x)W(x)\left( 1+\vert x\vert \right) ^{\alpha }\right\Vert
_{L_{\infty }(\mathbb{R}) }.
\end{equation}%
\textsl{Then for }$n\geq 1,$%
\begin{equation}
\left\Vert L_{n}\left[ \psi _{n}\right] (x)W(x) \left( 1+\vert x\vert \right) ^{-\Delta
}\right\Vert _{L_{p}\left[
-\frac{1}{2}a_{n},\frac{1}{2}a_{n}\right] }\leq C\Lambda ,
\end{equation}%
\textsl{where }$C$\textsl{\ is independent of }$f,n,P.$\end{lem}

\proof.  We let
\[
G(x)=W(x)^{-1}( 1+\left\vert
x\right\vert) ^{-\alpha }
\] %
and first show that
\begin{eqnarray}
&&\hskip-3cm\left\Vert L_{n}\left[ \psi _{n}\right] (x) W(x)\left( 1+\vert x\vert \right)
^{-\Delta }\right\Vert _{L_{p}\left[
-\frac{1}{2}a_{n},\frac{1}{2}a_{n}\right] } \nonumber\\[10pt]
&\leq &C\Lambda \sup_{\left\Vert h\right\Vert _{L_{\infty }(
\mathbb{R})} \leq 1}\left\Vert S_{n}\left[ hG\right] (x)
W(x)\left( 1+\vert x\vert \right) ^{-\Delta }\right\Vert
_{L_{p}\left[ -\frac{1}{2}a_{n},\frac{1}{2}a_{n}\right] }.
\end{eqnarray}%
Then we show that
\begin{equation}
\sup_{\left\Vert h\right\Vert _{L_{\infty }(\mathbb{R})} \leq
1}\left\Vert S_{n}\left[ hG\right] \left(
x\right) W(x)\left( 1+\vert x\vert \right) ^{-\Delta }\right\Vert _{L_{p}\left[ -%
\frac{1}{2}a_{n},\frac{1}{2}a_{n}\right] }\leq C,
\end{equation}%
giving the result.

\proof{(11.12)}. We use duality, but in a very precise
form. Let
\[
g_{n}(x):={\rm sign}\left\{ L_{n}\left[ \psi
_{n}\right] (x)\right\} \left\vert L_{n}\left[ \psi
_{n}\right] (x)\right\vert ^{p-1}W^{p-2}(x) \left( 1+\vert x\vert \right) ^{-\Delta p}\chi
_{\left[ -\frac{1}{2}a_{n},\frac{1}{2}a_{n}\right] }(x) .
\] %
Then%
\begin{eqnarray}
\left\Vert L_{n}\left[ \psi _{n}\right] (x)W(x) \left( 1+\vert x\vert \right) ^{-\Delta
}\right\Vert _{L_{p}\left[
-\frac{1}{2}a_{n},\frac{1}{2}a_{n}\right] }^{p} &=&\int_{-\infty
}^{\infty }L_{n}\left[ \psi _{n}\right] g_{n}(x) W^{2}(x)dx  \nonumber  \nonumber\\
&=&\int_{-\infty }^{\infty }L_{n}\left[ \psi _{n}\right] S_{n}\left[ g_{n}%
\right] (x)W^{2}(x)dx,
\end{eqnarray}%
by orthogonality of $g-S_{n}\left[ g\right] $ to polynomials of
degree $\leq $ $n-1$. Using the Gauss quadrature formula, and the
interpolatory
properties of $\psi _{n}$, we continue this as%
\begin{eqnarray*}
=\sum_{k=1}^{n}\lambda _{kn}\psi _{n}( x_{kn}) S_{n}\left[ g%
\right] ( x_{kn})
\leq \Lambda \sum_{\left\vert x_{kn}\right\vert \leq \frac{1}{4}%
a_{n}}\lambda _{kn}\left\vert S_{n}\left[ g\right] \left( x_{kn}\right) \right\vert W^{-1}( x_{kn}) \left(
1+\left\vert x_{kn}\right\vert \right) ^{-\alpha }.
\end{eqnarray*}%
We estimate the quadrature sum by an integral. These types of
inequalities will be discussed in Section 12. For the moment, we
note that $\lambda _{kn}$ behaves roughly like
$\frac{a_{n}}{n}W^{2}\left( x_{kn}\right) $ and
continue this as%
\begin{eqnarray*}
\leq C\Lambda \int_{-\infty }^{\infty }\left\vert S_{n}\left[
g_{n}\right] (x)\right\vert W(x)\left( 1+\vert x\vert \right) ^{-\alpha }dx &=&C\Lambda
\int_{-\infty }^{\infty }S_{n}\left[ g_{n}\right] (x) h_{n}(x)G(x)W^{2}(x)dx \\
&=&C\Lambda \int_{-\infty }^{\infty }g_{n}(x)S_{n}\left[ h_{n}G%
\right] (x)W^{2}(x)dx,
\end{eqnarray*}%
where $h_{n}(x):={\rm sign}\left( S_{n}\left[
g_{n}\right] (x) \right) $, and we have used the self-adjointness of $S_{n}$. Since $%
g_{n}$ vanishes outside $\left[
-\frac{1}{2}a_{n},\frac{1}{2}a_{n}\right] $,
and letting $q=\frac{p}{p-1}$, we continue this as%
\begin{eqnarray*}
&=&C\Lambda \int_{-\frac{1}{2}a_{n}}^{\frac{1}{2}a_{n}}g_{n}(x)
S_{n}\left[ h_{n}G\right] (x)W^{2}(x)dx \\
&\leq &C\Lambda \left\Vert g_{n}(x)W(x)
\left( 1+\vert x\vert \right) ^{\Delta }\right\Vert _{L_{q}\left[ -\frac{%
1}{2}a_{n},\frac{1}{2}a_{n}\right] }
 \left\Vert S_{n}\left[ h_{n}G\right] (x)W(x) \left( 1+\vert x\vert \right) ^{-\Delta
}\right\Vert
_{L_{p}\left[ -\frac{1}{2}a_{n},\frac{1}{2}a_{n}\right] } \\[10pt]
&=&C\Lambda \left\Vert L_{n}\left[ \psi _{n}\right] (x) W(x)\left( 1+\vert x\vert \right)
^{-\Delta }\right\Vert _{L_{p}\left[
-\frac{1}{2}a_{n},\frac{1}{2}a_{n}\right] }^{p-1}
 \left\Vert S_{n}\left[ h_{n}G\right] (x)W(x) \left( 1+\vert x\vert \right) ^{-\Delta
}\right\Vert _{L_{p}\left[
-\frac{1}{2}a_{n},\frac{1}{2}a_{n}\right] },
\end{eqnarray*}%
by the form of $g_{n}$. Cancelling the $(p-1) $th
power of
\[
\left\Vert L_{n}\left[ \psi _{n}\right] (x)W(x) \left( 1+\vert x\vert \right) ^{-\Delta
}\right\Vert \]
 in the
extreme left of (11.14), we obtain (11.12).

\proof{(11.13)}. We use the Christoffel--Darboux formula
much as we did in the section on mean convergence of orthogonal
expansions:
\[
S_{n}\left[ hG\right] (x)=\frac{\gamma _{n-1}}{\gamma _{n}}%
\left\{ p_{n}(x)H\left[ p_{n-1}hGW^{2}\right] (x) -p_{n-1}(x)H\left[ p_{n}hGW^{2}\right]
(x)\right\} .
\] %
Recalling the bound (10.6) on $p_{n}$ and the bound
\[
\frac{\gamma _{n-1}}{\gamma _{n}}\leq Ca_{n},
\] %
we obtain%
\[
\left\vert S_{n}\left[ hG\right] W\right\vert (x)\leq
Ca_{n}^{1/2}\sum_{j=n-1}^{n}\left\vert H\left[ p_{j}hGW^{2}\right]
(x)\right\vert .
\] %
We next let
\[
\chi _{n}^{\ast }=\chi _{\left[
-\frac{3}{4}a_{n},\frac{3}{4}a_{n}\right] },
\] %
and split%
\begin{eqnarray}
\left\vert S_{n}\left[ hG\right] W\right\vert (x)\leq Ca_{n}^{1/2}
\sum_{j=n-1}^{n}\left\vert H\left[ p_{j}hG\chi _{n}^{\ast
}W^{2}\right] (x)\right\vert \cr
+Ca_{n}^{1/2}\sum_{j=n-1}^{n}\left\vert H\left[ p_{j}hG\left(
1-\chi _{n}^{\ast }\right) W^{2}\right] (x) \right\vert .
\end{eqnarray}%
We next use weighted bounds for Hilbert transforms, of the type
established
by Muckenhoupt \cite{Muckenhoupt1970}. Let%
\[
\hat{\Delta}=\min \left\{ \Delta ,\frac{1}{p}-\delta \right\} ,
\] %
for some small $\delta >0$. Then by Muckenhoupt's bounds for
weighted Hilbert transforms,
\begin{eqnarray}
&&\hskip-1cm a_{n}^{1/2}\sum_{j=n-1}^{n}\left\Vert H\left[ p_{j}hG\chi _{n}^{\ast }W^{2}%
\right] (x)\left( 1+\vert x\vert \right) ^{-\hat{%
\Delta}}\right\Vert _{L_{p}\left[
-\frac{1}{2}a_{n},\frac{1}{2}a_{n}\right] }\nonumber
\\
&\leq &Ca_{n}^{1/2}\sum_{j=n-1}^{n}\left\Vert \left( p_{j}hG\chi
_{n}^{\ast
}W^{2}\right) (x)\left( 1+\vert x\vert \right) ^{-%
\hat{\Delta}}\right\Vert _{L_{p}(\mathbb{R}) }
\leq C\left\Vert (hGW) (x)\left( 1+\left\vert
x\right\vert \right) ^{-\hat{\Delta}}\right\Vert _{L_{p}\left[ -\frac{3}{4}%
a_{n},\frac{3}{4}a_{n}\right] }\nonumber \\
&\leq& C\left\Vert \left( 1+\vert x\vert \right) ^{-\alpha -\hat{%
\Delta}}\right\Vert _{L_{p}\left[
-\frac{3}{4}a_{n},\frac{3}{4}a_{n}\right] }\leq C,
\end{eqnarray}%
as $\left( \alpha +\hat{\Delta}\right) p>1$. (Recall our
hypotheses in
Theorem 11.4 --- they imply $\Delta +\alpha >\frac{1}{p}$.) Next, for $%
\vert x\vert \leq \frac{1}{2}a_{n}$, and $j=n-1,n,$%
\begin{eqnarray*}
\left\vert H\left[ p_{j}hG\left( 1-\chi _{n}^{\ast }\right)
W^{2}\right] (x)\right\vert
&=&\left\vert \int_{\left\{ t:\vert t\vert \geq \frac{3}{4}%
a_{n}\right\} }\frac{\left( p_{j}hGW^{2}\right) (t) }{x-t}%
dt\right\vert
\leq C\int_{\left\{ t:\vert t\vert \geq \frac{3}{4}%
a_{n}\right\} }\left\vert p_{j}W\right\vert (t)
t^{-1-\alpha }dt
\\[10pt]
&\leq &C\left[ \int_{-\infty }^{\infty }\left( p_{j}W\right)
^{2}\right] ^{1/2}\left[ \int_{\frac{3}{4}a_{n}}^{\infty
}t^{-2-2\alpha }dt\right] ^{1/2} \leq Ca_{n}^{-\frac{1}{2}-\alpha
}.
\end{eqnarray*}%
So%
\begin{eqnarray}
&& a_{n}^{1/2}\sum_{j=n-1}^{n}\left\Vert H\left[ p_{j}hG\left( 1-\chi _{n}^{\ast }\right) W^{2}\right] (x)\left(
1+\left\vert
x\right\vert \right) ^{-\Delta }\right\Vert _{L_{p}\left[ -\frac{1}{2}a_{n},%
\frac{1}{2}a_{n}\right] } \leq Ca_{n}^{-\alpha }\left\Vert \left( 1+\vert x\vert \right)
^{-\Delta }\right\Vert _{L_{p}\left[ -\frac{1}{2}a_{n},\frac{1}{2}a_{n}%
\right] } \nonumber\\
&&\hskip3cm\leq Ca_{n}^{-\alpha }\left\{
\begin{array}{ll}
1, & \Delta p>1, \\
\left( \log n\right) ^{\frac{1}{p}}, & \Delta p=1, \\
a_{n}^{\frac{1}{p}-\Delta }, & \Delta p<1.%
\end{array}%
\right. =o(1) ,
\end{eqnarray}%
since%
\[
-\alpha +\frac{1}{p}-\Delta <0.
\] %
Combining (11.15), (11.16), and (11.17), yields (11.13) and hence
the result.  \endproof

Of course, this is just part of the estimation, but is the most
difficult part. For the full details, see
\cite{LubinskyMatjila1995}.

Mean convergence of Lagrange interpolation associated with
Erd\H{o}s weights has been investigated by S.~Damelin and the author \cite%
{DamelinLubinsky1996A}, \cite{DamelinLubinsky1996B}. The author
has
investigated Lagrange interpolation associated with exponential weights on $%
\left[ -1,1\right] $ \cite{Lubinsky1998}, while D.~Kubayi and the
author investigated mean convergence associated with general
exponential weights \cite{KubayiLubinsky2003}. The convergence in
$L_{p}$, $p<1$, has been investigated by Matjila
\cite{Matjila1994B}.

One obvious question is whether other interpolation sets lead to
cleaner, or simpler, results. It was J.~Szabados
\cite{Szabados1997}, \cite{Vertesi2005} who came up with this
idea, the ``method of additional points'', in the context of
Lebesgue functions, which we shall discuss shortly. Szabados' idea
was to add two extra points of interpolation, near $\pm a_{n}$, to
damp
the growth of $p_{n}(x)$ (recall, roughly $\left( 1-\frac{%
\vert x\vert }{a_{n}}\right) ^{-1/4}$) near $\pm
a_{n}$. Suitable
points are those where $p_{n}W$ attains its sup norm on the real line:%
\[
\vert p_{n}W\vert( \xi _{n}) =\left\Vert
p_{n}W\right\Vert _{L_{\infty }(\mathbb{R}) }.
\] %
One can show that
\[
\left\vert 1-\xi _{n}/a_{n}\right\vert \leq Cn^{-2/3}
\] %
and that%
\[
\left\vert x_{1n}-\xi _{n}\right\vert /a_{n}\sim n^{-2/3}.
\] %
We let $L_{n}^{\ast }\left[ f\right] $ denote the polynomial of
degree $\leq
n+1$ that interpolates to $f$ at the $n+2$ points%
\[
\left\{ -\xi _{n},x_{nn},x_{n-1,n},\ldots,x_{2n},x_{1n},\xi
_{n}\right\} .
\] %
Inspired by Szabados' work on Lebesgue functions, G.~Mastroianni
and the author \cite{LubinskyMastroianni1999} proved:

\begin{thh}  \textsl{Let }$W\in \mathcal{F}^{\ast }$.
\textsl{Let }$p>1,$ $\Delta \in \mathbb{R}$, $\alpha >0$,
\textsl{and} $\hat{\alpha}:=\min \left\{ 1,\alpha
\right\} $. \textsl{For }%
\begin{equation}
\lim_{n\rightarrow \infty }\left\Vert \left( f-L_{n}^{\ast }\left[
f\right] \right) (x)W(x)\left( 1+\vert x\vert \right) ^{\Delta }\right\Vert
_{L_{p}(\mathbb{R}) }^{p}=0
\end{equation}%
\textsl{to hold for every continuous function satisfying }%
\begin{equation}
\lim_{\vert x\vert \rightarrow \infty }f(x)
\left( 1+\vert x\vert \right) ^{\alpha }W(x) =0
\end{equation}%
\textsl{it is necessary and sufficient that}%
\[
\Delta >\frac{1}{p}-\hat{\alpha}.
\]\end{thh}

What a difference the extra two points make! Let us explain a
little of the mechanics behind this. $L_{n}^{\ast }\left[ f\right]
$ admits the representation
\begin{eqnarray*}
L_{n}^{\ast }\left[ f\right] (x)=f\left( \xi
_{n}\right) \ell _{0n}^{\ast }(x)+f\left( -\xi
_{n}\right) \ell _{n+1,n}^{\ast }(x)
+\sum_{j=1}^{n}f\left( x_{jn}\right) \ell _{j,n}^{\ast }(x) ,
\end{eqnarray*}%
where for $1\leq j\leq n,$%
\[
\ell _{j,n}^{\ast }(x)=\ell _{jn}(x)
\frac{\xi _{n}^{2}-x^{2}}{\xi _{n}^{2}-x_{jn}^{2}},
\] %
and
\[
\ell _{jn}(x)=\frac{p_{n}(x)
}{p_{n}^{\prime }\left( x_{jn}\right) \left( x-x_{jn}\right) }.
\] %
The factor $\xi _{n}^{2}-x^{2}$ damps the growth of $p_{n}(x) $ near $\pm a_{n}$. Of course, there is an extra factor
of $\xi
_{n}^{2}-x_{jn}^{2}$ in the denominator, but this is large for most $x_{jn}$%
, so actually helps.

We already noted that convergence of $\left\{ L_{n}\left[ f\right]
\right\} $ associated with Freud weights in $L_{p}$, $p<1$, has
been investigated by Matjila \cite{Matjila1994B}. If one
sacrifices precision in this case, one can achieve very general
results, using distribution functions,
rearrangements of functions, and a classic lemma of Loomis \cite%
{BennettSharpley1988}. The latter asserts that for any $\left\{
c_{j}\right\} _{j=1}^{n}$,
\[
{\rm meas}\left\{ x:\left\vert
\sum_{j=1}^{n}\frac{c_{j}}{x-x_{j}}\right\vert
>\lambda \right\} \leq \frac{8}{\lambda }\sum_{j=1}^{n}\left\vert
c_{j}\right\vert .
\] %
(It is often stated in the case where all $c_{j}>0$ and without
absolute
values.) If one writes%
\[
L_{n}\left[ f\right] (x)=p_{n}(x)\sum_{j=1}^{n}%
\frac{f\left( x_{jn}\right) }{p_{n}^{\prime }\left( x_{jn}\right)
\left( x-x_{jn}\right) }=:p_{n}(x)g_{n}(x) ,
\] %
we see that for any positive function $\phi $ and real numbers $b,c,$%
\begin{equation}
\Vert L_{n}\left[ f\right] W\phi ^{b}\Vert _{L_{p}( \mathbb{R%
}) }\leq \Vert p_{n}W\phi ^{b+c}\Vert
_{L_{2p}(\mathbb{R}) }\Vert g_{n}\phi
^{-c}\Vert _{L_{2p}(\mathbb{R}) }.
\end{equation}%
Next, Loomis' Lemma shows that%
\begin{eqnarray*}
m_{g_{n}}( \lambda ) :={\rm meas}\left\{ x:\left\vert
g_{n}(x)\right\vert >\lambda \right\}
&\leq &\frac{8\left\Vert fW\phi ^{c}\right\Vert _{L_{\infty }( \mathbb{R%
}) }}{\lambda }\sum_{j=1}^{n}\frac{1}{\vert
p_{n}^{\prime }W\phi
^{c}\vert( x_{jn}) } \\
&=:&\frac{8\left\Vert fW\phi ^{c}\right\Vert _{L_{\infty }\left( \mathbb{R}%
\right) }}{\lambda }\Omega _{n}.
\end{eqnarray*}%
For the moment, let $h^{\ast }$ denote the decreasing
rearrangement of a function $h$. In particular, the decreasing
rearrangement of $g_{n}$ is
\begin{eqnarray*}
g_{n}^{\ast }(t) =\sup \left\{ \lambda
:m_{g_{n}}\left( \lambda \right) >t\right\} \leq \sup \left\{
\lambda :\frac{8\left\Vert fW\phi ^{c}\right\Vert _{L_{\infty
}(\mathbb{R}) }}{\lambda }\Omega _{n}>t\right\}
=\frac{8\left\Vert fW\phi ^{c}\right\Vert _{L_{\infty }\left( \mathbb{R}%
\right) }}{t}\Omega _{n}.
\end{eqnarray*}%
Inequalities for decreasing rearrangements give \cite{BennettSharpley1988}%
\begin{eqnarray}
\left\Vert g_{n}\phi ^{-c}\right\Vert _{L_{2p}(\mathbb{R}) }^{2p} &=&\int_{-\infty }^{\infty }\left\vert
g_{n}\phi ^{-c}\right\vert ^{2p}  \nonumber \leq \int_{0}^{\infty
}\left( \left\vert g_{n}\right\vert ^{2p}\right) ^{\ast }\left( \phi ^{-2pc}\right) ^{\ast }  \nonumber =\int_{0}^{\infty }\left(
g_{n}^{\ast }\right) ^{2p}\left( \left( \phi
^{-1}\right) ^{\ast }\right) ^{2pc}  \nonumber \\[10pt]
&\leq &\left( 8\left\Vert fW\phi ^{c}\right\Vert _{L_{\infty }\left( \mathbb{%
R}\right) }\Omega _{n}\right) ^{2p}\int_{0}^{\infty }t^{-2p}\left(( \phi ^{-1}) ^{\ast }(t) \right) ^{2pc}dt.
\end{eqnarray}%
Now if $p<\frac{1}{2}$, this integrand poses no problem at $0$. If
for example,
\[
\phi (t) =1+\vert t\vert \mbox{,}
\] %
then%
\[
( \phi ^{-1}) ^{\ast }(t) =\left( 1+\frac{t}{2}%
\right) ^{-1},\qquad t\geq 0,
\] %
and the integral will converge provided $2p+2pc>1$. Combining
this, (11.20)
and (11.21), we obtain%
\[
\Vert L_{n}\left[ f\right] W\phi ^{b}\Vert _{L_{p}( \mathbb{R%
}) }\leq C\Vert fW\phi ^{c}\Vert _{L_{\infty
}(\mathbb{R}) }\Omega _{n}\Vert p_{n}W\phi
^{b+c}\Vert _{L_{2p}(\mathbb{R}) }
\] %
with $C$ independent of $n,f$. Recalling the definition of $\Omega
_{n}$
above, we have the main part of \cite[Theorem 3, p.~155]{Lubinsky2002}:%

\begin{thh}  \textsl{Let }$W$\textsl{\ be even and
positive in }$[0,\infty )$\textsl{\ and }$\phi (t)
=1+\vert t\vert $. \textsl{The following are
equivalent:}
\begin{description}
\item[(I)] {There exist }$b,c\in \mathbb{R}$\textsl{\
and }$C,p>0$\textsl{\ such
that for every continuous }$f$\textsl{\ }%
\[
\Vert L_{n}\left[ f\right] W\phi ^{b}\Vert _{L_{p}( \mathbb{R%
}) }\leq C\Vert fW\phi ^{c}\Vert _{L_{\infty
}(\mathbb{R}) }.
\] %
\item[(II)]\textsl{\  There exist }$\beta ,\gamma \in \mathbb{R}$\textsl{\ and }$%
r>0 $\textsl{\ such that}
\[
\sup_{n\geq 1}\Vert p_{n}W\phi ^{\beta }\Vert
_{L_{r}(\mathbb{R}) }\sum_{j=1}^{n}\frac{1}{\vert
p_{n}^{\prime }W\phi ^{\gamma }\vert( x_{jn})
}<\infty .
\]\end{description}\end{thh}

From this one may easily obtain necessary and sufficient
conditions for mean
convergence of $\left\{ L_{n}\left[ f\right] \right\} $ in $L_{p}$, some $p<%
\frac{1}{2}$, with appropriate weights. The striking feature of
these results is the simplicity of the proofs, and their
generality. However, one sacrifices a great deal of precision, for
even the specific $L_{p}$ is not fixed in advance. For further
work on mean convergence of Lagrange
interpolation associated with exponential weights, see \cite{Damelinetal2001}%
, \cite{Damelinetal2001C}, \cite{Damelinetal2003}, \cite{Jung2003}, \cite%
{MastroianniOccorsio2001}, \cite{MastroianniOccorsio2002}, \cite%
{MastroianniVertesi2006}, \cite{OccorsioRusso2006}, \cite{Sakai1998}, \cite%
{Sakai1998B}.

\subsection{Lebesgue Functions and Pointwise Convergence}

\bigskip

Since $L_{n}\left[ P\right] =P$ for any polynomial $P$ of degree
$\leq n-1$, the convergence of $\left\{ L_{n}\right\} $ in almost
any norm is equivalent to the norms of these operators being
bounded independent of $n$. The \dword{Lebesgue function} is the
weighted norm of $L_{n}\left[ f\right] $ at a given point. When
working on infinite intervals, it makes more sense to investigate
weighted version of these. The idea is
\begin{eqnarray*}
\left\vert L_{n}\left[ f\right] (x)\right\vert
W(x)&\leq &\sum_{j=1}^{n}\left\vert f\left( x_{jn}\right) \right\vert \left\vert
\ell _{jn}(x)\right\vert W(x)\\
&\leq &\left\Vert fW\right\Vert _{L_{\infty }(\mathbb{R}) }W(x)\sum_{j=1}^{n}\left\vert \ell
_{jn}(x)\right\vert W^{-1}\left( x_{jn}\right)
=:\left\Vert fW\right\Vert _{L_{\infty }(\mathbb{R})
}\Lambda _{n}(x).
\end{eqnarray*}%
It is easy to see by suitable choice of $f$ that%
\[
\Lambda _{n}(x)=\sup \left\{ \left\vert L_{n}\left[
f\right] (x)W(x)\right\vert :\left\Vert
fW\right\Vert _{L_{\infty }(\mathbb{R}) }\leq
1\right\} ,
\] %
so the Lebesgue function $\Lambda _{n}(x)$ is indeed
a norm. The \dword{Lebesgue constant} of $L_{n}$ is
\[
\Lambda _{n}:=\left\Vert \Lambda _{n}\right\Vert _{L_{\infty }\left( \mathbb{%
R}\right) }.
\] %
Observe that for any function $f$, there is the pointwise error estimate%
\begin{eqnarray*}
W(x)\left\vert f-L_{n}\left[ f\right] \right\vert
(x)&\leq &\left( 1+\Lambda _{n}(x)
\right) \inf_{\deg (P) \leq n-1}\left\Vert
(f-P)W\right\Vert _{L_{\infty }(\mathbb{R}) } \\[10pt]
&=&\left( 1+\Lambda _{n}(x)\right) E_{n}\left[
f;W\right] _{\infty }
\end{eqnarray*}%
and hence the global estimate%
\[
\left\Vert W\left( f-L_{n}\left[ f\right] \right) \right\Vert
_{L_{\infty
}(\mathbb{R}) }\leq \left( 1+\Lambda _{n}\right) E_{n}\left[ f;W%
\right] _{\infty }.
\]

Amongst the contributors to estimation of $\Lambda _{n}(x) $ on infinite intervals are Damelin, Freud, Horv\'{a}th,
Kubayi, Matjila, Nevai, Sklyarov, Szabados, Szab\'{o}, Szili and
V\'{e}rtesi. The earliest results were due to Freud and Sklyarov.

The following result for general Freud weights is due to Szabados \cite%
{Szabados1997}. Indeed, it was in this context that Szabados
introduced his ``method of additional points". Matjila
\cite{Matjila1994}, \cite{Matjila1995} independently obtained the
upper bound implicit in (a).

\begin{thh}
\textsl{Let }$W=\exp (-Q) \in \mathcal{F}^{\ast }$.%
\begin{description}
\item[(a)] \textsl{Let }$\Lambda _{n}$\textsl{\ denote the Lebesgue constant for
the Lagrange interpolation polynomials }$L_{n}\left[ \cdot \right]
$\textsl{\ at
the zeros of }$p_{n}\left( W^{2},\cdot \right) $\textsl{. Then }%
\[
\Lambda _{n}\sim n^{1/6}\mbox{, }\qquad n\geq 1.
\] %
\item[(b)]\textsl{ Let }$\Lambda _{n}$\textsl{\ denote the Lebesgue
constant for
the Lagrange interpolation polynomials }$L_{n}^{\ast }\left[ \cdot \right] $%
\textsl{\ at the zeros of }$p_{n}\left( W^{2},\cdot \right)
$\textsl{, together with the points }$\pm \xi _{n}$\textsl{\ where
}$\left\vert
p_{n}W\right\vert $\textsl{\ attains its maximum. Then}%
\[
\Lambda _{n}^{\ast }\sim \log (n+1) \mbox{, }\qquad
n\geq 1.
\]\end{description}\end{thh}

The $n^{1/6}$ arises from the factor $\left( 1-\frac{\left\vert
x\right\vert }{a_{n}}\right) ^{-1/4}$ that gives the growth of
$p_{n}$ near $\pm a_{n}$, together with the fact that
\[
1-\frac{x_{1n}}{a_{n}}\sim n^{-2/3}.
\] %
It is remarkable that Szabados' addition of the two points $\pm
\xi _{n}$
gives the optimal factor of $\log n$. Yes, that is the same lower bound of $%
\log n$ for the sup norm of projection operators onto the space of
polynomials of degree $\leq n$, in the context of finite intervals \cite%
{SzabadosVertesi1990}. In the context of exponential weights, this
optimality was proved by Szabados \cite{Szabados1997}. There are
extensions of this result to Erd\H{o}s weights, and exponential
weights on $\left( -1,1\right) $ \cite{Damelin1998B},
\cite{Damelin1998C}, \cite{Kubayi2001}, \cite{Kubayi2002B},
\cite{Kubayi2002C}. Typically there one obtains the same result
for $\Lambda _{n}^{\ast }$, while
\[
\Lambda _{n}\sim \left( nT( a_{n}) \right)
^{1/6}\mbox{, }\qquad n\geq 1.
\]

Some inkling of the proofs comes from the representation for the
fundamental
polynomials%
\[
\ell _{jn}(x)=\lambda _{jn}K_{n}\left( x,x_{jn}\right) ,
\] %
when we interpolate at the zeros of $p_{n}$. Then%
\begin{eqnarray*}
\Lambda _{n}(x)&=&W(x)
\sum_{j=1}^{n}\left\vert \ell _{jn}(x)\right\vert
W^{-1}( x_{jn}) =W(x)\sum_{j=1}^{n}\lambda
_{jn}W^{-1}( x_{jn})
\left\vert K_{n}( x,x_{jn}) \right\vert \\
&\approx &W(x)\int_{-\infty }^{\infty }W(t) \left\vert K_{n}( x,t) \right\vert dt.
\end{eqnarray*}%
Recall that $\lambda _{jn}\sim \frac{a_{n}}{n}W^{2}( x_{jn}) $ at least for $x_{jn}$ well inside the interval
$\left[ -a_{n},a_{n}\right] $.

For general sets of interpolation nodes, lower bounds for Lebesgue
functions
have been obtained by V\'{e}rtesi and his coworkers \cite{Vertesi1998}, \cite%
{Vertesi1999}. These show that most of the time the Lebesgue
function is bounded below by $\log n$
\cite[p.~359]{Vertesi1999}:

\begin{thh}  \textsl{Let }$W\in \mathcal{F}^{\ast
}$\textsl{. For }$n\geq 1$\textsl{, let }$\Lambda _{n}(x) $\textsl{\ denote the Lebesgue function corresponding to
}$n$\textsl{\ distinct points (not necessarily the zeros of
}$p_{n}$\textsl{). Let }$\varepsilon >0$\textsl{. Then there is a set }$%
\mathcal{H}_{n}$\textsl{\ of linear Lebesgue measure }$\leq
2\varepsilon
a_{n}$\textsl{, such that }%
\[
\Lambda _{n}(x)\geq C\varepsilon \log n\mbox{,
}\qquad x\in \left[ -a_{n},a_{n}\right] \backslash
\mathcal{H}_{n}.
\] %
\textsl{The constant is independent of }$n,x,\varepsilon
$\textsl{, and the array of interpolation points.}\end{thh}

Szabados analyzed the influence of the distribution of the
interpolation points, and in particular the largest interpolation
point, on the size of the fundamental polynomials and Lebesgue
constant \cite{Szabados2000}. Kubayi \cite{Kubayi2002} and Damelin
\cite{Damelin2002} explored Szabados' result in the context of
more general weights. Summability of weighted
Lagrange interpolation has been investigated by Szili and V\'{e}rtesi \cite%
{SziliVertesi2006}, \cite{SziliVertesi2006B},
\cite{SziliVertesi2006C}.

When one increases the degree of the polynomial interpolating at
$n$ points from $n-1$ to $\floor{ n(1+\varepsilon )}$, some
fixed $\varepsilon >0$, then one can avoid the factor of $\log n$
above. For finite intervals, this was an old result of Erd\H{o}s.
For Freud weights, it is due to V\'{e}rtesi \cite{Vertesi2001},
and for exponential weights on $\left[ -1,1\right] $, due to Szili
and V\'{e}rtesi \cite{SziliVertesi2003}.

\begin{thh}  \textsl{Let }$W\in \mathcal{F}^{\ast
}$\textsl{. Let us be given an array of interpolation points
}$\left\{ x_{jn}\right\} _{n\geq 1,1\leq j\leq n}$ \textsl{with
associated fundamental polynomials of interpolation }$\left\{ \ell
_{jn}\right\} _{n\geq 1,1\leq j\leq n}$. \textsl{Assume that these
admit the bound}%
\[
\sup_{x\in \mathbb{R}}\left\vert \ell _{jn}(x)
\right\vert W^{-1}\left( x_{jn}\right) W(x)\leq A,
\] %
\textsl{for }$n\geq 1$\textsl{\ and }$1\leq j\leq n$\textsl{. Let }$%
\varepsilon >0$\textsl{\ and let }$f$\textsl{\ be a continuous
function such that }$fW^{1+\varepsilon }$\textsl{\ has limit
}$0$\textsl{\ at }$\pm \infty $\textsl{. Then there exist
polynomials }$\left\{ P_{n}\right\} $\textsl{\ such that}
\smallskip

\begin{description}
\item[(I)] $P_{n}$\textsl{\ has degree }$N_{n}\leq n\left( 1+\varepsilon +C\varepsilon n^{-2/3}\right) ;$
\item[(II)] $P_{n}\left( x_{jn}\right) =f\left( x_{jn}\right)
$\textsl{, }$1\leq j\leq n;$
\item[(III)] $\Vert W^{1+\varepsilon }( f-P_{n})
\Vert _{L_{\infty }(\mathbb{R}) }\leq
CE_{N_{n}}\left[ f;W^{1+\varepsilon }\right] _{\infty }.$
\end{description}\end{thh}

Note that the fundamental polynomials associated with the
orthogonal polynomials $p_n(W^2,x)$ satisfy the above bound
\cite{LevinLubinsky1992}.

Despite the many investigations of Lagrange interpolation, there
are relatively few results where local assumptions are imposed
that guarantee pointwise convergence of $\left\{ L_{n}\left[
f\right] \right\} $. It seems that the following interesting
result of Nevai has never been extended beyond the case of the
Hermite weight \cite{Nevai1974}:

\begin{thh}
\textsl{Let }$W(x)=\exp( -x^{2}) $\textsl{. Let }$f$%
\textsl{\ be Riemann integrable in each finite interval, and
assume that for
some }$c<\frac{1}{2},$%
\[
\left\Vert fW^{c}\right\Vert _{L_{\infty }(\mathbb{R})
}<\infty .
\] %
\textsl{Assume that on }$\left[ a,b\right] $\textsl{,
}$f$\textsl{\
satisfies the one-sided Dini condition }%
\[
f(x+t) -f(x)\geq -v(t)
\left\vert \log t\right\vert ^{-1}\mbox{, }\qquad a<x<x+t<b,
\] %
\textsl{where }$v(t) $\textsl{\ decreases to }$0$\textsl{\ as }$%
t $\textsl{\ decreases to }$0$\textsl{. Then}%
\begin{eqnarray}
\lim_{n\rightarrow \infty }L_{n}\left[ f\right] (x) =f(x)
\end{eqnarray} %
\textsl{at each point of continuity of }$f$\textsl{\ in }$\left( a,b\right) $%
\textsl{. The convergence is uniform in every closed subinterval
of }$\left( a,b\right) $\textsl{\ where }$f$\textsl{\ is
continuous.}\end{thh}

Note that any increasing function $f$ satisfies this Dini
condition, so (11.22) holds also when $f$ is of bounded variation
in $\left[ a,b\right] $. (For then $f$ is the difference of
monotone increasing functions.) Finally, we note that convergence
of Lagrange interpolation to entire functions of suitably
restricted growth leads to a geometric rate of convergence, see
\cite{Mhaskar1996}.

\sect{Marcinkiewicz--Zygmund Inequalities}

The classical Marcinkiewicz--Zygmund inequality has the form%
\[
C_{1}\leq \left( \frac{1}{2\pi }\int_{-\pi }^{\pi }\left\vert
R( \theta
) \right\vert ^{p}d\theta \right) ^{1/p}\Bigl/\left( \frac{1}{2n+1}%
\sum_{j=0}^{2n}\left\vert R\!\left( \frac{j\pi }{2n+1}\right)
\right\vert ^{p}\right) \leq C_{2},
\] %
where $n\geq 1$, $R$ is a trigonometric polynomial of degree $\leq n$, and $%
C_{1}$ and $C_{2}$ are independent of $n$ and $R$ \cite[Vol.~2, p.~28]%
{Zygmund2002}. These types of inequalities are useful, for
example, in discretisation problems, quadrature theory, and
Lagrange interpolation.

In the context of exponential weights, these have chiefly been
used for mean convergence of Lagrange interpolation. In fact, the
proof of the Erd\H{o}s--Tur\'{a}n theorem above was based on a
converse quadrature sum, namely,
\[
\int_{-\infty }^{\infty }L_{n}\left[ f\right] ^{2}W^{2}=\sum_{j=1}^{n}%
\lambda _{jn}f^{2}( x_{jn}) .
\] %
This of course is an immediate consequence of the Gauss quadrature
formula. In our outline of the proof of Lemma 11.5 above we used a
forward quadrature sum. Let us repeat this, in simplified form:
assume that we have an
inequality%
\begin{equation}
\sum_{j=1}^{n}\lambda _{jn}W^{-2}( x_{jn}) \vert
PW\vert( x_{jn}) \leq C\int_{-\infty }^{\infty
}\left\vert PW\right\vert .
\end{equation}%
If $q=\frac{p}{p-1}$, duality gives%
\[
\left\Vert L_{n}\left[ f\right] W\right\Vert _{L_{p}(\mathbb{R}) }=\sup_{g}\int_{-\infty }^{\infty }L_{n}\left[
f\right] gW^{2},
\] %
where the $\sup $ is taken over all $g$ with $\left\Vert
gW\right\Vert
_{L_{q}(\mathbb{R}) }\leq 1$. Using orthogonality of $g-S_{n}%
\left[ g\right] $ to polynomials of degree $\leq n-1$, and then
the Gauss
quadrature formula gives%
\begin{eqnarray*}
&=&\sup_{g}\int_{-\infty }^{\infty }L_{n}\left[ f\right]
S_{n}\left[ g\right] W^{2}
=\sup_{g}\sum_{j=1}^{n}\lambda _{jn}f( x_{jn}) S_{n}\left[ g%
\right] ( x_{jn}) \\
&\leq &\left\Vert fW\right\Vert _{L_{\infty }(\mathbb{R}) }\sum_{j=1}^{n}\lambda _{jn}W^{-2}(
x_{jn}) \left\vert S_{n}\left[
g\right]( x_{jn}) W( x_{jn}) \right\vert \\
&\leq &C\left\Vert fW\right\Vert _{L_{\infty }(\mathbb{R}) }\int_{-\infty }^{\infty }\left\vert S_{n}\left[
g\right] \right\vert W.
\end{eqnarray*}%
In the last line, we used the quadrature sum estimate (12.1). Setting $%
\sigma _{n}:={\rm sign}\left( S_{n}\left[ g\right] \right) $, and
then using
self-adjointness of $S_{n}$, we continue the above as%
\begin{eqnarray*}
=C\left\Vert fW\right\Vert _{L_{\infty }(\mathbb{R})
}\int_{-\infty }^{\infty }\sigma _{n}W^{-1}S_{n}\left[ g\right]
W^{2} &=&C\left\Vert fW\right\Vert _{L_{\infty }(\mathbb{R})
}\int_{-\infty }^{\infty }S_{n}\left[ \sigma_n W^{-1}\right] gW^{2} \\
&\leq &C\left\Vert fW\right\Vert _{L_{\infty }(\mathbb{R})
}\left\Vert S_{n}\left[\sigma_n W^{-1}\right] W\right\Vert
_{L_{p}(\mathbb{R}) }.
\end{eqnarray*}%
If $S_{n}$ is bounded in a suitable sense, we obtain%
\[
\left\Vert L_{n}\left[ f\right] W\right\Vert _{L_{p}(\mathbb{R}) }\leq C\left\Vert fW\right\Vert _{L_{\infty
}(\mathbb{R}) }
\] %
with $C$ independent of $n$. We emphasize that this is a
simplified version of what one needs. Typically there are extra
factors inside the quadrature sum. In subsequent sections, we
shall present several methods to establish such quadrature sums.

\subsection{Nevai's Method for Forward Quadrature Sums}

We let $u\in \left[ x_{jn},x_{j-1,n}\right] $, $p\geq 1,$ and
start with the
following consequence of the fundamental theorem of calculus:%
\[
\vert PW\vert ^{p}(x_{jn}) \leq\vert
PW\vert ^{p}(u) +p\int_{u}^{x_{j-1,n}}\vert
PW\vert ^{p-1}(s) \left\vert (PW)
^{\prime }(s) \right\vert ds.
\] %
We may assume that $u$ is the point in $\left[
x_{jn},x_{j-1,n}\right] $ where $\left\vert PW\right\vert ^{p}$
attains its minimum. We now need estimates for Christoffel
functions, and spacing of zeros of orthogonal polynomials. In fact
the former imply the latter via Markov--Stieltjes
inequalities. We use%
\[
\lambda _{jn}W^{-2}(x_{jn}) =\lambda _{n}( W^{2},x_{jn}) W^{-2}(x_{jn}) \sim \varphi
_{n}(x_{jn}) ,
\] %
where $\varphi _{n}$ is defined by (7.8), and%
\[
x_{jn}-x_{j-1,n}\geq C\lambda _{jn}W^{-2}(x_{jn}) .
\] %
Then%
\begin{eqnarray*}
\lambda _{jn}W^{-2}(x_{jn})\vert PW\vert
^{p}(x_{jn}) &\leq
&C\int_{x_{jn}}^{x_{j-1,n}}\vert PW\vert
^{p}(u) du \\
&&+C\varphi _{n}(x_{jn})
\int_{x_{jn}}^{x_{j-1,n}}\vert PW\vert ^{p-1}(s) \vert (PW) ^{\prime }(s)
\vert ds.
\end{eqnarray*}%
Summing over $j$, and using the fact that $\varphi _{n}$ does not
change much in $\left[ x_{jn},x_{j-1,n}\right] $, one obtains
\begin{eqnarray}
\sum_{j=1}^{n}\lambda _{jn}W^{-2}(x_{jn})\vert PW\vert
^{p}(x_{jn}) \leq C\int_{-\infty }^{\infty }\vert PW\vert ^{p}(u)
du\cr +C\int_{-\infty }^{\infty }\vert PW\vert ^{p-1}(s)\vert (PW)
^{\prime }(s) \vert \varphi _{n}(s) ds.
\end{eqnarray}%
We apply H\"{o}lder's inequality, and then the Bernstein
inequality Theorem
7.6 to the second term:%
\begin{eqnarray*}
\int_{-\infty }^{\infty }\vert PW\vert ^{p-1}(s)\vert (PW) ^{\prime }(s)
\vert \varphi _{n}(s) ds \leq \left\Vert
PW\right\Vert _{L_{p}(\mathbb{R}) }^{p-1}\Vert
(PW)^{\prime }\varphi _{n}\Vert _{L_{p}(\mathbb{R}) } \leq C\left\Vert PW\right\Vert _{L_{p}(
\mathbb{R}) }^{p-1}\left\Vert PW\right\Vert _{L_{p}(\mathbb{R}) }.
\end{eqnarray*}%
Together with (12.2), this gives%
\begin{equation}
\sum_{j=1}^{n}\lambda _{jn}W^{-2}(x_{jn})\vert
PW\vert ^{p}(x_{jn}) \leq C\left\Vert
PW\right\Vert _{L_{p}(\mathbb{R}) }^{p}.
\end{equation}

Here is a sample of what can be proved using this method \cite[Thm.~2, p.~287%
]{LubinskyMatjila1995B}:

\begin{thh}  \textsl{Let }$W\in \mathcal{F}^{\ast }$,
$1<p<\infty $, $-\infty <b\leq 2$\textsl{\ and }$a\in
\mathbb{R}$\textsl{. Then}
\begin{eqnarray*}
\sum_{j=1}^{n}\lambda _{jn}W^{-b}(x_{jn}) \left( 1+\left\vert x_{jn}\right\vert \right) ^{a}\vert
PW\vert ^{p}(x_{jn}) \leq C\int_{-\infty
}^{\infty }\vert PW\vert ^{p}(t) \left( 1+\vert t\vert \right) ^{a}W^{2-b}(t) dt.
\end{eqnarray*}\end{thh}

This method has also been used by H. K\"{o}nig in a vector valued
and Banach space setting, with $W$ as the Hermite weight \cite{Konig1994}, \cite%
{KonigNielsen1994}. It was used first for Hermite weights by Nevai \cite%
{Nevai1980}. Of course, it can be used for points other than zeros
of orthogonal polynomials --- all we need are suitable estimates
on the spacing between successive zeros.

\subsection{The Large Sieve Method}

As far as the author is aware, this also was developed by
P.~Nevai, based on the large sieve of number theory, and used by
him and his coworkers for exponential weights
\cite{Lubinskyetal1987}. We already presented the idea in the
context of Markov--Bernstein inequalities in Section 8.6. We start
with Lemma 8.3. For polynomials $S$ of degree $\leq n,$
\[
\vert SW\vert ^{p}(\xi) \leq
C\frac{\int_{-\infty }^{\infty }\vert SW\vert
^{p}(t) K_{Mn}( \xi ,t) ^{2}W^{2}(t) dt}{\int_{-\infty }^{\infty }K_{Mn}( \xi ,t)
^{2}W^{2}(t) dt},
\] %
which is equivalent to
\[
\vert SW\vert ^{p}(\xi) \leq C\lambda
_{Mn}(\xi) \int_{-\infty }^{\infty }\vert
SW\vert ^{p}(t) K_{Mn}( \xi ,t)
^{2}W^{2}(t) dt.
\] %
Here $p>0$, $S$ is a polynomial of degree $\leq n$, and
$\left\vert \xi \right\vert \leq 2a_{n}$, while $M$ is a fixed
large enough positive integer, independent of $n,P,\xi $. We now
simply take linear combinations of this, and use the bound
\[
\lambda _{Mn}( W^{2},\xi) \leq
C\frac{a_{n}}{n}W^{2}(\xi) ,\mbox{ }\qquad\left\vert
\xi \right\vert \leq 2a_{n}.
\] %
For example, if we use the zeros of $p_{n},$ we obtain%
\[
\sum_{j=1}^{n}\vert SW\vert ^{p}(x_{jn})
\leq C\int_{-\infty }^{\infty }\vert SW\vert ^{p}(t) \Sigma _{n}(t) dt,
\] %
where%
\[
\Sigma _{n}(t) =\frac{a_{n}}{n}W^{2}(t)
 \sum_{j=1}^{n}K_{Mn}^{2}\left( x_{jn},t\right)W^{2}(x_{jn}).
\] %
Here, by the Christoffel--Darboux formula, and our bounds for
$p_{n}$, and recalling our bound (8.31) for $K_{Mn}\left( \xi ,\xi
\right) $, we see that for $\left\vert \xi \right\vert ,\left\vert
t\right\vert \leq 2a_{n},$
\[
W(\xi) W(t) \left\vert K_{Mn}( \xi
,t) \right\vert \leq \min \left\{
\frac{n}{a_{n}},\frac{1}{\left\vert \xi -t\right\vert }\right\} .
\] %
Moreover for Freud weights, uniformly in $j,n$%
\[
x_{jn}-x_{j+1,n}\geq C\frac{a_{n}}{n}.
\] %
Hence%
\begin{eqnarray*}
\Sigma _{n}(t) \leq C\sum_{j=1}^{n}\left( x_{jn}-x_{j+1,n}\right)
\min \left\{ \frac{a_{n}}{n},\left\vert x_{jn}-t\right\vert
\right\} ^{-2} \leq C\int_{-\infty }^{\infty }\min \left\{
\frac{a_{n}}{n},\left\vert u-t\right\vert \right\} ^{-2}du\leq
C_{1}\frac{n}{a_{n}}.
\end{eqnarray*}%
Thus we arrive at the estimate
\[
\sum_{j=1}^{n}\frac{a_{n}}{n}\left\vert SW\right\vert ^{p}(x_{jn}) \leq C\int_{-\infty }^{\infty }\vert
SW\vert ^{p}(t) dt,
\] %
valid for $n\geq 1$ and polynomials $P$ of degree $\leq n$. We
note that this does not imply the estimate
\[
\sum_{j=1}^{n}\lambda _{jn}W^{-2}(x_{jn})\vert
SW\vert ^{p}(x_{jn}) \leq C\int_{-\infty
}^{\infty }\vert SW\vert ^{p}(t) dt
\] %
since
\[
\lambda _{jn}W^{-2}(x_{jn}) \sim \frac{a_{n}}{n}\left( \left\vert
1-\frac{\left\vert x_{jn}\right\vert }{a_{n}}\right\vert
+n^{-2/3}\right) ^{-1/2} \gg \frac{a_{n}}{n},
\] %
when $x_{jn}$ is close to $\pm a_{n}$. Thus this method is in some
respects not as powerful as Nevai's method. However, it can easily
yield estimates like
\[
\sum_{j=1}^{n}\frac{a_{n}}{n}\psi(\vert SW\vert
^{p}(x_{jn})) \leq C\int_{-\infty }^{\infty
}\psi(\vert SW\vert ^{p}(t))
dt
\] %
where $\psi $ is an increasing convex function that is positive in
$\left( 0,\infty \right) $ and has $\psi (0) =0$, and
works for any set of suitably spaced points.

\subsection{The Duality Method}

The idea here is to use duality of $L_{p}$ spaces to derive a
forward
quadrature sum estimate from a converse one. It was probably first used by K%
\"{o}nig \cite{Konig1994}, \cite{KonigNielsen1994}. Let us
illustrate this in the context of quadrature sums at zeros of
orthogonal polynomials. Let $n$ be fixed and $\mu _{n}$ be a pure
jump measure having mass $\lambda
_{jn}W^{-2}(x_{jn}) $ at $x_{jn}$. Let, as usual, $q=\frac{p}{p-1%
}$. If $P$ is a polynomial of degree $\leq n$,
\begin{eqnarray*}
\left( \sum_{j=1}^{n}\lambda _{jn}W^{-2}(x_{jn})
\vert PW\vert ^{p}(x_{jn}) \right) ^{1/p}
=\left( \int\vert PW\vert ^{p}d\mu _{n}\right) ^{1/p}
=\sup_{g}\int PW^{2}g\mbox{ }d\mu _{n},
\end{eqnarray*}%
where the sup is taken over all $g$ with
\[
\int \left\vert gW\right\vert ^{q}d\mu _{n}\leq 1.
\] %
Since $g$ needs to be defined only at $\left\{ x_{jn}\right\}
_{j=1}^{n}$,
we can assume that $g$ is a polynomial of degree $\leq n-1$. So%
\begin{eqnarray*}
\int PW^{2}g\mbox{ }d\mu _{n} =\sum_{j=1}^{n}\lambda _{jn}(Pg) (x_{jn}) =\int_{-\infty }^{\infty }PgW^{2}
\leq \left\Vert PW\right\Vert _{L_{p}(\mathbb{R})
}\left\Vert gW\right\Vert _{L_{q}(\mathbb{R}) }.
\end{eqnarray*}%
Now we assume a suitable converse quadrature sum inequality: for
$n\geq 1$ and all $S$ of degree $\leq n,$
\begin{equation}
\left\Vert SW\right\Vert _{L_{q}(\mathbb{R}) }\leq
C\left( \sum_{j=1}^{n}\lambda _{jn}W^{-2}(x_{jn})
\vert SW\vert ^{q}(x_{jn}) \right) ^{1/q}.
\end{equation}%
Applying this to $g$, gives%
\begin{eqnarray*}
\left\Vert gW\right\Vert _{L_{q}(\mathbb{R}) } \leq
C\left( \sum_{j=1}^{n}\lambda _{jn}W^{-2}(x_{jn})
\vert gW\vert ^{q}(x_{jn}) \right) ^{1/q}
=C\left( \int\vert gW\vert ^{q}d\mu _{n}\right)
^{1/q}\leq C.
\end{eqnarray*}%
Thus we obtain
\[
\left( \sum_{j=1}^{n}\lambda _{jn}W^{-2}(x_{jn})
\vert SW\vert ^{p}(x_{jn}) \right)
^{1/p}\leq C\left\Vert PW\right\Vert _{L_{p}(\mathbb{R}) }.
\] %
In K\"{o}nig's work, he first proved a converse inequality like (12.4) for $%
1<q<4$, and then deduced a forward estimate for
$\frac{4}{3}<p<\infty $. Note that one can also use Carleson
measures to derive forward quadrature estimates, just as we
derived Markov--Bernstein inequalities in Section
8.5. For this approach on finite intervals, see the survey \cite%
{Lubinsky1998B}.

\subsection{The Duality Method for Converse Quadrature Sum estimates}

This was perhaps the earliest method used for converse quadrature
sum estimates. For trigonometric polynomials, it appears in
Zygmund's treatise
\cite[Ch. X, pp.~28--29]{Zygmund2002}. Let $1<p<\infty $ and $q=\frac{p}{p-1}$%
. Let $P$ be a polynomial of degree $\leq n-1$. By duality,
\[
\left\Vert PW\right\Vert _{L_{p}(\mathbb{R})
}=\sup_{g}\int gPW^{2},
\] %
where the sup is taken over all $g$ with $\left\Vert gW\right\Vert
_{L_{q}(\mathbb{R}) }\leq 1$. Let $S_{n}\left[
g\right] $ denote
the $n$th partial sum of the orthonormal expansion of $g$ with respect to $%
\left\{ p_{j}\right\} _{j=0}^{\infty }$. By orthogonality, we
continue this as
\[
\left\Vert PW\right\Vert _{L_{p}(\mathbb{R}) }=\sup_{g}\int S_{n}%
\left[ g\right] PW^{2}.
\] %
Using the Gauss quadrature formula, and then H\"{o}lder's
inequality, we
continue this as%
\begin{eqnarray*}
&=&\sup_{g}\sum_{j=1}^{n}\lambda _{jn}S_{n}\left[ g\right] (x_{jn}) P(x_{jn}) \\
&\leq &\sup_{g}\left\{ \sum_{j=1}^{n}\lambda _{jn}W^{-2}(x_{jn}) \vert S_{n}\left[ g\right] W\vert
^{q}(x_{jn}) \right\} ^{1/q} \left\{
\sum_{j=1}^{n}\lambda _{jn}W^{-2}(x_{jn})
\vert PW\vert ^{p}(x_{jn}) \right\} ^{1/p} \\[10pt]
&=:&\sup_{g}T_{1}\times T_{2}.
\end{eqnarray*}%
If we have a suitable forward quadrature sum estimate, and mean
boundedness
of $S_{n}$ in a weighted $L_{q}$ norm, we can bound $T_{1}$ as follows:%
\[
T_{1}\leq C\left\Vert S_{n}\left[ g\right] W\right\Vert
_{L_{q}(\mathbb{R}) }\leq C\left\Vert gW\right\Vert _{L_{q}\left( \mathbb{R}%
\right) }\leq C.
\] %
Then for $n\geq 1$ and all polynomials $P$ of degree $\leq n$,%
\[
\left\Vert PW\right\Vert _{L_{p}(\mathbb{R}) }\leq
C\left\{ \sum_{j=1}^{n}\lambda _{jn}W^{-2}(x_{jn})
\vert PW\vert ^{p}(x_{jn}) \right\} ^{1/p}.
\] %
This method does yield estimates for full Gauss quadrature sums.
However, it requires mean boundedness of orthogonal expansions,
such as that in Theorem 10.3, which is difficult to prove. It also
requires a forward quadrature sum estimate, which is not all that
difficult to prove. It works best for $p\geq 4$. For $p<4$, the
next method gives better results. Here is a sample of what one can
achieve \cite[Thm.~1.2, p.~531]{Lubinsky1995B}. We note that
because we are considering a special case of results there, the
hypotheses there are simplified below:

\begin{thh}  \textsl{Let }$W\in \mathcal{F}^{\ast
}$\textsl{, and assume that the
orthonormal polynomials }$\left\{ p_{n}\right\} $\textsl{\ for }$W^{2}$%
\textsl{\ satisfy the bound $(10.7)$. Let }$p>4$\textsl{\ and }$r,R\in \mathbb{%
R}$ \textsl{satisfy}%
\begin{equation}
R>-\frac{1}{p}
\end{equation}%
\textsl{and}
\begin{equation}
a_{n}^{r-\min \left\{ R,1-\frac{1}{p}\right\} }n^{\frac{1}{6}\left( 1-\frac{4%
}{p}\right) }=\left\{
\begin{array}{ll}
O(1) , & R\neq 1-\frac{1}{p}, \\
O\left( \left( \log n\right) ^{-R}\right) , & R=1-\frac{1}{p}%
\end{array}%
\right. .
\end{equation}%
\textsl{Then for }$n\geq 1$\textsl{\ and polynomials }$P$\textsl{\
of degree
}$\leq n-1,$%
\begin{equation}
\left\Vert (PW) (x)\left( 1+\left\vert
x\right\vert \right) ^{r}\right\Vert _{L_{p}(\mathbb{R}) }\leq C\left\{ \sum_{j=1}^{n}\lambda
_{jn}W^{-2}(x_{jn})\vert PW\vert
^{p}(x_{jn}) \left( 1+\left\vert x_{jn}\right\vert
\right) ^{Rp}\right\} ^{1/p}.
\end{equation}%
\textsl{For }$p=4$\textsl{, this persists provided $(12.5)$ holds
and}
\[
r-\min \left\{ R,1-\frac{1}{p}\right\} <0.
\]\end{thh}

Weights that satisfy the bound (10.7) include exp$\left( -\vert x\vert ^{\alpha }\right) $, $\alpha >1$. In
\cite{Lubinsky1995B}, we showed that slightly weaker forms of
these conditions on $r,R$ are necessary.

\subsection{K\"{o}nig's Method}

K\"{o}nig's method \cite{Konig1994}, \cite{KonigNielsen1994} is
based on a clever estimate for Hilbert transforms of
characteristic functions of intervals. It is technically the most
difficult of all the methods we present (for forward or converse
estimates), but is very powerful, and relatively direct. It does
not depend on mean boundedness of orthogonal
expansions. For $n\geq 1$ and $1\leq j\leq n,$ let%
\begin{eqnarray*}
I_{jn} :=[x_{jn},x_{j-1,n}); \qquad \left\vert I_{jn}\right\vert
&=&x_{j-1,n}-x_{j,n};
 \qquad
\chi _{jn} :=\chi _{I_{jn}.}
\end{eqnarray*}%
Let $P$ be a polynomial of degree $\leq n-1$, and
\[
y_{jn}:=a_{n}^{-1/2}\frac{P(x_{jn}) }{p_{n}^{\prime
}(x_{jn}) }.
\] %
Then%
\begin{eqnarray*}
P(x)&=&L_{n}\left[ P\right] (x)
=a_{n}^{1/2}p_{n}(x)\sum_{j=1}^{n}\frac{y_{jn}}{x-x_{jn}} \\
&=&a_{n}^{1/2}p_{n}(x)\sum_{j=1}^{n}y_{jn}\left\{ \frac{1}{%
x-x_{jn}}-\frac{1}{\left\vert I_{jn}\right\vert }H\left[ \chi
_{jn}\right] (x)\right\}
+a_{n}^{1/2}p_{n}(x)H\left[ \sum_{j=1}^{n}\frac{y_{jn}}{%
\left\vert I_{jn}\right\vert }\chi _{jn}\right] (x)\\
&=:&J_{1}(x)+J_{2}(x).
\end{eqnarray*}%
Here (as above),
\[
H\left[ f\right] (x)=\lim_{\varepsilon \rightarrow
0+}\int_{\left\vert x-t\right\vert \geq \varepsilon }\frac{f(t) }{t-x}dt
\] %
is the Hilbert transform of $f$. The term $J_{2}$ is easier, so
let us deal with it first. We use the bound
\[
\vert p_{n}W\vert (x)\leq Ca_{n}^{-1/2}\left\vert 1-%
\frac{\vert x\vert }{a_{n}}\right\vert ^{-1/4},\qquad
x\in \mathbb{R},
\] %
valid when $W\in \mathcal{F}^{\ast }$, and the Hilbert transform bound \cite%
{Muckenhoupt1970}%
\[
\Vert H\left[ g\right] (x)\left\vert
1-\frac{\vert x\vert }{a_{n}}\right\vert
^{-1/4}\Vert _{L_{p}( \vert x\vert \leq
2a_{n}) }\leq C\Vert g(x)\left\vert
1-\frac{\vert x\vert }{a_{n}}\right\vert
^{-1/4}\Vert _{L_{p}( \vert x\vert \leq
2a_{n}) },
\] %
valid for $1<p<4$. Here $C\neq C\left( n,g\right) $. Then%
\begin{eqnarray*}
\left\Vert J_{2}W\right\Vert _{L_{p}\left[ -2a_{n},2a_{n}\right] }
&=&a_{n}^{1/2}\Vert p_{n}WH\left[ \sum_{j=1}^{n}\frac{y_{jn}}{%
\left\vert I_{jn}\right\vert }\chi _{jn}\right]\Vert
_{L_{p}\left[
-2a_{n},2a_{n}\right] } \\
&\leq &C\Vert \left\vert 1-\frac{\vert x\vert }{a_{n}}%
\right\vert ^{-1/4}\sum_{j=1}^{n}\frac{y_{jn}}{\left\vert I_{jn}\right\vert }%
\chi _{jn}(x)\Vert _{L_{p}\left[
-2a_{n},2a_{n}\right] }
\\[10pt]
&=&C\left\{ \sum_{j=1}^{n}\left( \frac{y_{jn}}{\left\vert I_{jn}\right\vert }%
\right) ^{p}\int_{I_{jn}}\left\vert 1-\frac{\vert x\vert }{a_{n}}%
\right\vert ^{-p/4}dx\right\} ^{1/p} \\[10pt]
&\leq &C\left\{ \sum_{j=1}^{n}\left( \frac{y_{jn}}{\left\vert
I_{jn}\right\vert }\right) ^{p}\left\vert I_{jn}\right\vert \psi
_{n}^{-p/4}(x_{jn}) \right\} ^{1/p}.
\end{eqnarray*}%
Here%
\[
\psi _{n}(x)=\max \left\{ n^{-2/3},1-\frac{\left\vert
x\right\vert }{a_{n}}\right\} .
\] %
It can be shown that uniformly in $j,n$,%
\[
\frac{y_{jn}}{\left\vert I_{jn}\right\vert }\sim \left\vert
PW\right\vert (x_{jn}) \psi _{n}^{1/4}(x_{jn}) ;\quad a_{n}^{1/2}\left\vert p_{n}^{\prime
}W\right\vert (x_{jn}) \sim \left\vert
I_{jn}\right\vert ^{-1}\psi _{n}^{-1/4}(x_{jn}) ;
\quad \left\vert I_{jn}\right\vert \sim \lambda _{jn}W^{-2}(x_{jn}) .
\] %
Hence%
\[
\left\Vert J_{2}W\right\Vert _{L_{p}\left[ -2a_{n},2a_{n}\right]
}\leq C\left\{ \sum_{j=1}^{n}\lambda _{jn}W^{-2}(x_{jn})\vert PW\vert ^{p}(x_{jn})
\right\} ^{1/p}.
\] %
The estimation of $J_{1}$ is difficult. One shows that for
$\left\vert x-x_{jn}\right\vert >2\left\vert I_{jn}\right\vert ,$
\begin{eqnarray*}
\left\vert \frac{1}{x-x_{jn}}-\frac{1}{\left\vert I_{jn}\right\vert }H%
\left[ \chi _{jn}\right] (x)\right\vert
a_{n}^{1/2}\vert p_{n}W\vert (x)
&\leq& C\psi _{n}^{-1/4}(x)\frac{\left\vert I_{jn}\right\vert }{%
\left\vert x-x_{jn}\right\vert }\left\{ \frac{1}{\left\vert
x-x_{jn}\right\vert }+\frac{1}{1+\left\vert x_{jn}\right\vert }\right\}\\
&=:&Cf_{jn}(x),
\end{eqnarray*}%
and for $\left\vert x-x_{jn}\right\vert \leq 2\left\vert I_{jn}\right\vert ,$%
\begin{eqnarray*}
\left\vert \frac{1}{x-x_{jn}}-\frac{1}{\left\vert I_{jn}\right\vert }H%
\left[ \chi _{jn}\right] (x)\right\vert
a_{n}^{1/2}\vert p_{n}W\vert (x)\leq
C\frac{\psi _{n}^{-1/4}(x)}{\left\vert
I_{jn}\right\vert }=:Cf_{jn}(x),
\end{eqnarray*}%
so that
\[
\left\vert J_{1}(x)\right\vert \leq
C\sum_{j=1}^{n}\left\vert y_{jn}\right\vert f_{jn}(x)
.
\] %
As each $f_{jn}$ does not change much in $I_{kn}$, one obtains%
\begin{eqnarray*}
\left\Vert J_{1}W\right\Vert _{L_{p}\left[ x_{nn},x_{1n}\right] }
\leq C\left\{ \sum_{k=2}^{n}\left\vert I_{kn}\right\vert \left[
\sum_{j=1}^{n}\left\vert y_{jn}\right\vert f_{jn}(x_{kn}) \right] ^{p}\right\} ^{1/p}.
\end{eqnarray*}%
One now splits this sum into three pieces. The main part is
estimated by showing that certain $n$ by $n$ matrices are
uniformly bounded in norm. These in turn depend on clever
consequences of H\"{o}lder's inequality.
After a lot of technical work, one obtains%
\[
\left\Vert J_{1}W\right\Vert _{L_{p}(\mathbb{R}) }\leq
C\left\{ \sum_{j=1}^{n}\lambda _{jn}W^{-2}(x_{jn})
\vert PW\vert ^{p}(x_{jn}) \right\} ^{1/p}.
\] %
Here is a sample of what can be proved \cite[Thm.~1.2, p.~531]{Lubinsky1995B}:

\begin{thh}
\textsl{Let }$W\in \mathcal{F}^{\ast }$\textsl{\ and }$1<p<4$\textsl{. Let }%
\begin{equation}
r<1-\frac{1}{p};\qquad r\leq R;\qquad R>-\frac{1}{p}.
\end{equation}%
\textsl{Then for }$n\geq 1$\textsl{\ and polynomials }$P$\textsl{\
of degree
}$\leq n$\textsl{,}%
\begin{eqnarray}
\left\Vert (PW) (x)\left( 1+\left\vert
x\right\vert \right) ^{r}\right\Vert _{L_{p}(\mathbb{R}) } \leq C\left\{ \sum_{j=1}^{n}\lambda
_{jn}W^{-2}(x_{jn}) \vert PW\vert
^{p}(x_{jn})( 1+\left\vert x_{jn}\right\vert
) ^{Rp}\right\} ^{1/p}.
\end{eqnarray}%
Here $C\neq C\left( n,P\right) $.\end{thh}

There it was also shown that the
first two conditions on $r,R$ are necessary. For $p\geq 4$, the
author and Mastroianni
used K\"{o}nig's method to prove \cite[Thm.~1.2, p.~149]%
{LubinskyMastroianni2002}.

\begin{thh}  \textsl{Let }$W\in \mathcal{F}^{\ast
}$\textsl{\ and }$p\geq 4$\textsl{. Let
(12.8) hold, and in addition, assume that for some }$\delta >0,$\textsl{\ }%
\begin{equation}
n^{\frac{1}{6}\left( 1-\frac{4}{p}\right) }a_{n}^{r-\min \left\{ R,1-\frac{1%
}{p}\right\} }=O( n^{-\delta }) .
\end{equation}%
\textsl{Then $(12.9)$ holds for }$n\geq 1$\textsl{\ and polynomials }$P$%
\textsl{\ of degree }$\leq n$\textsl{.}\end{thh}

For further work on quadrature sum estimates in the context of
exponential
weights, see \cite{Damelin2003}, \cite{Damelinetal2002}, \cite{Lubinsky1999}%
. In particular, in the latter paper, an attempt was made to
provide a very general framework for K\"{o}nig's method.

\sect{Hermite and Hermite--F\'{e}jer Interpolation}

The Hermite interpolation polynomial to a function $f$ at the
zeros $\left\{
x_{jn}\right\} _{j=1}^{n}$ of $p_{n}(x)$ is a polynomial $H_{n}%
\left[ f\right] $ of degree $\leq 2n-1$ such that for $1\leq j\leq
n,$
\begin{eqnarray*}
H_{n}\left[ f\right] (x_{jn}) &=&f(x_{jn})
\mbox{, }
\\
H_{n}^{\prime }\left[ f\right] (x_{jn}) &=&f^{\prime
}(x_{jn}) .
\end{eqnarray*}%
Of course, we have to assume $f^{\prime }$ is defined at the zeros
$\left\{ x_{jn}\right\} $. To avoid this restriction, the
Hermite--F\'{e}jer interpolant $H_{n}^{\ast }\left[ f\right] $ is
often used instead. It is a polynomial of degree $\leq 2n-1$ such
that for $1\leq j\leq n,$
\begin{eqnarray*}
H_{n}^{\ast }\left[ f\right] (x_{jn}) &=&f(x_{jn})
\mbox{, } \\
H_{n}^{\ast \prime }\left[ f\right] (x_{jn}) &=&0.
\end{eqnarray*}%
Hermite--F\'{e}jer polynomials first came to prominence when
F\'{e}jer used them (with interpolation at zeros of Chebyshev
polynomials) to provide another proof of Weierstrass'
approximation theorem. A close relative of the
Hermite--F\'{e}jer interpolation polynomial is the Gr\"unwald operator%
\[
Y_{n}\left[ f\right] (x)=\sum_{j=1}^{n}f(x_{jn}) \ell _{jn}^{2}(x),
\] %
where $\left\{ \ell _{jn}\right\} $ are the fundamental
polynomials of Lagrange interpolation taken at the zeros of
$p_{n}(x)=p_{n}\left( W^{2},x\right) $. $Y_{n}$ has
the advantage of being a positive operator, which makes
convergence more general, though possibly at a slower rate than
interpolatory operators.

Convergence of Hermite and Hermite--F\'{e}jer polynomials
associated with
exponential weights has been investigated in \cite{Damelinetal2001B}, \cite%
{Damelinetal2001D}, \cite{Jung2005}, \cite{JungKwon1998}, \cite%
{KanjinSakai1994}, \cite{KanjinSakai1995}, \cite{KasugaSakai1999}, \cite%
{KasugaSakai2004}, \cite{KasugaSakai2004B}, \cite{KasugaSakai2005}, \cite%
{Lubinsky1992}, \cite{LubinskyRabinowitz1992}, \cite{LubinskyRabinowitz1992B}%
, \cite{MhaskarXu1992}, \cite{Szabo1999}. In some cases, processes
involving higher order derivatives have been studied, and the
derivatives of the interpolating polynomials have also been
investigated. Here we note just one result of Szab\'{o}
\cite{Szabo1999}:

\begin{thh}  \textsl{Let }$W=\exp (-Q) \in
\mathcal{F}^{\ast }$\textsl{, with }$A$\textsl{\ as in $(4.18)$.
Let}
\[
W_{2}(x)=W^{2}(x)\left( 1+\left\vert
Q^{\prime }(x)\right\vert \right) ^{\frac{1}{3\left( 1-1/A\right) }}\log \left( 2+\left\vert Q^{\prime }(x)
\right\vert \right) .
\] %
\textsl{Then for every continuous} $f:\mathbb{R\rightarrow R}$
\textsl{with}
\[
\lim_{\vert x\vert \rightarrow \infty }fW_{2}(x) =0,
\] %
\textsl{we have}
\[
\lim_{n\rightarrow \infty }\left\Vert W^{2}\left( f-H_{n}^{\ast
}\right) \right\Vert _{L_{\infty }(R) }=0.
\] \end{thh}

Note the square of the weight in the condition on $f$ and in the
approximating norm. Szab\'{o} showed that the factor
\[\left( 1+\left\vert
Q^{\prime }(x)\right\vert \right) ^{\frac{1}{3\left( 1-1/A\right) }}\log \left( 2+\left\vert Q^{\prime }(x)
\right\vert \right) \]
 is essential and best possible.

\sect{Acknowledgement}

The author is grateful for corrections and comments from Steven
Damelin, Zeev Ditzian, Hrushikesh Mhaskar, Ryozi Sakai, Ying Guang
Shi, Misha Sodin, Jozsef Szabados, Vili Totik, and Peter Vertesi.
Undoubtedly many shortcomings remain, for which the author takes
full responsibility. The research in this paper was supported by
NSF Grant DMS--0400446 and Israel US--BSF Grant 2004353.

\sect{Summary of Notation}

In this section, we repeat the main notation introduced earlier in
this paper, and repeat some definitions. Throughout, $C, C_{1},
C_{2}, \ldots$ denote positive constants independent of $n,x,t$ and
polynomials $P$ of degree $\leq n$. The same symbol does not
necessarily denote the same
constant in different occurrences. We write $C=C( \alpha ) $ or $%
C\neq C( \alpha ) $ to respectively show that $C$ depends on $%
\alpha $, or does not depend on $\alpha $.

We use the notation $\sim $ in the following sense: given
sequences of real numbers $\left\{ c_{n}\right\} $ and $\left\{
d_{n}\right\} $, we write
\[
c_{n}\sim d_{n}
\] %
if for some positive constants $C_{1},C_{2}$ independent of $n$,
we have
\[
C_{1}\leq c_{n}/d_{n}\leq C_{2}.
\] %
If $x$ is a real number, then $\floor{x}$ denotes the
greatest integer $\leq x$.

Throughout, $W=\exp (-Q) $ is an even exponential
weight.
Associated with $W$ are the Mhaskar--Rakhmanov--Saff numbers $%
a_{n}=a_{n}(Q) $, the positive root of the equation%
\[
n=\frac{2}{\pi }\int_{0}^{1}\frac{a_{n}tQ^{\prime }(a_{n}t) }{%
\sqrt{1-t^{2}}}dt.
\] %
If $xQ^{\prime }(x)$ is strictly increasing, with
limits $0$ and $\infty $ at $0$ and $\infty $ respectively, then
$a_{n}$ is uniquely
defined. Corresponding to the weight $W$, we define orthonormal polynomials%
\[
p_{n}(x)=p_{n}\left( W^{2},x\right) =\gamma
_{n}x^{n}+\cdots ,\qquad\gamma _{n}>0,
\] %
satisfying
\[
\int_{-\infty }^{\infty }p_{n}p_{m}W^{2}=\delta _{mn.}
\] %
Note that the weight is $W^{2}$, not $W$. We denote the zeros of
$p_{n}$ by
\[
x_{nn}<x_{n-1,n}< \cdots <x_{2n}<x_{1n}.
\]

The $n$th Christoffel function for $W^{2}$ is%
\begin{equation}
\lambda _{n}( W^{2},x) =\inf_{\deg (P) \leq n-1}\frac{%
\int_{-\infty }^{\infty }(PW) ^{2}}{P^{2}(x) }.
\end{equation}%
It also satisfies%
\[
\lambda _{n}( W^{2},x)
=1\Bigl/\;\sum_{j=0}^{n-1}p_{j}^{2}(x).
\] %
The reproducing kernel of order $m$ is
\[
K_{m}( x,t) :=\sum_{j=0}^{m-1}p_{j}(x)
p_{j}(t) =\frac{\gamma _{m-1}}{\gamma
_{m}}\frac{p_{m}(x)p_{m-1}(t)
-p_{m-1}(x)p_{m}(t) }{x-t}.
\]

Our four main classes of weights follow:

\begin{deff} {\bf (Freud Weights $\mathcal{F}$):}
\textsl{Let }$W=\exp (-Q) $\textsl{, where }$Q:\mathbb{R}%
\rightarrow \mathbb{R}$\textsl{\ is even, }$Q^{\prime }$\textsl{\
exists and
is positive in }$\left( 0,\infty \right) $\textsl{. Moreover, assume that }$%
xQ^{\prime }(x)$\textsl{\ is strictly increasing,
with right
limit }$0$\textsl{\ at }$0$\textsl{, and for some }$\lambda ,A,B>1,C>0,$%
\[
A\leq \frac{Q^{\prime }( \lambda x) }{Q^{\prime }(x)}%
\leq B,\qquad x\geq C.
\] %
\textsl{Then we write }$W\in \mathcal{F}$.\end{deff}

\begin{deff} {\bf (Freud Weights $\mathcal{F}^{\ast }$):}
\textsl{Let }$W=\exp (-Q) $\textsl{, where }$Q^{\prime \prime }$%
\textsl{\ exists and is positive in }$\left( 0,\infty \right)
$\textsl{, while }$Q^{\prime }$\textsl{\ is positive there, with
limit }$0$\textsl{\ at
}$0$\textsl{, and for some }$A,B>1,$%
\[
A-1\leq \frac{xQ^{\prime \prime }(x)}{Q^{\prime
}(x)}\leq B-1,\qquad x\in \left( 0,\infty \right) .
\] %
\textsl{Then we write }$W\in \mathcal{F}^{\ast }$.\end{deff}

The canonical examples of $W$ in $\mathcal{F}^{\ast }$ are
\[
W_{\alpha }(x)=\exp( -\vert x\vert
^{\alpha }) ,\qquad \alpha >0,
\] %
while $\mathcal{F}^{\ast }\subset \mathcal{F}.$

\begin{deff} {\bf (Erd\H{o}s Weights $\mathcal{E}$):}
\textsl{Let }$W=\exp (-Q) $\textsl{, where }$Q:\mathbb{R}%
\rightarrow \mathbb{R}$\textsl{\ is even, }$Q^{\prime }$\textsl{\
exists and is positive in }$\left( 0,\infty \right) $\textsl{.
Assume that }$xQ^{\prime
}(x)$\textsl{\ is strictly increasing, with right limit }$0$%
\textsl{\ at }$0$\textsl{, and the function}%
\begin{equation}
T(x):=\frac{xQ^{\prime }(x)}{Q(x) }
\end{equation}%
\textsl{is quasi-increasing in the sense that}
\[
0\leq x<y\quad \Longrightarrow \quad T(x)\leq
CT(y) ,
\] %
\textsl{while}
\begin{equation}
\lim_{x\rightarrow \infty }T(x)=\infty .
\end{equation}%
\textsl{Assume, moreover, that for some }$C_{1},C_{2}$\textsl{\ and }$%
C_{3}>0, $%
\[
\frac{yQ^{\prime }(y) }{xQ^{\prime }(x)
}\leq C_{1}\!\!\left( \frac{Q(y) }{Q(x)
}\right) ^{C_{2}},\qquad y\geq x\geq C_{3}.
\] %
\textsl{Then we write }$W\in \mathcal{E}$.\end{deff}

\begin{deff} {\bf (General Exponential Weights
$\mathcal{F}_{even}( C^{2}) $):}
\textsl{Let }$I=\left( -d,d\right) $\textsl{\ where }$0<d\leq \infty $%
\textsl{. Let }$Q:I\rightarrow \lbrack 0,\infty )$\textsl{\
satisfy the following properties:}
\begin{description}
\item[(a)] $Q^{\prime }$\textsl{\ is continuous and positive in }$I$\textsl{\ and }%
$Q(0) =0;$
\item[(b)] $Q^{\prime \prime
}$\textsl{\ exists and is positive in }$\left( 0,d\right)
;$

\item[(c)]
\[
\lim_{t\rightarrow d-}Q(t) =\infty ;
\] %
\item[(d)] { The function }%
\[
T(t) =\frac{tQ^{\prime }(t) }{Q(t) }
\] %
\textsl{is quasi-increasing in }$\left( 0,d\right) $\textsl{, in
the sense
that }%
\[
0\leq x<y<d\quad \Longrightarrow \quad T(x)\leq
CT(y) .
\] %
\textsl{Moreover, assume for some }$\Lambda >1$\textsl{, }%
\[
T(t) \geq \Lambda >1,\qquad t\in \left( 0,d\right) .
\] %
\item[(e)] {There exists }$C_{1}>0$\textsl{\ such that }%
\[
\frac{\left\vert Q^{\prime \prime }(x)\right\vert
}{Q^{\prime }(x)}\leq C_{1}\frac{Q^{\prime }(x) }{Q(x)},\qquad x\in \left( 0,d\right)
\] %
\end{description}
\textsl{Then we write }$W\in \mathcal{F}_{even}( C^{2}) $\textsl{%
. If also, there exists }$c\in \left( 0,d\right) $\textsl{\ such
that}
\[
\frac{\left\vert Q^{\prime \prime }(x)\right\vert
}{Q^{\prime }(x)}\leq C_{2}\frac{Q^{\prime }(x) }{Q(x)},\qquad x\in \left( c,d\right) ,
\] %
\textsl{then we write }$W\in \mathcal{F}_{even}( C^{2}+) $.\end{deff}

For $k\geq 1,$ the $k$th iterated exponential is
\[
\exp _{k}=\exp( \exp( \cdots\exp ({})
))
\]

Define the symmetric differences%
\begin{eqnarray*}
\Delta _{h}f(x)=f( x+\frac{h}{2}) -f( x-\frac{%
h}{2}) ; \qquad \Delta _{h}^{k}f(x)=\Delta
_{h}\left( \Delta _{h}^{k-1}f(x)\right) ,\quad k\geq
1,
\end{eqnarray*}%
so that
\[
\Delta _{h}^{k}f(x)=\sum_{i=0}^{k}{k\choose i}(-1) ^{i}f( x+k\frac{h}{2}-ih) .
\] %
Given $r\geq 1$ and a Freud weight $W=\exp (-Q) $ with
Mhaskar--Rakhmanov--Saff numbers $\left\{ a_{n}\right\} $, we
define the $r$th order modulus of continuity as follows: First
define the decreasing function
of $t,$%
\[
\sigma (t) :=\inf \left\{ a_{n}:\frac{a_{n}}{n}\leq
t\right\},\qquad t>0.
\] %
The $r$th order modulus of continuity for $W$ is%
\begin{eqnarray*}
\omega _{r,p}(f,W,t) =\sup_{0<h\leq t}\norm{W\left( \Delta
_{h}^{r}f\right)}_{L_{p}[-\sigma (h) ,\sigma (h) ]} +\inf_{\deg (P) \leq r-1}
\norm{(
f-P) W}_{L_{p}\left( \mathbb{R}\backslash \lbrack
-\sigma (t) ,\sigma (t) ]\right) }.
\end{eqnarray*}%
The $r$th order $K$-functional associated with $W$ is%
\[
K_{r,p}\left( f,W,t^{r}\right) :=\inf_{g}\left\{ \left\Vert (f-g) W\right\Vert _{L_{p}(\mathbb{R})
}+t^{r}\Vert g^{(r) }W\Vert _{L_{p}(\mathbb{R}) }\right\} ,\qquad t\geq 0.
\]

For Erd\H{o}s weights $W\in \mathcal{E}$, the modulus of
continuity involves
the function%
\[
\Phi _{t}(x):=\sqrt{1-\frac{\vert x\vert
}{\sigma (t) }}+T( \sigma (t) )
^{-1/2},
\] %
where $\sigma (t) $ is as above. The $r$th order
modulus is then
\begin{eqnarray*}
\omega _{r,p}(f,W,t) &=&\sup_{0<h\leq t}\norm{W(x)
\left( \Delta _{h\Phi _{t}(x)}^{r}f(x)
\right)}
_{L_{p}[-\sigma (2t) ,\sigma (2t) ]} \\[10pt]
&&\hskip1cm +\inf_{\deg (P) \leq r-1} \|(f-P) W\| _{L_{p}\left( \mathbb{R}\backslash \lbrack -\sigma
(4t) ,\sigma (4t) ]\right) }.
\end{eqnarray*}

The error in weighted approximation of $f$ by polynomials of degree $\leq $ $%
n$ in the $L_{p}$ norm is%
\[
E_{n}\left[ f;W\right] _{p}:=\inf_{\deg (P) \leq n}
\norm{(f-P) W}_{L_{p}( \mathbb{R}) }.
\] %
If $f$ is such that $fW^{2}\in L_{1}(\mathbb{R}) $, we
define
its orthonormal expansion%
\[
f\sim \sum_{j=0}^{\infty }c_{j}p_{j},
\] %
where%
\[
c_{j}:=\int_{-\infty }^{\infty }fp_{j}W^{2},\qquad j\geq 0.
\] %
The $n$th partial sum of the orthonormal expansion is%
\[
S_{n}\left[ f\right] =\sum_{j=0}^{n-1}c_{j}p_{j},
\] %
while the $n$th de la Vall\'ee Poussin operator is
\[
V_{n}\left[ f\right] =\frac{1}{n}\sum_{j=n+1}^{2n}S_{j}\left[
f\right] .
\]

Given a function $f$, we define the $n$th Lagrange interpolation
polynomial to $f$ at $\left\{ x_{jn}\right\} _{j=1}^{n}$ to be the
unique polynomial of
degree $\leq n-1$ that satisfies%
\[
L_{n}\left[ f\right] (x_{jn}) =f(x_{jn}) \mbox{, }\qquad%
1\leq j\leq n\mbox{.}
\] %
One formula for $L_{n}\left[ f\right] $ is
\[
L_{n}\left[ f\right] =\sum_{j=1}^{n}f(x_{jn}) \ell
_{jn}
\] %
where $\left\{ \ell _{jn}\right\} _{j=1}^{n}$ are fundamental
polynomials of Lagrange interpolation, satisfying
\[
\ell _{jn}(x_{kn}) =\delta _{jk}.
\] %
The $\left\{ \ell _{jn}\right\} $ admit the representation%
\[
\ell _{jn}(x)=\prod\limits_{k=1,k\neq j}^{n}\frac{x-x_{kn}}{%
x_{jn}-x_{kn}}.
\]

\noindent School of Mathematics\\
Georgia Institute of Technology\\
Atlanta\\
GA 30332-0160\\
lubinsky@math.gatech.edu

\endddoc